%% file: main-r1.tex
\documentclass[conf]{new-aiaa}
\usepackage[utf8]{inputenc}
\usepackage{textcomp}

\usepackage{graphicx}
\usepackage{amsmath}
\usepackage[version=4]{mhchem}
\usepackage{siunitx}
\usepackage{longtable,tabularx}
\usepackage{afterpage}
\usepackage{booktabs}
\usepackage{changepage}
\usepackage[capitalize]{cleveref}
\usepackage[useregional]{datetime2}
\usepackage[acronyms]{glossaries-extra}
\usepackage{mathtools}
\usepackage{multicol}
\usepackage{multirow}
\usepackage{orcidlink}
\usepackage{physics}
\usepackage{subcaption}
\usepackage{tikz,tikz-3dplot}
\usepackage{xfrac}
\usepackage{xcolor}

\AtBeginDocument{\RenewCommandCopy\qty\SI}
\ExplSyntaxOn
\msg_redirect_name:nnn { siunitx } { physics-pkg } { none }
\ExplSyntaxOff

\graphicspath{{./figs/}}

\definecolor{isaeaquamarine}{HTML}{26b0a0}
\definecolor{isaecoral}{HTML}{e95e40}
\definecolor{isaeslategray}{HTML}{516d6d}
\definecolor{isaeyellowgreen}{HTML}{97bf0d}
\definecolor{RoyalBlue}{HTML}{4169E1}

\hypersetup{colorlinks=true,allcolors=RoyalBlue}
\glsdisablehyper

\setabbreviationstyle[s]{short}
\setabbreviationstyle[ls]{long-short}
\setabbreviationstyle[l]{long}
\setabbreviationstyle[sl]{short-long}
\loadglsentries[acronym]{acronyms.tex}

\newcommand{\vc}[1]{\vb*{#1}}
\newcommand{\da}[1]{\left[#1\right]}
\newcommand{\ty}[3]{\mathcal{T}^{(#1)}_{#2}\left(#3\right)}
\newcommand{\ex}[1]{\mathbb{E}\left[#1\right]}

\makeatletter
\let\oldappendix\appendix
\renewcommand\appendix{\oldappendix\gdef\thesection{\appendixname~\@Alph\c@section}}
\makeatother

\title{Efficient Multifidelity Uncertainty Propagation in the Presence of Process Noise}

\author{Alberto Foss\`a\footnote{PhD Candidate, Department of Aerospace Vehicles Design and Control (DCAS). Currently Postdoctoral Fellow, Oden Institute for Computational Engineering and Sciences, The University of Texas at Austin, Austin, Texas 78712; \href{mailto:alberto.fossa@utexas.edu}{alberto.fossa@utexas.edu} (Corresponding Author)}\orcidlink{0000-0002-0756-4998}}
\affil{Higher Institute of Aeronautics and Space (ISAE-SUPAERO), Toulouse, 31055, France}
\author{Roberto Armellin\footnote{Full Professor, Te P\=unaha \=Atea -- Space Institute.}\orcidlink{0000-0002-3516-6428}}
\affil{The University of Auckland, Auckland, 1010, New Zealand}
\author{Emmanuel Delande\footnote{Space Surveillance Specialist, Space Surveillance and Tracking.}\orcidlink{0000-0001-8819-285X}}
\affil{National Center for Space Studies (CNES), Toulouse, 31401, France}
\author{Francesco Sanfedino\footnote{Associate Professor, Department of Aerospace Vehicles Design and Control (DCAS).}\orcidlink{0000-0002-3732-5357}}
\affil{Higher Institute of Aeronautics and Space (ISAE-SUPAERO), Toulouse, 31055, France}

\begin{document}

\maketitle

\begin{abstract}

A \gls{mf} method for the nonlinear propagation of uncertainties in the presence of stochastic accelerations is presented. The proposed algorithm treats the \gls{up} problem by separating the propagation of the uncertainty in the \glspl{ic} and static model parameters from that of the \gls{pn}. The first two types of uncertainty are propagated using an adaptive \gls{gmm} method which exploits a \gls{lf} dynamical model to minimize the computational costs. The effects of \gls{pn} are instead computed using the \gls{plasma} technique, which considers a \gls{hf} model of the stochastic dynamics. The main focus of the paper is on the latter and on the key idea to approximate the \gls{pdf} of the solution by a polynomial representation of its moments, which are efficiently computed using \gls{da} techniques. The two estimates are finally combined to restore the accuracy of the \gls{lf} surrogate and account for all sources of uncertainty. The proposed approach is applied to the problem of nonlinear orbit \gls{up} and its performance compared to that of \gls{mc} simulations.

\end{abstract}

\glsresetall

\section{Introduction}\label{sec:introduction}

\lettrine{A}{ccurate} propagation of uncertainties is key for \gls{ssa} applications. Ranging from \gls{sst} to \gls{cam} planning, these activities require the estimation of the spacecraft state at a future epoch given the \glspl{ic} at an earlier time. These \glspl{ic} are usually the result of \gls{od} techniques, which are stochastic in nature. Since the output of the latter is a full state estimate - i.e., some \gls{pdf} - it is necessary to propagate as much, rather than a nominal state, to any epoch of interest. Given the intrinsic nonlinearity of the dynamics, no analytical solution exists to map the aforementioned \gls{pdf}. On the contrary, several methods have been developed to approximate this quantity with increasing level of accuracy, usually at the expenses of higher computational efforts. These approaches tackle the \gls{up} problem in different ways, and can be grouped into two main categories: linear methods and nonlinear methods. The second category can then be further subdivided into sample-based, dynamics-based, and \gls{pdf}-based techniques~\cite{Luo2017}.

Linear methods build on the assumptions that a linearized dynamical model is sufficient to capture the dynamics of neighboring trajectories and that the uncertainty can be completely characterized by a multivariate Gaussian distribution. Under these hypotheses, only the mean and the covariance matrix need to be propagated, and the problem reduces to the integration of the \gls{stm}, which is known as \gls{lc} propagation. However, if the dynamics is highly nonlinear, this approach fails to accurately describe the time evolution of the state \gls{pdf}. An extension of the \gls{lc} method is that of \glspl{stt}, which are higher-order Taylor expansions of the dynamics about the nominal trajectory~\cite{Park2006a}. Yet, this method is not ideal for complex dynamics at high orders, since it requires the derivation of analytical expressions for the  variational equations and their numerical integration along with the reference trajectory.

On the other side, \gls{mc} simulations are widely used in operational scenarios since they provide highly accurate estimation of the state uncertainty and are easy to implement. Yet their accuracy comes with a major computational cost, which makes \gls{mc} techniques not applicable for real-time applications or maintenance of very large space objects catalogs. Sample based methods also include \gls{ut}~\cite{Julier2004} and \gls{cut}~\cite{Adurthi2015,Adurthi2018} techniques, which are deterministic in nature as opposed to stochastic \gls{mc} simulations. These methods are based on the idea that is \textit{easier to approximate a state \gls{pdf} than an arbitrary nonlinear transformation}~\cite{Julier2004}. Assuming a Gaussian distribution for the initial uncertainty, they enforce the corresponding \glspl{mce} up to a given order to solve for a deterministic set of samples, denoted as sigma points, whose weighted statistics correctly capture the first statistical moments of the state \gls{pdf}. To reduce the search space, these points are constrained to lie on a carefully selected set of axes, and a fully symmetric set is chosen to automatically satisfy the \glspl{mce} for any odd moment of the Gaussian distribution. If the \gls{ut} can only capture the statistical moments up to third order, \gls{cut} sigma points sets were obtained to match these moments up to ninth order. However, this approach suffers from a rapid increase in the number of samples needed to satisfy the \glspl{mce} for increasing orders and state dimension.

\Gls{pce}, in its non-intrusive formulation, is another sample-based technique developed for \gls{up}~\cite{Jones2013}. This method aims at obtaining a functional representation of the propagated uncertainty with respect to the input random variables, thus providing accurate information on higher-order statistical moments of the state \gls{pdf}. Inputs and outputs of the considered transformation are approximated as series expansions of orthogonal polynomials for which the coefficients are sought. This technique requires a number of \gls{pce} terms which varies exponentially with both the order of the polynomial basis and the dimensions of the input random variables, leading to the curse of dimensionality for large problems.

Moving further, a different approach to the \gls{up} problem is that of dynamics-based methods, which include the previously mentioned \gls{lc} and \glspl{stt}, intrusive \gls{pce}, and \gls{da}. In contrast to its non-intrusive counterpart, the intrusive \gls{pce} method does not require sampling the initial distribution. Instead, it performs a Galerkin projection of the governing dynamics onto a subspace to transform a \gls{sode} in low-dimension into an equivalent deterministic \gls{ode} in higher dimension. However, given the tedious modifications to the governing equations that are needed to obtain the higher-dimensional \gls{ode}, this approach is usually less appealing than the non-intrusive one for complex dynamics. \Gls{da} is instead equivalent to \glspl{stt}, as they both provide the $k^{th}$-order Taylor expansion of the flow of the dynamics around the nominal solution. The main advantage of \gls{da} over \glspl{stt} is that the former does not require the derivation of the variational equations nor the numerical integration of an augmented state, thus making this technique more flexible and computationally efficient as already demonstrated in the context of nonlinear orbit \gls{up}~\cite{Valli2013,Armellin2018}. A key difference between \gls{da} (or \glspl{stt}) and \gls{pce} is that the former represent the uncertainty using Taylor polynomials, while the latter uses orthogonal polynomials. Although Taylor expansions are highly accurate near the expansion point, their accuracy decreases rapidly with the distance from the latter. In contrast, orthogonal polynomials provide a more uniform accuracy across the entire domain of interest while begin less accurate near the expansion point. \Gls{ads} techniques were thus developed to mitigate the limitations of a single Taylor polynomial~\cite{Wittig2015a}. The idea is to continuously monitor the accuracy of the Taylor expansion by estimating the magnitude of the truncated terms and react consequently when the imposed error threshold is crossed. When this happens, the single polynomial is split into two new expansions, each of them covering only half of the initial domain, and the process is repeated until the truncation error is below the select threshold in each subdomain. This algorithm requires an expansion order $k\geq 3$ and its accuracy increases for increasing orders. More recently, a novel splitting algorithm, named \gls{loads}, was instead developed specifically for second-order Taylor expansions~\cite{Losacco2024}. A merging scheme was also proposed to reduce the number of domains whenever possible.

The last category of nonlinear \gls{up} methods is that of \gls{pdf}-based techniques, which provide a functional representation of the uncertainty as opposed to a set of samples or statistical moments that characterize the \gls{pdf}. The first example is the solution of the \gls{fpe}. For any dynamical system that satisfies the It\^o \gls{sde}, the \gls{fpe} describes the time evolution of the \gls{pdf} associated to the state of the system~\cite{Fuller1969}. It is a \gls{pde} that satisfies the propagation of a \gls{pdf}, and its solution provides a complete statistical description of the trajectory governed by the system's \gls{sde}~\cite{Luo2017}. However, this \gls{pde} has no analytical solution in the general case, and a numerical approximation of the latter requires the discretization of the probability space, which leads to the curse of dimensionality for high-dimensional problems. Solution methods based on meshless discretization~\cite{Kumar2010} and tensor decomposition~\cite{Sun2014,Sun2016} were proposed to overcome this limitation, but they remain unfeasible in most practical scenarios. Further computational efficiency is still obtained in those problems where the diffusion term of the \gls{fpe} can be neglected. In these cases, the \gls{fpe} reduces to the Liouville's equation, which is a quasi-linear \gls{pde} that describes the time evolution of the phase space distribution function~\cite{Weisman2016}. Using the transformation of variables technique, closed form solutions to this equation were obtained for the unperturbed Keplerian motion~\cite{Majji2012} and the $J_2$ perturbed problem~\cite{Weisman2016}. An equivalent technique, the method of characteristics, was also applied to study the behavior of debris clouds~\cite{Letizia2015,Letizia2016}, high-altitude and high-area-to-mass ratio satellites~\cite{Sun2022a}, and reentry uncertainties~\cite{Trisolini2021a}. However, if the underlying system of \glspl{ode} does not admit an analytical solution, the density must be computed numerically for a number of discrete trajectories, and further processing is needed to retrieve a continuous representation of the \gls{pdf}. Although techniques that provide accurate results were developed~\cite{Trisolini2021a}, they also add an extra layer of complexity that makes this approach better suited for long-term density propagation rather than short-term or real-time \gls{up}. The second example of \gls{pdf}-based methods are those that represent the uncertainty using \glspl{gmm}. \Glspl{gmm} can approximate any \gls{pdf} of practical concern using a weighted sum of Gaussian distributions, and this approximation converges to the true \gls{pdf} as the number of components approaches infinity~\cite{Sorenson1971}. Building on this observation, several hybrid methods were proposed to propagate the kernels' statistics using the \gls{ekf}, the \gls{ukf}~\cite{DeMars2013}, \glspl{pce}~\cite{Vittaldev2016}, and \gls{da}~\cite{Sun2019a,Servadio2021a,Servadio2021c,Losacco2024}. The main challenge, however, is the estimation of the optimal number of mixture components for which different approaches have been investigated. \Citet{Horwood2011} generate a large number of kernels \emph{a priori}, while the \gls{aegis} method developed by \citet{DeMars2013} uses the error between the linearized and nonlinear differential entropy as a measure to adaptively trigger new splits. On the contrary, \citet{Terejanu2008} keep the number of kernels constant, but propose two approaches to continuously update their weights. The first one is based on the integral square difference between the true forecast \gls{pdf} and its \gls{gmm} approximation, while the second uses the \gls{fpke} error as feedback. This method was then extended to prune those components with a negligible impact on the overall state \gls{pdf} and merge kernels whose combined contributions resemble a Gaussian~\cite{Terejanu2011}. Similarly, \citet{Vishwajeet2018} use the Kolmogorov equation error to trigger new splits and merge redundant components. All these methods require the solution of an augmented problem, as additional equations are needed to estimate the nonlinearity of the transformation. Depending on the implementation, this may result in a computational overhead that cannot be afforded in practice. To overcome this drawback, \citet{Losacco2024} proposed an adaptive \gls{gmm} method that exploits a novel measure of nonlinearity to inform the splitting and merging processes. This method is particularly resource-effective, as the second-order Taylor expansion of the transformation is also used to estimate the kernels' statistics via efficient polynomial evaluation. Although not being the only method that leverages \gls{da} to propagate the kernels' means and covariances~\cite{Sun2019a,Servadio2021a,Servadio2021c}, this approach is the first one to exploit the information contained in the Taylor expansion of the transformation to generate new kernels on demand. Indeed, \citet{Sun2019a} split the initial distribution into a fix number of components and use a single polynomial expansion of the dynamics to map their statistics over time. \Citet{Servadio2021a} generate instead a mixture representation of the initial Gaussian distribution in analogy with the \gls{ut} theory~\cite{Julier2004}. Each kernel mean is then propagated independently in the \gls{da} framework to map the \gls{gmm} means and covariances through the system's dynamics. Similarly, \citet{Servadio2021c} leverage the Taylor expansions of the kernels' means to map a set of samples drawn from the initial distribution. A \gls{gmm} representation of the final \gls{pdf} is then obtained from these samples using a combination of $k$-means and \gls{exmx} algorithms. All three methods, however, maintain the number of kernels constant throughout the propagation. Their behavior is thus in sharp contrast with the adaptive method proposed by \citet{Losacco2024} and extended in this work. In essence, the latter aligns more closely in spirit with the \gls{aegis} method~\cite{DeMars2013} than with the aforementioned \gls{da}-based techniques.

All methods cited above describe the evolution of the orbit state using a deterministic dynamical model. The focus is thus solely on the propagation of the uncertainty in the \glspl{ic} and static model parameters. Nevertheless, the uncertainty due to poorly modeled or unmodeled forces acting on the \gls{so} can have a substantial impact on the evolution of its state. If not accurately accounted for, this impact can cause a loss of uncertainty realism. These perturbations are usually introduced in the dynamics as stochastic accelerations, and can be modeled as random processes such as Wiener processes or Gauss-Markov processes. A widespread algorithm for the treatment of \gls{pn} is the \gls{snc} method~\cite{Tapley2004}. This algorithm, however, relies on the hypothesis that stochastic accelerations are modeled as uncorrelated white noise, and provides a solution to the \gls{up} problem in the form of an expected state and a covariance matrix. Alternatively, the \gls{dmc} method can be used to account for other forms of \gls{pn}, but this algorithm reduces to the \gls{snc} method in an augmented state space~\cite{Tapley2004}.

Another key aspect in nonlinear \gls{up} is the trade-off between accuracy and computational cost. Indeed, the \gls{up} algorithm and the dynamical model that describes the evolution of the system must be tailored to the available computational resources, which could result in a loss of realism for the propagated uncertainty. Similarly, applications that require highly accurate results might face excessive costs due to the important computational power allocated to \gls{up}. Recent works thus focused on \gls{mf} methods to reduce the execution time of this task while limiting the accuracy losses~\cite{Peherstorfer2016,Jones2018}. These methods assume that multiple mathematical models are available to describe the system of interest, and that each model has different evaluation cost and fidelity. The key idea is to delegate most evaluations to a \gls{lf} model, which has a low accuracy but is also cheap to evaluate, while performing only few evaluations of the more expensive \gls{hf} one to improve the accuracy of the propagated estimate. \Citet{Peherstorfer2016} provides a broad survey of \gls{mf} methods with applications in \gls{up}, inference and optimization. Methods based on control variates, importance sampling, and \gls{sc} are mentioned as alternatives to efficiently solve the \gls{up} problem. Moreover, \citet{Jones2018} propose a \gls{mf} method for orbit \gls{up} that uses \gls{sc} to propagate the state \gls{pdf} represented as a particle ensemble or as a \gls{gmm}.

This work addresses the aforementioned challenges with two main contributions: a novel method for the propagation of uncertainties in the presence of stochastic accelerations, and a \gls{mf} technique for efficient nonlinear \gls{up}. The first algorithm, named \gls{plasma}, maintains a polynomial approximation of the moments of the state \gls{pdf} subject to the effects of \gls{pn}. Although this method builds upon the work of \citet{LopezYela2023}, it overcomes two major limitations of the original algorithm. Firstly, it uses \gls{da} instead of \glspl{stt} to compute an approximation of the \gls{pdf} moments, thus removing the need for the analytical derivation of the Taylor expansion of the dynamics with respect to the state and the \gls{pn}. Secondly, it does not assume that the solution to the \gls{sde} that governs the dynamical system of interest is approximated using the \gls{em} method. Instead, the \gls{plasma} algorithm is extended to any explicit numerical integration method for nonlinear \glspl{sde} including adaptive ones. Previous works~\cite{Servadio2021a,Servadio2021b} also exploited \gls{da} to obtain a Taylor expansion of the final state function of the \glspl{ic} and the \gls{pn}. In this case, however, the propagated dynamics are deterministic, and the \gls{pn} is accounted for with a linear term added to the final state. Moreover, the computation of the \gls{pdf} moments assumes a Gaussian distribution for all independent \gls{da} variables, which is not the case in this work. The performance of the \gls{plasma} method can be further enhanced by precomputing a \gls{pw} solution to the deterministic dynamics. This solution can then be utilized as a reference trajectory, around which the \gls{pdf} moments are estimated using a simplified model for the deterministic dynamics. This is especially useful when the deterministic part of the dynamics is considerably more expensive to evaluate than the stochastic accelerations. The second contribution of this work is a \gls{bf} method for nonlinear \gls{up} that exploits the adaptive \gls{gmm} method proposed by \citet{Losacco2024} as workhorse for the propagation of the uncertainty in the \glspl{ic} and static model parameters. The latter considers a \gls{lf} model of the dynamics, and it is followed by a \gls{hf}, \gls{pw} propagation of the polynomials' constant parts to restore the accuracy of the \gls{lf} solution. This technique is also extended to the treatment of \gls{pn} effects by replacing the \gls{pw} propagations with the \gls{plasma} algorithm just introduced. In this case, the \gls{lf} solution is not only corrected for the drift due to the mismatch between the \gls{lf} and \gls{hf} models, but also for the effects of stochastic accelerations estimated by the \gls{plasma} method.

The reminder of this paper is organized as follows. \Cref{sec:background} sets the mathematical background on which the subsequent material is based. The \gls{plasma} algorithm and its \gls{bf} formulation are discussed in \cref{sec:plasma}. Moving further, the \gls{mf} techniques that couple the adaptive \gls{gmm} method with either \gls{pw} propagations in \gls{hf} dynamics or the \gls{plasma} algorithm are presented in \cref{sec:multifidelity}. The performance of the proposed methods is finally assessed in \cref{sec:applications} by applying these techniques to the problem of orbit \gls{up} and comparing their results to those of \gls{mc} simulations. Conclusions are drawn in \cref{sec:conclusions}. Key properties of the \gls{plasma} algorithm are highlighted in \ref{apdx:duffing}, where it is applied to the discrete Duffing oscillator.

\section{Mathematical Background}\label{sec:background}

This section provides the mathematical background on which the \gls{plasma} algorithm and the \gls{mf} \gls{up} methods are based. This includes the \gls{da} techniques in \cref{sec:background:da}, the \gls{da}-based \gls{nli} in \cref{sec:background:nli}, and the adaptive \gls{gmm} method for \gls{up} in \cref{sec:background:gmm}.

\subsection{Differential Algebra}\label{sec:background:da}

\Gls{da} is a computing technique that stems from the idea that it is possible to extract more information from a function $f$ than its mere value $f(x)$ at a point $x$. Given any function $\vc{f}:\mathbb{R}^n\to\mathbb{R}^m$ that is $\mathcal{C}^{k+1}$ in the domain of interest $\mathcal{D}=[-1,1]^n$, the algebra of \gls{fp} numbers is replaced by a new algebra of Taylor polynomials to compute the $k^{th}$-order expansion of $\vc{f}$~\cite{Berz1999}. This algebra is denoted as $_kD_n$, where $k$ is the order of expansion and $n$ the number of independent variables. Its dimension, i.e. the number of monomials of order $k$ in $n$ variables, is given by
\begin{equation}
    \mathrm{dim}(_kD_n) = \binom{n+k}{k} = \binom{n+k}{n} = \dfrac{(n+k)!}{n!\,k!}.
    \label{eq:da_dim}
\end{equation}
Throughout the paper, the $k^{th}$-order Taylor expansion of $\vc{f}$ in terms of the $n$ independent \gls{da} variables $\var{\vc{x}}=\{\var{x_1},\ldots,\var{x_n}\}^T$ is denoted by
\begin{equation}
    \vc{f}\approx\da{\vc{f}}=\ty{k}{\vc{f}}{\var{\vc{x}}},
\end{equation}
where the superscript ``$(k)$'' indicates a $k^{th}$-order Taylor expansion in $\var{\vc{x}}$.

The four arithmetic operations ($+,-,\ast,\divisionsymbol$), some elementary functions such as exponential, logarithm and trigonometric functions, as well as derivation, integration, map composition and inversion, are all well defined in \gls{da}. These basic operations can then be combined to derive powerful algorithms for several applications. Common examples are the solution of implicit equations, the computation of the Taylor expansion of the flow of the dynamics in terms of their \glspl{ic}~\cite{Valli2013}, and the solution of \glspl{bvp}~\cite{Armellin2018h}.

When working with physical quantities described by some state $\vc{x}$, it is convenient to introduce scaling factors $\vc{\beta}$ such that the deviations $\var{\vc{x}}\in\mathcal{D}$ are mapped onto the physical domain of interest for $\vc{x}$. In this case, $\vc{x}$ is initialized in the \gls{da} framework as
\begin{equation}
    \da{\vc{x}} = \bar{\vc{x}} + \vc{\beta}\odot\var{\vc{x}}
    \label{eq:da_scl_var_def},
\end{equation}
where $\bar{\vc{x}}$ is the nominal value of $\vc{x}$, $\vc{\beta}\in\mathbb{R}^n_{\geq 0}$ are the nonnegative scaling factors, $\odot$ denotes the Hadamard (or element-wise) product, and $\var{\vc{x}}\in\mathcal{D}$ are the $n$ independent \gls{da} variables.

The parameters $\vc{\beta}$ may have different physical meanings. For instance, if $\vc{X}\in\mathbb{R}^n$ is a multivariate random variable with \gls{pdf} $p_{\vc{X}}$, it is of interest to define $\da{\vc{x}}$ such that it represents a domain that captures a sizeable share of the probability mass. Assuming that the first two moments of $\vc{X}$, i.e. its mean $\vc{\mu}_{\vc{X}}$ and covariance matrix $\vc{P}_{\vc{X}}$, are finite and known, a natural choice to is to initialize the independent \gls{da} variables in the eigenspace of $\vc{P}_{\vc{X}}$, and to center the Taylor expansion in $\vc{\mu}_{\vc{X}}$. In this case, the polynomial $\da{\vc{x}}$ is defined as
\begin{equation}
    \da{\vc{x}} = \vc{\mu}_{\vc{X}}+\vc{V}\left[\zeta\sqrt{\vc{\lambda}}\odot\var{\vc{x}}\right]
    \label{eq:da_mu_cov_def},
\end{equation}
where $\vc{V}$ and $\vc{\lambda}$ are the eigenvectors and eigenvalues of $\vc{P}_{\vc{X}}$, respectively, and $\zeta>0$ is a nonnegative coefficient. The scaling factors $\vc{\beta}$ are thus defined as $\vc{\beta}=\zeta\sqrt{\vc{\lambda}}$ in this case. \Cref{eq:da_mu_cov_def} does not require any knowledge of $\vc{X}$ besides its first two moments. In particular, it does not require $\vc{X}$ to be Gaussian distributed. Indeed, the objective of \cref{eq:da_mu_cov_def} is to capture a significant portion of the probability mass of $\vc{X}$ within the region defined by a first-order polynomial, and not to provide an accurate representation of the \gls{pdf} $p_{\vc{X}}$. Nevertheless, if $\vc{X}$ is Gaussian distributed, then $\da{\vc{x}}$ represents a domain spanning $\zeta$ standard deviations around the mean $\vc{\mu}_{\vc{X}}$.

\subsection{Nonlinearity Index}\label{sec:background:nli}

To monitor the transition between linear and nonlinear regimes, a \gls{da}-based measure of nonlinearity was firstly proposed by \citet{Losacco2024}. The \gls{nli} measure $\nu$ exploits the second-order Taylor expansion of the transformation $\vc{f}$ to quantify its deviation from linearity, and is computed as follows.

Consider a nonlinear transformation $\vc{f}:\mathbb{R}^n\to\mathbb{R}^m$ that is $\mathcal{C}^2$ over its domain and an uncertainty set $\da{\vc{x}}\in{_2D_n}$. The image of $\da{\vc{x}}$ through $\vc{f}$ is firstly obtained in the \gls{da} framework as
\begin{equation}
    \da{\vc{y}} = \vc{f}\left(\da{\vc{x}}\right) = \ty{2}{\vc{f}}{\var{\vc{x}}},
    \label{eq:da_image}
\end{equation}
where $\var{\vc{x}}$ are the $n$ independent \gls{da} variables. The Jacobian of $\vc{f}$ with respect to $\da{\vc{x}}$ is then given by
\begin{equation}
    \da{\vc{J}} = \pdv{\da{\vc{y}}}{\vc{\beta}\odot\var{\vc{x}}} = \overline{\vc{J}} + \var{\vc{J}},
    \label{eq:jacobian}
\end{equation}
which is an $m\times n$ matrix of first-order polynomials in $\var{\vc{x}}$. The \gls{nli} is finally obtained as
\begin{equation}
    \nu = \dfrac{\norm{\vc{B}_{\var{\vc{J}}}}_2}{\norm{\overline{\vc{J}}}_2},
    \label{eq:nli}
\end{equation}
where $\overline{\vc{J}}$ is the constant part of $\da{\vc{J}}$, $\vc{B}_{\var{\vc{J}}}$ are the upper bounds of $\var{\vc{J}}$, and $\norm{\cdot}_2$ denotes the Frobenius norm. Given that each entry of $\var{\vc{J}}$ is a first-order polynomial with zero constant part, its upper bound $B_{\var{J_{i,j}}}$ is given by
\begin{equation}
    B_{\var{J_{i,j}}} = \sum_{p=1}^n \abs{c_{i,j}^p},
    \label{eq:jacobian_bound}
\end{equation}
with $c_{i,j}^p$ the coefficients of $\var{J_{i,j}}$. Moreover, matrix norms are most meaningful if both $\vc{x}$ and $\vc{y}$ are nondimensional quantities such that their constant coefficients are all of order \num{1}. Then \cref{eq:nli} is itself a nondimensional measure of nonlinearity, it is in fact a matrix relative error of the largest variation of the Jacobian with respect to its constant part. An estimate of the nonlinearity of $\vc{f}$ over $\da{\vc{x}}$ is thus obtained by inspecting the value of $\nu$. This index is equal to zero for linear transformations, and grows with nonlinearity effects.

\subsection{Adaptive Gaussian Mixture Models}\label{sec:background:gmm}

Consider a nonlinear transformation $\vc{f}:\mathbb{R}^n\to\mathbb{R}^m$ that is $\mathcal{C}^2$ over its domain and a multivariate random variable $\vc{X}$ defined on $\mathbb{R}^n$ with \gls{pdf} $p_{\vc{X}}$. Moreover, assume that a \gls{gmm} representation of $p_{\vc{X}}$ is available as
\begin{equation}
    p_{\vc{X}}(\vc{x}) \approx \sum_{l=1}^{L_0} \alpha^{(l)} p_g\left(\vc{x};\vc{\mu}_{\vc{X}}^{(l)},\vc{P}_{\vc{X}}^{(l)}\right),
    \label{eq:gmm_pdf}
\end{equation}
where $\alpha_l$ are the mixture weights, $\vc{\mu}_{\vc{X}}^{(l)},\vc{P}_{\vc{X}}^{(l)}$ are the kernels' means and covariances, respectively, and $L_0$ is the initial number of kernels. To compute a \gls{gmm} representation of the propagated \gls{pdf} $p_{\vc{Y}}$, the algorithm proposed by \citet{Losacco2024} starts by initializing a Taylor polynomial for each kernel in \cref{eq:gmm_pdf} using \cref{eq:da_mu_cov_def}. The image of each polynomial through $\vc{f}$ is then obtained in the \gls{da} framework, and the \gls{nli} is computed to estimate the truncation error. If $\nu$ is below a fixed threshold $\varepsilon_\nu$, the resulting polynomial is kept; otherwise, it is discarded. In this second case, the initial kernel and associated polynomial are split into three components such that the domain of each new polynomial is a fraction of the original one. The procedure is then repeated until the condition $\nu\leq\varepsilon_\nu$ is satisfied within each subdomain. These splits are carried out according to predefined splitting libraries that are designed to minimize the error between original and split \glspl{pdf}~\cite{DeMars2013}. Once the recursive procedure has converged, the outputs are two manifolds that represent the initial and propagated distributions, respectively. They are denoted as
\begin{subequations}
    \begin{align}
        \mathcal{M}_{\vc{X}} &= \left\{\left(\da{\vc{x}^{(l)}},\alpha^{(l)},\vc{\mu}_{\vc{X}}^{(l)},\vc{P}_{\vc{X}}^{(l)}\right)\right\}_{l=1}^L \label{eq:man_init} \\
        \mathcal{M}_{\vc{Y}} &= \left\{\left(\da{\vc{y}^{(l)}}\right)\right\}_{l=1}^L \label{eq:man_out},
    \end{align}
\end{subequations}
with $L$ the number of kernels after splitting. \Cref{eq:man_init} already contains all the information required to obtain a refined \gls{gmm} representation of the initial \gls{pdf} $p_{\vc{X}}$ as
\begin{equation}
    p_{\vc{X}}(\vc{x}) \approx \sum_{l=1}^L \alpha^{(l)} p_g\left(\vc{x};\vc{\mu}_{\vc{X}}^{(l)},\vc{P}_{\vc{X}}^{(l)}\right).
    \label{eq:gmm_pdf_refined}
\end{equation}
For each kernel in the mixture approximation of $p_{\vc{Y}}$, its mean $\vc{\mu}_{\vc{Y}}^{(l)}$ and covariance $\vc{P}_{\vc{Y}}^{(l)}$ are then estimated by drawing a set of \gls{ut} sigma points from the Gaussian distribution with mean $\vc{\mu}_{\vc{X}}^{(l)}$ and covariance $\vc{P}_{\vc{X}}^{(l)}$, evaluating the corresponding polynomial $\da{\vc{y}^{(l)}}$ on these points, and computing the weighted mean and covariance of the resulting set according to the \gls{ut} formulae~\cite{Losacco2024,Fossa2024b}. The weights $\alpha^{(l)}$ are instead kept constant in this work. Indeed, although techniques to update the kernels' weights during propagation have been proposed~\cite{Terejanu2008,Terejanu2011a}, these methods are often computationally expensive and might degrade the quality of the estimate~\cite{Horwood2011}. The \gls{gmm} representation of $p_{\vc{Y}}$ is thus given by
\begin{equation}
    p_{\vc{Y}}(\vc{y}) \approx \sum_{l=1}^L \alpha^{(l)} p_g\left(\vc{y};\vc{\mu}_{\vc{Y}}^{(l)},\vc{P}_{\vc{Y}}^{(l)}\right),
    \label{eq:gmm_pdf_prop}
\end{equation}
where $\vc{\mu}_{\vc{Y}}^{(l)},\vc{P}_{\vc{Y}}^{(l)}$ is the \gls{ut} of $\vc{\mu}_{\vc{X}}^{(l)},\vc{P}_{\vc{X}}^{(l)}$ through $\da{\vc{y}^{(l)}}$, and the weights $\alpha^{(l)}$ are those in \cref{eq:man_init}. This method was demonstrated to be particularly effective in the context of nonlinear orbit \gls{up}~\cite{Losacco2024}, and it constitutes one of the pillars of the \gls{mf} technique presented in \cref{sec:multifidelity}.

\section{The PLASMA Method}\label{sec:plasma}

The first contribution of this paper, namely, the \gls{plasma} method, is introduced hereafter. The key idea is to describe the \gls{pdf} of the stochastic process solution to a \gls{sde} by maintaining a polynomial representation of the moments of said \gls{pdf} (e.g. mean, covariance, skewness, kurtosis) up to a preselected order $N$. To achieve this objective, the state at time $t_k$ is expanded in Taylor series function of the state at time $t_{k-1}$ and the realization of the \gls{pn} at step $k$. Combinatorics and the linearity of expectations are then exploited to compute the moments of the state at $t_k$ given those of the state at $t_{k-1}$ and that of the \gls{pn}. By applying this algorithm recursively, the moments of the final state are thus obtained from the sole knowledge of the initial moments at $t_0$, the \gls{pn} model, and the dynamical model. The theoretical foundations of the algorithm are firstly presented in \cref{sec:plasma:moments}, followed by its \gls{da} formulation in \cref{sec:plasma:da}. The use of \gls{mf} techniques within the \gls{plasma} method itself are finally discussed in \cref{sec:plasma:multifidelity}.

\subsection{Polynomial expansion of the effective noise}\label{sec:plasma:moments}

Consider the \gls{ivp} for the following \gls{sde}
\begin{equation}
    \begin{dcases}
        \dd{\vc{X}(t)} &= \vc{u}(\vc{X},t)\dd{t} + \vc{G}(\vc{X},t)\dd{\vc{W}(t)} \\
        \vc{X}(t_0) &= \vc{X}_0,
    \end{dcases}
    \label{eq:sde_def}
\end{equation}
with $\vc{u}:\mathbb{R}^n\times\mathbb{R}\to\mathbb{R}^n$ the drift coefficient, $\vc{G}:\mathbb{R}^n\times\mathbb{R}\to\mathbb{R}^n\times\mathbb{R}^m$ the diffusion coefficient, $\vc{W}(t)$ the $m$-dimensional Wiener process with independent increments, and $\vc{X}_0$ the $n$-dimensional random variable representing the \glspl{ic}. The \gls{em} scheme can be employed to obtain an approximate solution to \cref{eq:sde_def} over a predefined time grid $\left(t_k\right)_{k=1}^M$, where $M$ is the number of steps. The discretization of the \gls{sde} defined by \cref{eq:sde_def} leads to
\begin{equation}
    \begin{dcases}
        \hat{\vc{X}}_k &= \hat{\vc{X}}_{k-1} + h\vc{u}(\hat{\vc{X}}_{k-1},t_{k-1}) + \vc{G}(\hat{\vc{X}}_{k-1},t_{k-1})\Delta\vc{W}_k\\
        t_{k} &= t_{k-1} + h,
    \end{dcases}
    \label{eq:em_scheme}
\end{equation}
with $h$ the integration step size, and $\Delta\vc{W}_k$ the increment of the stochastic process $\vc{W}(t)$ in the time interval $[t_{k-1},t_k]$. The random sequence in \cref{eq:em_scheme} can then be rewritten as
\begin{equation}
    \hat{\vc{X}}_k = \hat{\vc{x}}^C_k + \Delta\hat{\vc{W}}_k,
    \label{eq:split_sequence}
\end{equation}
where $\hat{\vc{x}}^C_k$ is the deterministic sequence (or central part) obtained by neglecting the diffusion term in \cref{eq:em_scheme}, and $\Delta\hat{\vc{W}}_k$ is the effective noise random process. At each time step, $\Delta\hat{\vc{W}}_k$ can be approximated as a polynomial expansion of order $N$ in $\left\{\Delta\hat{\vc{W}}_{k-1},\Delta\vc{W}_k\right\}$ leading to the following recursive equation $\forall\,i\in[1,n]$~\cite{LopezYela2023}
\begin{equation}
    \begin{split}
            \Delta\hat{W}^{(i)}_{k,N}(\vc{x}_0) &= \Delta\hat{W}^{(i)}_{k-1,N}(\vc{x}_0)\\ &+ h\sum_{\abs{\vc{r}}=1}^N\dfrac{1}{\vc{r}!}\dfrac{\partial^{\abs{\vc{r}}}}{\partial\vc{x}_{k-1}^{\vc{r}}}u^{(i)}\left(\hat{\vc{x}}^C_{k-1}(\vc{x}_0),t_{k-1}\right)\Delta\hat{\vc{W}}^{\vc{r}}_{k-1,N}(\vc{x}_0)\\ &+ \sum_{\abs{\vc{r}}=0}^{N-1}\sum_{j=1}^m\dfrac{1}{\vc{r}!}\dfrac{\partial^{\abs{\vc{r}}}}{\partial\vc{x}_{k-1}^{\vc{r}}}G^{(i,j)}\left(\hat{\vc{x}}^C_{k-1}(\vc{x}_0),t_{k-1}\right)\Delta W_k^{(j)}\Delta\hat{\vc{W}}^{\vc{r}}_{k-1,N}(\vc{x}_0),
    \end{split}
    \label{eq:effective_noise_expansion}
\end{equation}
with $\Delta\hat{\vc{W}}=\left(\Delta\hat{W}^{(1)},\ldots,\Delta\hat{W}^{(n)}\right)$, $\Delta \vc{W}=\left(\Delta W^{(1)},\ldots,\Delta W^{(m)}\right)$ and the multi-index notation introduced in \ref{apdx:indices} is used. The conditional moments of $\Delta\hat{\vc{W}}_{k,N}(\vc{x}_0)$ are then obtained by taking the expectation of \cref{eq:effective_noise_expansion} as
\begin{equation}
    \ex{\Delta\hat{\vc{W}}^{\vc{r}}_{k,N}(\vc{x}_0)} = \sum_{\abs{\vc{s}}+\abs{\vc{r}'}=1}^N a_{\vc{r},k,N}^{\vc{s},\vc{r}'}\left(\hat{\vc{x}}^C_{k-1}(\vc{x}_0),t_{k-1}\right)\ex{\Delta\vc{W}^{\vc{s}}_{k}}\ex{\Delta\hat{\vc{W}}^{\vc{r}'}_{k-1,N}(\vc{x}_0)},
    \label{eq:effective_noise_moments}
\end{equation}
where $a_{\vc{r},k,N}^{\vc{s},\vc{r}'}\left(\hat{\vc{x}}^C_{k-1}(\vc{x}_0),t_{k-1}\right)$ are the Taylor coefficients obtained from the polynomial expansion of $\Delta\hat{W}^{(i)}_{k,N}(\vc{x}_0)$. Using the binomial theorem, the raw moments of $\hat{\vc{X}}_k$ are then computed from \cref{eq:split_sequence} as
\begin{equation}
    \ex{\hat{\vc{X}}^{\vc{r}}_{k,N}(\vc{x}_0)} = \sum_{\vc{r}'=\vc{0}}^{\vc{r}} \binom{\vc{r}}{\vc{r}'} \hat{\vc{x}}^C_{k}(\vc{x}_0)^{(\vc{r}-\vc{r}')}\ex{\Delta\hat{\vc{W}}^{\vc{r}'}_{k,N}(\vc{x}_0)},
    \label{eq:poly_moments}
\end{equation}
for all $\vc{r}$ such that $0\leq\abs{\vc{r}}\leq N$.

\subsection{Differential algebra-based expansion}\label{sec:plasma:da}

The \gls{plasma} method implements \cref{eq:effective_noise_expansion,eq:effective_noise_moments,eq:poly_moments} in the \gls{da} framework to automate the computation of the partial derivatives and that of the expectations found in \cref{sec:plasma:moments}. Previous work~\cite{Servadio2020} already applied a similar technique in the context of sequential state estimation. In that case, the prior at step $k$ is expressed as a Taylor polynomial function of the posterior at step $k-1$ and of the \gls{pn}. The dynamics, however, are written in integral form without the contribution of the \gls{pn}, and the latter is included with a linear term added to the propagated state. In this work, the dynamics are written in differential form, and the effects of \gls{pn} are integrated numerically together with those of the deterministic accelerations. The main advantage is a more accurate model of the stochastic accelerations, which now have an immediate effect on the trajectory that is being integrated. Indeed, as pointed out by the same authors~\cite{Servadio2020}, their formulation does not account for the effects that the \gls{pn} has on the mean trajectory, as it is assumed a zero-mean random sequence that is linearly added to the predicted state. In contrast, this work does account for its effects on the mean, and complex \gls{pn} models can be introduced by augmenting the state vector to also describe the time evolution of the noise itself. Moreover, the \gls{plasma} method is applicable to any explicit numerical integration method for \glspl{sde}. To derive this algorithm, the differential equation in \cref{eq:sde_def} is firstly rewritten in a suitable form such that it can be evaluated by the selected numerical integrator. This new ``\gls{sde}'' is defined as
\begin{equation}
    \da{\dot{\vc{x}}} = \vc{u}\left(\da{\vc{x}},t\right) + \vc{G}\left(\da{\vc{x}},t\right)\var{\vc{w}},
    \label{eq:pseudo_ode}
\end{equation}
where $\var{\vc{w}}=\left\{\var{w_1},\ldots,\var{w_m}\right\}^T$ is an $m$-dimensional vector of \gls{da} identities representing the $m$ independent noise increments. Note that, to be consistent with \cref{eq:sde_def}, the diffusion term in \cref{eq:pseudo_ode} should be divided by $\dd{t}$. However, it is more convenient to rescale the monomials that depend on $\var{\vc{w}}$ at a later stage as shown in \cref{eq:poly_exp_effective_noise}.

Secondly, a time step $h_k$ at which the effective noise moments are updated must be selected. This can be either a constant time step defined by the user, i.e., $h_k = h$ for any $t_k$, or any step size imposed by the numerical integration method. If an adaptive scheme is used, this step size will generally be time-dependent. In both cases, the objective is to compute the polynomial expansion of the effective noise moments at time $t_k$, i.e. $\da{\Delta\hat{\vc{W}}_{k,N}(\vc{x}_0)}$, given the central part $\hat{\vc{x}}^C_{k-1}(\vc{x}_0)$ at time $t_{k-1}$ (recall that $h_k=t_k-t_{k-1}$). This is achieved by initializing $\hat{\vc{x}}^C_{k-1}(\vc{x}_0)$ as the \gls{da} vector
\begin{equation}
    \da{\hat{\vc{x}}_{k-1}} = \hat{\vc{x}}^C_{k-1}(\vc{x}_0) + \var{\vc{x}}
    \label{eq:rk_ic},
\end{equation}
with $\var{\vc{x}}=\left\{\var{x_1},\ldots,\var{x_n}\right\}^T$ an $n$-dimensional vector of \gls{da} identities that represents small deviations of the state around its central part. Using the specified integration method, \cref{eq:pseudo_ode} is then integrated in $[t_{k-1},t_k]$ with \glspl{ic} given by \cref{eq:rk_ic}. If $h_k$ is selected by the integrator, a single integration step is performed. Several steps may instead be required if $h_k$ is user-defined and larger than the integration time step. In both cases, the output is a Taylor expansion of the state function of the deviations $\var{\vc{x}}$ and the \gls{pn} increments $\var{\vc{w}}$, i.e.
\begin{equation}
    \da{\hat{\vc{x}}_{k}} = \ty{N}{\hat{\vc{x}}_k}{\var{\vc{x}},\var{\vc{w}}}.
    \label{eq:unscaled_state_k}
\end{equation}
To account for the inconsistency between \cref{eq:pseudo_ode,eq:sde_def} mentioned above, the monomials of \cref{eq:unscaled_state_k} that depend on $\var{\vc{w}}$ must now be rescaled by a factor $h_k^{-1}$, where $h_k$ is the time step at which the effective noise moments are updated and is possibly different for each $k$. This scaling is performed by composing the polynomials $\da{\hat{\vc{x}}_{k}}$ in \cref{eq:unscaled_state_k} with the scaling vector $\var{\vc{\kappa}_k}$, which is obtained by stacking the identity vector $\var{\vc{x}}$ with the scaled vector $h_k^{-1}\cdot\var{\vc{w}}$, i.e.
\begin{equation}
    \var{\vc{\kappa}_k}=\left\{\var{\vc{x}}^T, h_k^{-1}\cdot\var{\vc{w}}^T\right\}^T.
    \label{eq:scaling_vector}
\end{equation}
Composing any Taylor polynomial $\ty{N}{\vc{z}}{\var{\vc{x}},\var{\vc{w}}}$ with \cref{eq:scaling_vector} will then result in a new polynomial $\ty{N}{\bar{\vc{z}}}{\var{\vc{x}},\var{\vc{w}}}$ function of the same \gls{da} variables $\left(\var{\vc{x}},\var{\vc{w}}\right)$, but where each monomial is now scaled by $\left(h_k^{-1}\right)^s$, where $s$ is the sum of the exponents with which the components of $\var{\vc{w}}$ appear in the monomial itself. The polynomials $\da{\Delta\hat{\vc{W}}_{k,N}(\vc{x}_0)}$ are obtained by applying this scaling to the current state $\da{\hat{\vc{x}}_{k}}$, and then subtracting its constant part $\hat{\vc{x}}_{k}$ from the result, i.e.
\begin{equation}
    \da{\Delta\hat{\vc{W}}_{k,N}(\vc{x}_0)} = \left(\da{\hat{\vc{x}}_{k}}\circ\var{\vc{\kappa}_k}\right) - \hat{\vc{x}}_{k}.
    \label{eq:poly_exp_effective_noise}
\end{equation}
After this step, the central sequence $\hat{\vc{x}}^C_{k-1}(\vc{x}_0)$ is also updated by setting $\hat{\vc{x}}^C_{k}(\vc{x}_0)=\hat{\vc{x}}_{k}$.

To compute the moments of the effective noise, the polynomials $\da{a_{\vc{r},k}}$ are firstly built for all permutations $\vc{r}=(r_1,\ldots,r_n)$ such that $0\leq\abs{\vc{r}}\leq N$ as
%
% \begin{equation}
%     \da{a_{\vc{r},k}} = \ty{N}{a}{\var{\vc{x}},\var{\vc{w}}} = \prod_{i=1}^n\prod_{j=1}^{r_i} \da{\Delta\hat{W}^{(i)}_{k,N}(\vc{x}_0)},
%     \label{eq:a_expansion}
% \end{equation}
\begin{equation}
    \da{a_{\vc{r},k}} = \ty{N}{a}{\var{\vc{x}},\var{\vc{w}}} = \prod_{i=1}^n \left(\da{\Delta\hat{W}^{(i)}_{k,N}(\vc{x}_0)}\right)^{r_i},
    \label{eq:a_expansion}
\end{equation}
where the coefficients of $\da{a_{\vc{r},k}}$ correspond to $a_{\vc{r},k,N}^{\vc{s},\vc{r}'}\left(\hat{\vc{x}}^C_{k-1}(\vc{x}_0),t_{k-1}\right)$ in \cref{eq:effective_noise_moments}, $\da{\Delta\hat{W}^{(i)}_{k,N}(\vc{x}_0)}$ is the $i^{th}$ component of $\da{\Delta\hat{\vc{W}}_{k,N}(\vc{x}_0)}$ which corresponds to \cref{eq:effective_noise_expansion} $\forall\,i\in[1,n]$, and $\left(\da{\Delta\hat{W}^{(i)}_{k,N}(\vc{x}_0)}\right)^{r_i}$ denotes the $r_i^{th}$ power of $\da{\Delta\hat{W}^{(i)}_{k,N}(\vc{x}_0)}$, i.e. the latter multiplied by itself $r_i$ times. The effective noise moments are then obtained by iterating over the Taylor coefficients of \cref{eq:a_expansion} as
\begin{equation}
    \begin{aligned}
            \ex{\Delta\hat{\vc{W}}^{\vc{r}}_{k,N}(\vc{x}_0)} &= \sum_{\abs{\vc{s}}+\abs{\vc{r}'}=1}^N a_{\vc{r},k}^{\vc{s},\vc{r}'}\cdot\ex{\Delta\vc{W}^{\vc{s}}_{k}\cdot\Delta\hat{\vc{W}}^{\vc{r}'}_{k-1,N}(\vc{x}_0)}\\
            &= \sum_{\abs{\vc{s}}+\abs{\vc{r}'}=1}^N a_{\vc{r},k}^{\vc{s},\vc{r}'}\cdot\ex{\Delta\vc{W}^{\vc{s}}_{k}}\cdot\ex{\Delta\hat{\vc{W}}^{\vc{r}'}_{k-1,N}(\vc{x}_0)},
    \end{aligned}
    \label{eq:effective_noise_moments_da}
\end{equation}
where the last equality holds since all $\Delta W_k^{(j)}$ are uncorrelated by hypothesis. In \cref{eq:effective_noise_moments_da} $a_{\vc{r},k}^{\vc{s},\vc{r}'}$ are the Taylor coefficients corresponding to the deviation $\var{w_1^{s_1}}\cdots\var{w_m^{s_m}}\var{x_1^{r'_1}}\cdots\var{x_n^{r'_n}}$ and $\ex{\Delta\hat{\vc{W}}^{\vc{r}'}_{k-1,N}(\vc{x}_0)}$ are the effective noise moments at the previous step. Finally, $\ex{\Delta\vc{W}^{\vc{s}}_{k}}$ are the $\vc{s}^{th}$-order moments of the multivariate Gaussian distribution with mean $\vc{\mu}_g=\vc{0}_{m\times 1}$ and covariance $\vc{P}_g=h_k\cdot\vc{I}_{m\times m}$. These are easily computed from $\vc{\mu}_g$ and $\vc{P}_g$ using the formulae derived by \citet{Kan2008}.

The conditional moments $\ex{\hat{\vc{X}}^{\vc{r}}_{k,N}(\vc{x}_0)}$ might be then computed at each step using \cref{eq:poly_moments}. The covariance matrix $\hat{\vc{P}}_{k,N}$ of $\hat{\vc{X}}_{k,N}(\vc{x}_0)$ is also obtained using the binomial theorem since its entries are by definition the second-order moments about the mean, i.e.
\begin{equation}
    \begin{aligned}
        \hat{P}_{k,N}^{(i,j)}(\vc{x}_0) &= \ex{\hat{X}^{(i)}_{k,N}(\vc{x}_0)\cdot\hat{X}^{(j)}_{k,N}(\vc{x}_0)} - \ex{\hat{X}^{(i)}_{k,N}(\vc{x}_0)} \cdot \ex{\hat{X}^{(j)}_{k,N}(\vc{x}_0)}\\
        &= \ex{\Delta\hat{W}^{(i)}_{k,N}(\vc{x}_0)\cdot\Delta\hat{W}^{(j)}_{k,N}(\vc{x}_0)} - \ex{\Delta\hat{W}^{(i)}_{k,N}(\vc{x}_0)} \cdot \ex{\Delta\hat{W}^{(j)}_{k,N}(\vc{x}_0)},
    \end{aligned}
    \label{eq:poly_cov}
\end{equation}
where the superscript ``$(i,j)$'' indicates the $i^{th}$ row and $j^{th}$ column of $\hat{\vc{P}}_{k,N}$. To better understand the properties of this algorithm, an application to the discrete Duffing oscillator is presented in \ref{apdx:duffing}.

\subsection{Multifidelity expansion}\label{sec:plasma:multifidelity}

The \gls{plasma} algorithm as described in \cref{sec:plasma:da} might become computationally expensive for high expansion orders $N$ or complex dynamical models $\vc{u}$. A \gls{mf} method is thus developed to reduce the computational cost of this algorithm while maintaining an accuracy similar to its \gls{hf} counterpart. This technique computes the Taylor expansion of the flow of the dynamics relative to a precomputed \gls{hf} trajectory as described hereafter.

Consider a \gls{hf} deterministic dynamical model $\vc{u}_{\mathrm{HF}}$ and its \gls{lf} counterpart $\vc{u}_{\mathrm{LF}}$, where the latter can be derived from $\vc{u}_{\mathrm{HF}}$ by neglecting those terms which are expensive to evaluate and have a minor impact on the propagated trajectory. To start with, the following \gls{ivp} is solved for $t\in[t_0,t_f]$ to obtain the \gls{hf} reference trajectory in the absence of \gls{pn}
\begin{equation}
    \begin{dcases}
        \dot{\vc{x}}(t) = \vc{u}_{\mathrm{HF}}(\vc{x},t)\\
        \vc{x}(t_0) = \ex{\vc{X}_0},
    \end{dcases}
    \label{eq:ivp_hf_ref_bf_sde}
\end{equation}
with $\ex{\vc{X}_0}$ the expected value of $\vc{X}_0$ in \cref{eq:sde_def}, and $t_0,\,t_f$ initial and final integration times, respectively. The so-called dense output provided by the numerical integrator is also stored while solving \cref{eq:ivp_hf_ref_bf_sde} such that continuous models for both the trajectory $\vc{\xi}(t)$ and its first time derivative $\dot{\vc{\xi}}(t)$ are available on exit. The Greek letter $\vc{\xi}$ in the former indicates that these models are piecewise interpolants consistent with the integration scheme used to solve \cref{eq:ivp_hf_ref_bf_sde}. They can thus be evaluated at any time $t\in[t_0,t_f]$ to obtain an approximation of both $\vc{x}(t)$ and $\dot{\vc{x}}(t)$ whose accuracy is controlled by the order and tolerances of the numerical integration method. A characterization of the uncertainty is then maintained through the set of moments $\ex{\Delta\hat{\vc{W}}^{\vc{r}}_{k,N}(\vc{x}_0)}$ centered on $\vc{\xi}(t_k)$, while the central sequence $\hat{\vc{x}}^C_k(\vc{x}_0)$ is now identically zero $\forall\,k$. This is achieved by redefining \cref{eq:pseudo_ode} as
\begin{equation}
    \begin{aligned}
        \var{\dot{\vc{x}}} &= \left[\vc{u}_{\mathrm{LF}}\left(\da{\vc{x}},t\right)-\vc{u}_{\mathrm{LF}}\left(\bar{\vc{x}},t\right)\right] + \vc{G}\left(\da{\vc{x}},t\right)\var{\vc{w}} \\
        &= \var{\vc{u}}_{\mathrm{LF}}\left(\da{\vc{x}},t\right) + \vc{G}\left(\da{\vc{x}},t\right)\var{\vc{w}},
    \end{aligned}
    \label{eq:pseudo_ode_rel}
\end{equation}
where $\da{\vc{x}}$ is defined as $\da{\vc{x}}=\vc{\xi}(t)+\var{\hat{\vc{x}}(t)}$, $\bar{\vc{x}}=\vc{\xi}(t)$ is the constant part of $\da{\vc{x}}$, and $\var{\hat{\vc{x}}(t)}$ is the relative state solution to \cref{eq:pseudo_ode_rel} for $t\in [t_{k-1},t_k]$ subject to \glspl{ic} $\var{\hat{\vc{x}}_{k-1}}=\var{\vc{x}}$. The term $\var{\vc{u}}_{\mathrm{LF}}\left(\da{\vc{x}},t\right)$ is thus the nilpotent part of the \gls{lf} dynamical model evaluated at $\left(\da{\vc{x}},t\right)$. The Taylor expansion in \cref{eq:poly_exp_effective_noise} is then given by
\begin{equation}
    \da{\Delta\hat{\vc{W}}_{k,N}(\vc{x}_0)} = \var{\hat{\vc{x}}_{k}}\circ\var{\vc{\kappa}_k},
    \label{eq:poly_exp_effective_noise_rel}
\end{equation}
while \cref{eq:a_expansion,eq:effective_noise_moments_da} remain the same. The conditional moments $\ex{\hat{\vc{X}}^{\vc{r}}_{k,N}(\vc{x}_0)}$ are finally computed from \cref{eq:poly_moments} where the central sequence $\hat{\vc{x}}^C_{k}(\vc{x}_0)$ is replaced by $\vc{\xi}(t_k)$.

\section{Multifidelity Uncertainty Propagation}\label{sec:multifidelity}

This section presents the second contribution of this work, namely a \gls{mf} method for the efficient propagation of uncertainties in nonlinear dynamical systems. Two variants of this algorithm are discussed: the first one is tailored to deterministic dynamics and is presented in \cref{sec:multifidelity:deterministic}. The second one is instead designed for stochastic dynamics and is discussed in \cref{sec:multifidelity:stochastic}. Both techniques leverage the adaptive \gls{gmm} method presented in \cref{sec:background:gmm} to efficiently propagate the uncertainty in the \glspl{ic} and static model parameters through a \gls{lf} dynamical model. They then differ in the way the \gls{lf} solution is corrected to account for the \gls{hf} dynamics, with the second method leveraging the \gls{plasma} algorithm discussed in \cref{sec:plasma} to account for the effects of \gls{pn}.

\subsection{Deterministic Dynamics}\label{sec:multifidelity:deterministic}

The \gls{mf} method for \gls{up} in deterministic dynamics leverages the adaptive \gls{gmm} algorithm described in \cref{sec:background:gmm} to obtain a \gls{lf} representation of the transformed state \gls{pdf} and corrects this estimate by performing a \gls{pw} propagation of the initial kernels' means~\cite{Fossa2022}. Its workflow is detailed hereafter.

Assume that a \gls{hf} description of the dynamical system of interest is given by the following \gls{ivp}
\begin{equation}
    \begin{dcases}
        \dot{\vc{x}}(t)&=\vc{g}_{\mathrm{HF}}(\vc{x},t)\\
        \vc{x}(t_0)&=\vc{x}_0,
    \end{dcases}
    \label{eq:hf_ivp}
\end{equation}
with $\vc{x}\in\mathbb{R}^n$ the state vector, $\vc{g}_{\mathrm{HF}}(\vc{x},t)$ the first-order \gls{ode} that describe the system's dynamics, and $\vc{x}_0$ the \glspl{ic}. The solution to \cref{eq:hf_ivp} is given by
\begin{equation}
    \vc{\varphi}_{\mathrm{HF}}(t;\vc{x}_0)=\vc{x}_0+\int\limits_{t_0}^{t}\vc{g}_{\mathrm{HF}}(\vc{x},t)\dd{t},
    \label{eq:hf_flow}
\end{equation}
such that $\vc{\varphi}_{\mathrm{HF}}(t_0;\vc{x}_0)=\vc{x}_0$. This \gls{hf} solution is denoted as $\vc{x}_{\mathrm{HF}}(\cdot)=\vc{\varphi}_{\mathrm{HF}}(\cdot;\vc{x}_0)$. At the same time, suppose that a \gls{lf} version of the system \cref{eq:hf_ivp} is also available, with corresponding solution denoted by
\begin{equation}
    \vc{x}_{\mathrm{LF}}(\cdot)=\vc{\varphi}_{\mathrm{LF}}(\cdot;\vc{x}_0).
    \label{eq:lf_model}
\end{equation}
Finally, assume that the initial \gls{pdf} $p_{\vc{X}}(\vc{x}_0)$ is approximated by a known \gls{gmm} as
\begin{equation}
    p_{\vc{X}}(\vc{x}_0)\approx\sum_{l=1}^{L_0}\alpha^{(l)}p_{g}\left(\vc{x}_0;\vc{\mu}^{(l)}_{0},\vc{P}^{(l)}_{0}\right).
    \label{eq:mf_initial_gmm}
\end{equation}

The \gls{mf} method starts by propagating $p_{\vc{X}}(\vc{x}_0)$ to the target time $t_f$ using the adaptive \gls{gmm} algorithm, thus obtaining~\cite{Losacco2024}
\begin{subequations}
    \label{eq:lf_gmm}
    \begin{align}
        p_{\vc{X}}(\vc{x}_0)&\approx\sum_{l=1}^{L}\alpha^{(l)}p_{g}\left(\vc{x}_0;\vc{\mu}^{(l)}_{\mathrm{LF}}(t_0),\vc{P}^{(l)}_{\mathrm{LF}}(t_0)\right) \label{eq:lf_initial_gmm} \\
        p_{\vc{X}}(\vc{x}_f)&\approx\sum_{l=1}^{L}\alpha^{(l)}p_{g}\left(\vc{x}_f;\vc{\mu}^{(l)}_{\mathrm{LF}}(t_f),\vc{P}^{(l)}_{\mathrm{LF}}(t_f)\right) \label{eq:lf_final_gmm},
    \end{align}
\end{subequations}
where the subscript ``LF'' indicates that the \gls{lf} model in \cref{eq:lf_model} was used to propagate each kernel. Moreover, the manifold $\mathcal{M}_f=\left\{\da{\vc{x}_{\mathrm{LF}}^{(l)}(t_f)}\right\}$, which corresponds to \cref{eq:man_out}, is also available at this stage. This manifold contains the polynomial expansions of the state at target time used to estimate $p_{\vc{X}}(\vc{x}_f)$ from $p_{\vc{X}}(\vc{x}_0)$. Given \cref{eq:lf_initial_gmm}, the algorithm then proceeds by propagating each kernel mean $\vc{\mu}^{(l)}_{\mathrm{LF}}(t_0)$ in \gls{hf} dynamics to obtain the \gls{hf} trajectories $\vc{\mu}^{(l)}_{\mathrm{HF}}(t_f)$ as
\begin{equation}
    \vc{\mu}^{(l)}_{\mathrm{HF}}(t_f)=\vc{\varphi}_{\mathrm{HF}}\left(t_f;\vc{\mu}^{(l)}_{\mathrm{LF}}(t_0)\right).
    \label{eq:mf_hf_propagation}
\end{equation}
The difference between the kernel means $\vc{\mu}^{(l)}_{\mathrm{LF}}(t_f)$ in \cref{eq:lf_final_gmm} and the states $\vc{\mu}^{(l)}_{\mathrm{HF}}(t_f)$ in \cref{eq:mf_hf_propagation} might be non-negligible, thus introducing a bias in the mean of the propagated \gls{pdf}. The higher order terms in $\da{\vc{x}_{\mathrm{LF}}^{(l)}(t_f)}$, on the other hand, still provide a good approximation of the relative dynamics in a neighborhood of these trajectories~\cite{Fossa2024}. This is the core idea of the \gls{mf} method, which is similar to the concept of a linearized Kalman filter~\cite{Gelb1974,Zarchan2000}. In the latter, the partials used to update the state estimate are computed around a nominal trajectory that does not incorporate the newly available information. In the \gls{mf} method, the \gls{lf} Taylor expansions are treated as if they where computed around the \gls{hf} trajectories rather than the \gls{lf} ones. In practice, this is achieved by firstly separating $\da{\vc{x}_{\mathrm{LF}}^{(l)}(t_f)}$ into its constant part and nilpotent part as
\begin{equation}
    \da{\vc{x}_{\mathrm{LF}}^{(l)}(t_f)}=\bar{\vc{x}}_{\mathrm{LF}}^{(l)}(t_f)+\ty{2}{\var{\vc{x}}_{\mathrm{LF}}^{(l)}(t_f)}{\var{\vc{x}_0}}
    \label{eq:mf_expansion_splitted}
\end{equation}
and then substituting $\bar{\vc{x}}_{\mathrm{LF}}^{(l)}(t_f)$ in \cref{eq:mf_expansion_splitted} with $\vc{\mu}^{(l)}_{\mathrm{HF}}(t_f)$ to obtain
\begin{equation}
    \da{\vc{x}_{\mathrm{MF}}^{(l)}(t_f)}=\vc{\mu}^{(l)}_{\mathrm{HF}}(t_f)+\ty{2}{\var{\vc{x}}_{\mathrm{LF}}^{(l)}(t_f)}{\var{\vc{x}_0}}.
    \label{eq:mf_expansion_splitted_corrected}
\end{equation}
Finally, a more accurate estimate of the transformed state \gls{pdf} is obtained by evaluating the polynomials $\left\{\da{\vc{x}_{\mathrm{MF}}^{(l)}(t_f)}\right\}$ on \gls{ut} sigma points drawn from the kernels of \cref{eq:lf_initial_gmm} as
\begin{equation}
    p_{\vc{X}}(\vc{x}_f)\approx\sum_{l=1}^{L}\alpha^{(l)}p_{g}\left(\vc{x}_f;\vc{\mu}^{(l)}_{\mathrm{MF}}(t_f),\vc{P}^{(l)}_{\mathrm{MF}}(t_f)\right),
    \label{eq:mf_final_gmm}
\end{equation}
where this last operation is performed as described for \cref{eq:gmm_pdf_prop} and more extensively in \citet{Fossa2024b,Losacco2024}.

\subsection{Stochastic Dynamics}\label{sec:multifidelity:stochastic}

The \gls{mf} method for the propagation of uncertainties in the presence of stochastic accelerations is now presented. This algorithm consists of three main steps: to start with, the uncertainty in the \glspl{ic} and static model parameters is propagated in a \gls{lf} dynamical model using the adaptive \gls{gmm} algorithm discussed in \cref{sec:background:gmm}. Then, for each kernel previously generated, the effects of \gls{pn} are computed using the \gls{plasma} algorithm introduced in \cref{sec:plasma}. The results of the first two steps are finally combined to obtain a \gls{gmm} approximation of the transformed state \gls{pdf} that accounts for all aforementioned sources of uncertainty~\cite{Fossa2023b}. These steps are detailed hereafter and depicted in \cref{fig:mf_plasma_graphical}.

Consider some \glspl{ic} modeled as a multivariate random variable $\vc{X}_0\in\mathbb{R}^n$ with \gls{pdf} $p_{\vc{X}}(\vc{x}_0)$, and assume that a \gls{gmm} approximation of the latter is available at time $t_0$ as
\begin{equation}
    p_{\vc{X}}\left(\vc{x}_0\right)\approx\sum_{l=1}^{L_0}\alpha^{(l)}p_g\left(\vc{x}_0;\vc{\mu}_0^{(l)},\vc{P}_0^{(l)}\right)
    \label{eq:p_x0_gmm_pn}
\end{equation}
with $\alpha^{(l)}$ the mixture weights, $\vc{\mu}_0^{(l)}$ the means, $\vc{P}_0^{(l)}$ the covariance matrices and $L_0$ the initial number of kernels. These \glspl{ic} are firstly propagated in a \gls{lf} dynamical model using the adaptive \gls{gmm} method of \cref{sec:background:gmm}. The outputs thus include a refined \gls{gmm} approximation of $\vc{X}_0$, denoted $f_{\vc{X}}(\vc{x}_0)$, the \gls{gmm} approximation of $\vc{X}$ at target time $t_f$, denoted $f_{\vc{X}}(\vc{x}_f)$, and the two manifolds $\mathcal{M}_0,\mathcal{M}_f^{\mathrm{LF}}$ containing the polynomial expansions of the state at times $t_0,t_f$, respectively. The superscript ``LF'' marks the fact that $\mathcal{M}_f^{\mathrm{LF}}$ was obtained by propagating the \glspl{ic} in a \gls{lf} dynamical model. The resulting \glspl{gmm} are depicted in \cref{fig:mf_plasma_graphical_step1}.

The uncertainty due to stochastic accelerations is then computed by propagating each mean $\vc{\mu}_{0}^{(l)}$ of the kernels in $f_{\vc{X}}(\vc{x}_0)$ using a stochastic dynamical model and one of the two implementations of the \gls{plasma} method described in \cref{sec:plasma:da,sec:plasma:multifidelity}. To start with, the central sequence is initialized as $\vc{x}^C_{0}(\vc{x}_0)=\vc{\mu}_0^{(l)}$ if \cref{eq:pseudo_ode} is used. It is otherwise set to zero if \cref{eq:pseudo_ode_rel} is employed. In both cases, the effective noise moments at $t_0$ are given by
\begin{equation}
    \ex{\Delta\hat{\vc{W}}^{\vc{r}}_{0,N}(\vc{x}_0)} =
    \begin{dcases}
        1 \qquad \abs{\vc{r}} = 0\\
        0 \qquad \abs{\vc{r}} \neq 0
    \end{dcases}
\end{equation}
since the contribution of the initial uncertainty is already taken into account in $\mathcal{M}_f^{\mathrm{LF}}$. The outputs of this step are the $N^{th}$-order approximations of the means $\vc{\mu}^{(l)}_{\mathrm{PN}}(t_f)$ and covariances $\vc{P}^{(l)}_{\mathrm{PN}}(t_f)$, where the subscript ``PN'' marks the fact that these quantities were obtained by propagating deterministic \glspl{ic} in a stochastic dynamical model. In this case, the means $\vc{\mu}^{(l)}_{\mathrm{PN}}(t_f)$ are computed from the set of effective noise moments $\ex{\Delta\hat{\vc{W}}^{\vc{r}}_{k,N}(\vc{x}_0)}$ using \cref{eq:poly_moments}, while the covariances $\vc{P}^{(l)}_{\mathrm{PN}}(t_f)$ are obtained from the same set of moments using \cref{eq:poly_cov}. As a result, at every time step $t_k>t_0$ these means will generally be different from the central sequences $\hat{\vc{x}}^C_k(\vc{x}_0)$ due to the nonlinearities in the dynamics that are captured by the expectations $\ex{\Delta\hat{\vc{W}}^{\vc{r}}_{k,N}(\vc{x}_0)}$. Moreover, this solution provides \gls{hf} reference trajectories for the polynomial expansions in $\mathcal{M}_f^{\mathrm{LF}}$, thus playing the role of \cref{eq:mf_hf_propagation} in \cref{sec:multifidelity:deterministic}. A graphical representation of this step is given in \cref{fig:mf_plasma_graphical_step2}.

Two main assumptions are now introduced. Firstly, for each kernel in $f_{\vc{X}}$, the distance between the mean trajectory estimated by the \gls{plasma} method and the noiseless trajectory output of the adaptive \gls{gmm} is sufficiently small. This hypothesis, which reflects the core idea of the \gls{mf} technique discussed in \cref{sec:multifidelity:deterministic}, justifies the use of $\vc{\mu}^{(l)}_{\mathrm{PN}}(t_f)$ as the new expansion points for the polynomials in $\mathcal{M}_f^{\mathrm{LF}}$. Indeed, if $\vc{\mu}^{(l)}_{\mathrm{PN}}(t)$ and $\vc{\mu}^{(l)}_{\mathrm{LF}}(t)$ are sufficiently close, the relative dynamics in a neighborhood around the former can be approximated by the \gls{lf} Taylor expansions computed around the latter. Secondly, for each kernel in the \gls{gmm} approximation of the state \gls{pdf}, the \gls{pn} is assumed independent from the initial uncertainty, such that the effects of either one do not affect the propagation of the other. This allows to consider their covariances as if they were uncorrelated, and to neglect the effects of the \gls{pn} on the mixture weights $\alpha^{(l)}$. Following these hypotheses, the Taylor expansions in $\mathcal{M}_f^{\mathrm{LF}}$ are firstly re-centered on the \gls{hf} trajectories $\vc{\mu}^{(l)}_{\mathrm{PN}}(t)$ computed above. Since the latter were obtained by propagating the \glspl{ic} in a \gls{hf} dynamical model, this shift compensates for the drift due to the \gls{lf} approximation of the dynamics used to compute $\mathcal{M}_f^{\mathrm{LF}}$. At the same time, $\vc{\mu}^{(l)}_{\mathrm{PN}}(t)$ also accounts for the effects of \gls{pn} on the \gls{pw} trajectories. The result is a new manifold $\mathcal{M}_f^{\mathrm{MF}}$ whose domains are given by
\begin{equation}
    \da{\vc{x}^{(l)}_{\mathrm{MF}}(t_f)} = \vc{\mu}^{(l)}_{\mathrm{PN}}(t_f) + \left\{\da{\vc{x}^{(l)}_{\mathrm{LF}}(t_f)} - \bar{\vc{x}}^{(l)}_{\mathrm{LF}}(t_f)\right\}
    \label{eq:shifted_poly}
\end{equation}
with $\bar{\vc{x}}^{(l)}_{\mathrm{LF}}(t_f)$ the constant part of $\da{\vc{x}^{(l)}_{\mathrm{LF}}(t_f)}$. Similarly to \cref{eq:mf_final_gmm}, a \gls{gmm} approximation of the transformed state \gls{pdf} is finally obtained as
\begin{equation}
    p_{\vc{X}}(\vc{x}_f) \approx \sum_{l=1}^{L} \alpha^{(l)} p_{g}\left(\vc{x}_f;\vc{\mu}_{\mathrm{MF}}^{(l)}(t_f),\vc{P}_{\mathrm{MF}}^{(l)}(t_f)+\vc{P}_{\mathrm{PN}}^{(l)}(t_f)\right)
    \label{eq:inflated_pdf}
\end{equation}
where $\vc{\mu}_{\mathrm{MF}}^{(l)}(t_f)$ and $\vc{P}_{\mathrm{MF}}^{(l)}(t_f)$ are computed by evaluating $\da{\vc{x}^{(l)}_{\mathrm{MF}}(t_f)}$ on \gls{ut} sigma points~\cite{Losacco2024,Fossa2024b} while the weights $\alpha^{(l)}$ are those determined by the adaptive \gls{gmm} method in the first step. This updated approximation is depicted in \cref{fig:mf_plasma_graphical_step3}.

\newcommand{\drawplasmamultifidelityfirststep}{%
    \draw[fill=none,draw=none] (-7,-1) rectangle (18,10);

    % initial uncertainty
    \fill[isaeyellowgreen!20] (0,0) ellipse (2 and 1);
    \draw[draw=isaeyellowgreen, line width=2pt] (0,0) ellipse (2 and 1);
    \fill (0,0) circle (3pt) node[below,xshift=-2.5cm,yshift=0.5cm,text width=2cm,align=right] {initial\\ uncertainty};

    % initial mixture
    \fill[isaeslategray!20] (0,0) ellipse (0.5 and 0.9);
    \draw[draw=isaeslategray, line width=2pt] (0,0) ellipse (0.5 and 0.9);
    \fill[isaeslategray!20] (1,0) ellipse (0.4 and 0.75);
    \draw[draw=isaeslategray, line width=2pt] (1,0) ellipse (0.4 and 0.75);
    \fill[isaeslategray!20] (-1,0) ellipse (0.4 and 0.75);
    \draw[draw=isaeslategray, line width=2pt] (-1,0) ellipse (0.4 and 0.75);

    % initial means
    \fill (0,0) circle (3pt);
    \fill (1,0) circle (3pt);
    \fill (-1,0) circle (3pt);

    % low-fidelity mixture
    \begin{scope}[shift={(4,7.5)}, rotate=-45]

        \begin{scope}
            \clip (-2.1,-2.1) rectangle (2.1,2.1) (0,-1.5) circle (2);
            \fill[isaecoral!20] (0,0) circle (2);
            \draw[draw=isaecoral, line width=2pt] (0,0) circle (2);
        \end{scope}
        \begin{scope}
            \clip (0,0) circle (2);
            \draw[draw=isaecoral, line width=2pt] (0,-1.5) circle (2);
        \end{scope}

        \begin{scope}[rotate=30]
            \fill[isaeslategray!20] (1.5,0) ellipse (0.35 and 0.7);
            \draw[draw=isaeslategray, line width=2pt] (1.5,0) ellipse (0.35 and 0.7);
            \fill (1.5,0) circle (3pt);
        \end{scope}

        \begin{scope}[rotate=-30]
            \fill[isaeslategray!20] (-1.5,0) ellipse (0.35 and 0.7);
            \draw[draw=isaeslategray, line width=2pt] (-1.5,0) ellipse (0.35 and 0.7);
            \fill (-1.5,0) circle (3pt);
        \end{scope}

        \fill[isaeslategray!20] (0,1.25) ellipse (0.9 and 0.6);
        \draw[draw=isaeslategray, line width=2pt] (0,1.25) ellipse (0.9 and 0.6);
        \fill (0,1.25) circle (3pt) node[left,xshift=-1.75cm,yshift=0cm,text width=2.5cm,align=right] {\glsfmtlong{lf} estimate};
    \end{scope}
    \draw[->, line width=1.5pt] (0.5,1.5) -- (4,7.5);
}

\newcommand{\drawplasmamultifidelitysecondstep}{%
    \begin{scope}[shift={(14,3)}, rotate=-45]
        \begin{scope}[rotate=40]
            \fill[isaeaquamarine!20] (2.25,-0.25) ellipse (0.3 and 0.6);
            \draw[draw=isaeaquamarine, line width=2pt] (2.25,-0.25) ellipse (0.3 and 0.6);
            \fill (2.25,-0.25) circle (3pt);
        \end{scope}

        \begin{scope}[rotate=-40]
            \fill[isaeaquamarine!20] (-2.25,-0.25) ellipse (0.3 and 0.6);
            \draw[draw=isaeaquamarine, line width=2pt] (-2.25,-0.25) ellipse (0.3 and 0.6);
            \fill (-2.25,-0.25) circle (3pt);
        \end{scope}

        \fill[isaeaquamarine!20] (0,2) ellipse (0.8 and 0.4);
        \draw[draw=isaeaquamarine, line width=2pt] (0,2) ellipse (0.8 and 0.4);
        \fill (0,2) circle (3pt) node[left,xshift=1.75cm,yshift=1.5cm,text width=2cm] {\glsfmtlong{hf} moments};
    \end{scope}
    \draw[->, line width=1.5pt] (2.5,0.5) -- (14.5,3.5);
}

\newcommand{\drawplasmamultifidelitythirdstep}{%
    \begin{scope}[shift={(14,3)}, rotate=-45]

        \begin{scope}
            \clip (-3.1,-3.1) rectangle (3.1,3.1) (0,-2) circle (3);
            \fill[isaecoral!20] (0,0) circle (3);
            \draw[draw=isaecoral, line width=2pt] (0,0) circle (3);
        \end{scope}
        \begin{scope}
            \clip (0,0) circle (3);
            \draw[draw=isaecoral, line width=2pt] (0,-2) circle (3);
        \end{scope}

        \begin{scope}[rotate=40]
            \fill[isaeslategray!20] (2.25,-0.25) ellipse (0.55 and 1.1);
            \draw[draw=isaeslategray, line width=2pt] (2.25,-0.25) ellipse (0.55 and 1.1);
            \fill (2.25,-0.25) circle (3pt);
        \end{scope}

        \begin{scope}[rotate=-40]
            \fill[isaeslategray!20] (-2.25,-0.25) ellipse (0.55 and 1.1);
            \draw[draw=isaeslategray, line width=2pt] (-2.25,-0.25) ellipse (0.55 and 1.1);
            \fill (-2.25,-0.25) circle (3pt);
        \end{scope}

        \fill[isaeslategray!20] (0,2) ellipse (1.2 and 0.85);
        \draw[draw=isaeslategray, line width=2pt] (0,2) ellipse (1.2 and 0.85);
        \fill (0,2) circle (3pt) node[below,xshift=1.25cm,yshift=-2.25cm,text width=2.1cm] {\glsfmtlong{mf} estimate};

    \end{scope}
    \draw[->, line width=1.5pt] (2.5,0.5) -- (14.25,3.25);
    \draw[->, line width=1.5pt] (6.5,7.25) -- (11.5,5.5);
}

\begin{figure}[p]
    \centering
    \tikzset{>=stealth}
    \subcaptionbox{Propagation of the initial uncertainty using the adaptive \glsfmtshort{gmm} method\label{fig:mf_plasma_graphical_step1}}{%
        \begin{tikzpicture}[scale=0.57]
            \drawplasmamultifidelityfirststep
            \node[text width=4cm,align=right] at (-2,4.5) {\phantom{low-fi}Adaptive \glsfmtshort{gmm} \glsfmtlong{lf} dynamics};
        \end{tikzpicture}
    }
    \vspace*{\fill}
    \subcaptionbox{Propagation of the effects of stochastic accelerations using the \glsfmtshort{plasma} method\label{fig:mf_plasma_graphical_step2}}{%
        \begin{tikzpicture}[scale=0.57]
            \drawplasmamultifidelityfirststep
            \drawplasmamultifidelitysecondstep
            \node[text width=5cm] at (11.5,0.5) {\Glsfmtshort{plasma} method\\ \glsfmtlong{hf} dynamics};
        \end{tikzpicture}
    }
    \vspace*{\fill}
    \subcaptionbox{Correction of the polynomial expansions and computation of the kernels' statistics\label{fig:mf_plasma_graphical_step3}}{%
        \begin{tikzpicture}[scale=0.57]
            \drawplasmamultifidelityfirststep
            \drawplasmamultifidelitythirdstep
            \node[text width=4cm] at (12,7) {Polynomial correction};
        \end{tikzpicture}
    }
    \caption{Graphical representation of the \glsfmtlong{mf} \glsfmtlong{up} method with stochastic accelerations}
    \label{fig:mf_plasma_graphical}
\end{figure}

\section{Numerical Applications}\label{sec:applications}

This section demonstrates the performance of the proposed \gls{mf} \gls{up} techniques applied to the problem of orbit \gls{up}. Three test cases are presented: \cref{sec:applications:deterministic} studies the performance of the algorithm described in \cref{sec:multifidelity:deterministic} for different Earth orbit regimes, while \cref{sec:applications:keplerian_motion} and \cref{sec:applications:low_thrust} applies the method described in \cref{sec:multifidelity:stochastic} to the perturbed Keplerian dynamics and a low-thrust orbit raising maneuver, respectively. The algorithms were coded in Java and make use of the \gls{cnes}'s \gls{da} library \gls{pace}. All simulations were run on a cluster equipped with Skylake Intel\textsuperscript{\textregistered} Xeon\textsuperscript{\textregistered} Gold 6126 @ \qty{2.6}{\giga\hertz} and \qty{96}{\giga\byte} of \gls{ram} running Red Hat\textsuperscript{\textregistered} Enterprise Linux\textsuperscript{\textregistered} 7.8.

\subsection{Deterministic orbital dynamics}\label{sec:applications:deterministic}

The \gls{mf} method discussed in \cref{sec:multifidelity:deterministic} is firstly applied to the propagation of uncertainties in three orbital regimes: a \gls{leo}, a \gls{heo} and a \gls{meo}. The first two cases are adapted from \citet{DeMars2013}, while the third one reflects the parameters of the Galileo constellation. The nominal \glspl{ic} are expressed in \gls{coe} and summarized in \cref{tab:det_ics} together with the corresponding $1\sigma$ uncertainties in position and velocity.

\begin{table}[!ht]
    \caption{Nominal \glsfmtlongpl{ic} and $1\sigma$ uncertainties}
    \label{tab:det_ics}
    \centering
    \sisetup{round-mode=places,exponent-mode=input}
    \begin{tabular}{r S[round-precision=3,table-format=5.3] S[exponent-mode=scientific,round-precision=0,table-format=1e-1] S[round-precision=1,table-format=2.1] *3{S[round-precision=1,table-format=1.1]}}
    \toprule
    Regime & {$a$, \unit{\km}} & {$e$} & {$i$, \unit{\deg}} & {$\Omega$, \unit{\deg}} & {$\omega$, \unit{\deg}} & {$M$, \unit{\deg}} \\ \midrule
    \gls{heo} & 35000.0 & 2e-1 & 0.0 & 0.0 & 0.0 & 0.0 \\
    \gls{leo} & 6678.0 & 1e-2 & 0.0 & 0.0 & 0.0 & 0.0 \\
    \gls{meo} & 29600.135 & 1e-4 & 56.0 & 0.0 & 0.0 & 0.0 \\ \midrule
    & {$\sigma_x$, \unit{\km}} & {$\sigma_y$, \unit{\km}} & {$\sigma_z$, \unit{\km}} & {$\sigma_{vx}$, \unit{\m\per\s}} & {$\sigma_{vy}$, \unit{\m\per\s}} & {$\sigma_{vz}$, \unit{\m\per\s}} \\ \midrule
    \gls{heo} & 1.0 & 1.0 & 0.0 & 1.0 & 1.0 & 0.0 \\
    \gls{leo} & 1.3 & 0.5 & 0.0 & 2.5 & 5.0 & 0.0 \\
    \gls{meo} & 0.5 & 1.0 & 1.0 & 0.5 & 0.5 & 0.5 \\
    \bottomrule
    \end{tabular}
\end{table}

The initial \gls{pdf} in Cartesian space is assumed Gaussian with mean $\vb*{\mu_0}$ and diagonal covariance matrix $\vb*{P}_0=\mathrm{diag}\left(\sigma_i^2\right)$, where $\vb*{\mu_0}$ is obtained by transforming the \glspl{coe} in \cref{tab:det_ics} to Cartesian parameters. The \gls{gmm} representation of the initial \gls{pdf} input to the \gls{mf} \gls{up} method is thus composed of a single kernel with unitary weight, mean $\vb*{\mu}_0$, and covariance matrix $\vb*{P}_0$.

The \gls{hf} dynamics selected for the \gls{pw} propagation of the initial kernels' means and for the reference \gls{mc} simulation is the perturbed Keplerian motion that accounts for the following forces: Earth non-uniform gravity field with $8 \times 8$ spherical harmonics potential~\cite{Holmes2002}, Sun and Moon third-body attraction with positions of the celestial bodies obtained from \gls{jpl} \gls{de}440 ephemerides~\cite{Park2021}, isotropic \gls{srp} force including Earth's umbra and penumbra transitions~\cite{Hubaux2012} and, for the \gls{leo} case, isotropic drag force with modified Harris-Priester static atmosphere model~\cite{Harris1962}. The numerical parameters for these perturbations are summarized in \cref{tab:forces_params}.

\begin{table}[!ht]
    \caption{Force models parameters}
    \label{tab:forces_params}
    \centering
    \begin{tabular}{ccl}
        \toprule
        Parameter & Value & Description\\ \midrule
        $\mu_{\oplus}$ & \qty{398600.4355}{\cubic\km\per\square\s} & Earth standard gravitational parameter\\
        $m$ & \qty{500.0}{\kg} & satellite mass\\
        $A_n$ & \qty{1.0}{\square\m} & satellite cross section area\\
        $C_D$ & 2.0 & atmospheric drag coefficient\\
        $C_R$ & 1.5 & \gls{srp} reflection coefficient\\
        \bottomrule
    \end{tabular}
\end{table}

Instead, the \gls{lf} dynamics used by the adaptive \gls{gmm} method is the analytical \gls{sgp4} model~\cite{Vallado2006a,Vallado2013}. As the \gls{lf} dynamics needs to be evaluated in \gls{da}, this work uses the implementation of \gls{sgp4} available in Orekit~\cite{Maisonobe2010}, namely the \texttt{FieldTLEPropagator} class\footnote{\url{https://www.orekit.org/site-orekit-12.2.1/apidocs/org/orekit/propagation/analytical/tle/FieldTLEPropagator.html}}, which is already compatible with \gls{da} operations. However, in its original form \gls{sgp4} requires the initial state to be represented in \gls{tle} format, and returns an osculating state in Cartesian coordinates. It also expresses these states in a non-standard frame, the \gls{teme} frame, which should be only used in conjunction with the specific \gls{sgp4} implementation~\cite{Vallado2013}. To overcome these limitations, a wrapper was built around the Orekit implementation to allow the use of osculating states expressed in a standard inertial reference frame as both input and output of the propagation~\cite{Fossa2024b}. As the conversion from osculating states to \glspl{tle} involves an iterative process, only the constant part of the polynomial expansion of the state is considered to assess the convergence of the algorithm. Moreover, as the internal representation of the state uses mean \glspl{oe}, the correction described by \cref{eq:mf_expansion_splitted_corrected} could be performed in either osculating \glspl{oe} or mean ones. The first approach is denoted as ``\Glsfmtshort{mf} osc'', while the second is named ``\Glsfmtshort{mf} mean''. If not explicitly stated, the correction is performed in osculating elements space.

All simulations consider an initial epoch set on \DTMdisplaydate{2021}{1}{1}{-1} at \DTMdisplaytime{0}{0}{0} \gls{utc} and a time span of two revolutions of the corresponding nominal orbit. Orbit states are expressed in \gls{gcrf}. The adaptive \gls{gmm} algorithm uses the three-components splitting library with a penalty factor for the components' variances (denoted $\lambda$ in \citet{DeMars2013}) equal to $\num{e-3}$, an uncertainty scaling factor $\zeta=3$, a maximum number of splits $N_{max}=20$, and a minimum components weight $\alpha_{min}=\num{e-6}$. The nonlinearity threshold $\varepsilon_\nu$ is instead tuned depending on the orbital regime as described below.

\begin{figure}[!ht]
    \centering
    \subcaptionbox{$x$-$y$ plane projection\label{fig:heo_xy_contour}}[0.49\textwidth]{\includegraphics[width=0.49\textwidth]{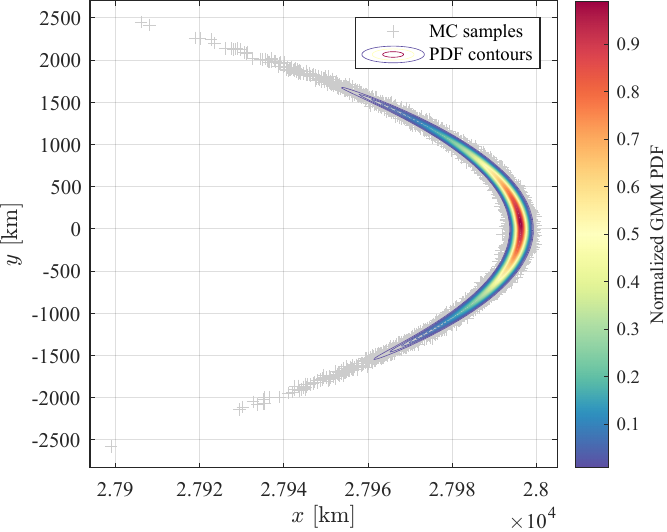}}%
    \hfill
    \subcaptionbox{$v_x$-$v_y$ plane projection\label{fig:heo_vxvy_contour}}[0.49\textwidth]{\includegraphics[width=0.49\textwidth]{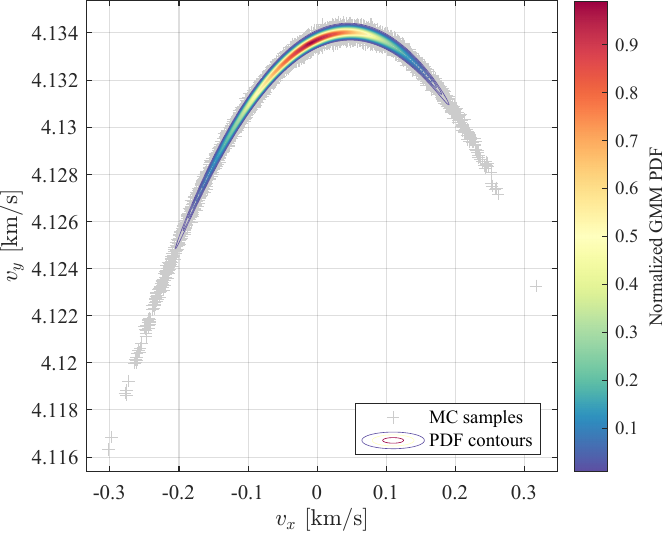}}
    \caption{Propagated \glsfmtshort{mc} samples and mixture \glsfmtshort{pdf} in Cartesian parameters for \glsfmtshort{heo} case}
    \label{fig:heo_contours}
\end{figure}

\begin{figure}[!ht]
    \centering
    \subcaptionbox{$x$-$y$ plane projection\label{fig:leo_xy_contour}}[0.49\textwidth]{\includegraphics[width=0.49\textwidth]{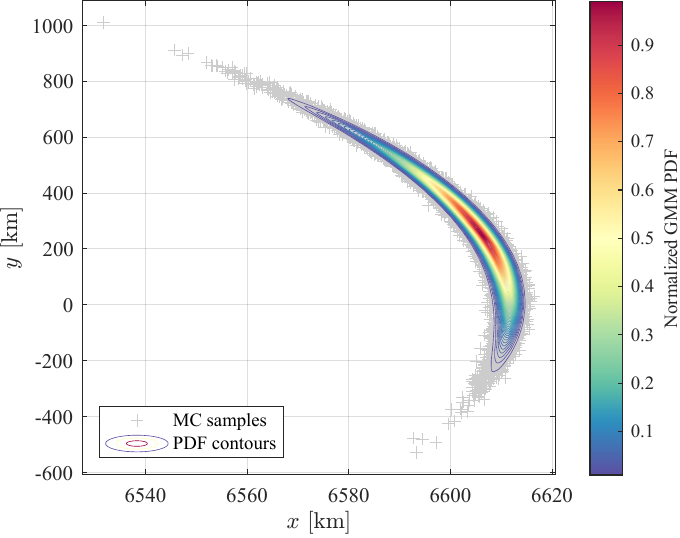}}%
    \hfill
    \subcaptionbox{$v_x$-$v_y$ plane projection\label{fig:leo_vxvy_contour}}[0.49\textwidth]{\includegraphics[width=0.49\textwidth]{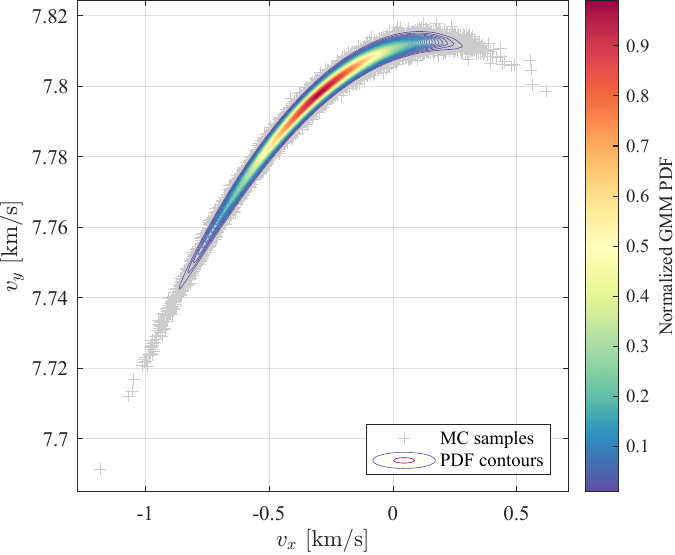}}
    \caption{Propagated \glsfmtshort{mc} samples and mixture \glsfmtshort{pdf} in Cartesian parameters for \glsfmtshort{leo} case}
    \label{fig:leo_contours}
\end{figure}

\begin{figure}[!ht]
    \centering
    \subcaptionbox{$x$-$z$ plane projection\label{fig:meo_xz_contour}}[0.49\textwidth]{\includegraphics[width=0.49\textwidth]{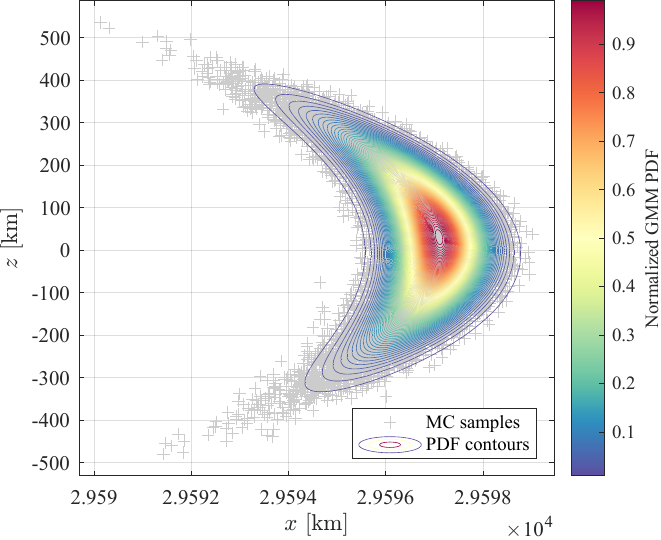}}%
    \hfill
    \subcaptionbox{$v_x$-$v_z$ plane projection\label{fig:meo_vxvz_contour}}[0.49\textwidth]{\includegraphics[width=0.49\textwidth]{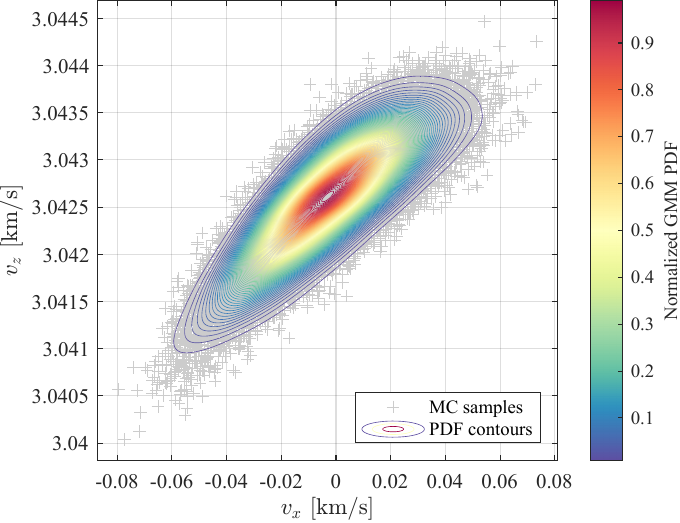}}
    \caption{Propagated \glsfmtshort{mc} samples and mixture \glsfmtshort{pdf} in Cartesian parameters for \glsfmtshort{meo} case}
    \label{fig:meo_contours}
\end{figure}

A first analysis is performed by propagating the \gls{hf} trajectories using the Cowell's method~\cite{Curtis2021}. The results are shown in \cref{fig:heo_contours,fig:leo_contours,fig:meo_contours}, where the iso-probability lines of the transformed state \gls{pdf} are plotted on top of the particle ensemble resulting from a \gls{hf} \gls{mc} simulation with \num{5e4} samples.

To start with, outputs of the \gls{heo} case obtained with $\varepsilon_\nu=0.012$ are displayed in \cref{fig:heo_contours}. \Cref{fig:heo_xy_contour} shows the propagated uncertainty in position projected onto the $x-y$ plane, while \cref{fig:heo_vxvy_contour} represents the uncertainty in velocity also projected onto the same plane. Values of the estimated \gls{pdf} are normalized with respect to their maximum located in correspondence of the expected state. In total, \num{2187} kernels are generated by the adaptive algorithm to accurately capture the curvature of the transformed state \gls{pdf}. The accuracy of the \gls{mf} method is then demonstrated by the close match between the shape of the \gls{gmm} approximation and the distribution of the reference \gls{mc} samples.

The second scenario is the equatorial \gls{leo} for which the estimated \gls{pdf} is displayed in \cref{fig:leo_contours}. \Cref{fig:leo_xy_contour} shows a projection onto the $x-y$ plane of the position uncertainty, while \cref{fig:leo_vxvy_contour} depicts a projection onto the same plane of the uncertainty in velocity. As before, the effects of nonlinearities are visible in the curvature of the initially Gaussian \gls{mc} particles ensemble, which is closely followed by the contour lines of the estimated \gls{pdf}. Despite a higher nonlinearity threshold, now equal to $\varepsilon_\nu=0.024$, and a lower eccentricity with respect to the previous scenario, the same number of kernels is needed to approximate the transformed \gls{pdf}. This behavior is mainly due to the larger uncertainty in velocity with respect to the \gls{heo} case, but is also contributed by the lower orbit altitude that exacerbates the effects of gravitational perturbations and introduces those of the atmospheric drag.

The last test case is the quasi-circular, inclined orbit of the Galileo constellation. The propagated uncertainties are shown in \cref{fig:meo_contours} together with the reference \gls{mc} samples. \Cref{fig:meo_xz_contour} represents the uncertainty in position projected onto the $x-z$ plane, while \cref{fig:meo_vxvz_contour} shows the corresponding uncertainty in velocity. In this case, setting $\varepsilon_\nu=0.012$ leads to \num{81} components needed to approximate the propagated \gls{pdf}. The great reduction in the number of kernels with respect to the previous cases is mainly due to the lower eccentricity (compared to the \gls{heo} case), the longer orbital period, and the absence of atmospheric drag (compared to the \gls{leo} case). The lower number of components is thus a hint that nonlinearities are weaker in this scenario, as demonstrated by the shape of the contour lines in \cref{fig:meo_vxvz_contour}. The latter are in fact closer to a Gaussian ellipsoid than any other projection shown in \cref{fig:heo_contours,fig:leo_contours}. This behavior highlights the advantages of an adaptive splitting scheme with respect to more conventional \glspl{gmm} with a fixed number of components: the final number of Gaussian kernels is the minimum required to satisfy the imposed nonlinearity threshold and new splits only occur if $\nu>\varepsilon_\nu$. Moreover, no \emph{a priori} information on the nonlinearity of the dynamics is needed to generate a Gaussian mixture before the propagation takes place.

Two performance indices are now computed to assess the accuracy of the \gls{gmm} methods compared to that of the \gls{mc} simulation. The first index is the largest component in the relative error between the estimated mean $\vc{\mu}_\mathrm{GMM}$ and the ensemble mean $\vc{\mu}_\mathrm{MC}$, i.e.
\begin{equation}
    \varepsilon_{\vc{\mu}} = \norm{\frac{\vc{\mu}_\mathrm{GMM}-\vc{\mu}_\mathrm{MC}}{\vc{\mu}_\mathrm{MC}}}_\infty,
    \label{eq:mean_err_idx}
\end{equation}
where the subtraction and division are performed element-wise. The second one is instead the relative error in the largest eigenvalue of the estimated covariance $\vc{P}_\mathrm{GMM}$ and the sample covariance $\vc{P}_\mathrm{MC}$, i.e.
\begin{equation}
    \varepsilon_{\vc{\lambda}} = \abs{\frac{\lambda_{\max}(\vc{P}_\mathrm{GMM})-\lambda_{\max}(\vc{P}_\mathrm{MC})}{\lambda_{\max}(\vc{P}_\mathrm{MC})}}.
    \label{eq:cov_err_idx}
\end{equation}
These indices are reported in \cref{tab:err_idx}, where each block of two subsequent columns corresponds to a different orbital regime. The three rows refer instead to a \gls{lf} propagation, where no correction is applied, and to the two \gls{mf} solutions, which are obtained by applying \cref{eq:mf_expansion_splitted_corrected} in osculating and mean elements, respectively.

\begin{table}[ht]
    \caption{Relative errors for the three regimes}
    \label{tab:err_idx}
    \centering
    \small
    \sisetup{round-mode=places,round-precision=3}
    \begin{tabular}{l *6{S[scientific-notation=true,table-format=1.3e-1]}}
        \toprule
        & \multicolumn{2}{c}{\Glsfmtlong{heo}} & \multicolumn{2}{c}{\Glsfmtlong{leo}} & \multicolumn{2}{c}{\Glsfmtlong{meo}} \\ \midrule
        & $\varepsilon_{\vc{\mu}}$ & $\varepsilon_{\vc{\lambda}}$ & $\varepsilon_{\vc{\mu}}$ & $\varepsilon_{\vc{\lambda}}$ & $\varepsilon_{\vc{\mu}}$ & $\varepsilon_{\vc{\lambda}}$ \\ \midrule
        \glsfmtshort{lf} & 0.3505439450645066 & 5.9635609673003645e-5 & 0.012559820082498058 & 0.07178863297049216 & 0.00012941186244583324 & 3.622012109113002e-6 \\
        \glsfmtshort{mf} osc & 0.013949009146991848 & 5.9825440073983335e-5 & 0.001165587347724993 & 0.07727151458075145 & 1.5212801637560385e-5 & 3.626263988231959e-6 \\
        \glsfmtshort{mf} mean & 0.013949923572607196 & 5.982510075026718e-5 & 0.0011653209588304088 & 0.07725709124576786 & 1.515918234769949e-5 & 3.5903836369694257e-6 \\ \bottomrule
    \end{tabular}
\end{table}

Few trends are immediately identified by inspecting the results in \cref{tab:err_idx}. Firstly, both \gls{mf} methods yields very similar results in all three regimes. They also reduce the error in the estimated mean $\varepsilon_{\vc{\mu}}$ by one order of magnitude with respect to the \gls{lf} solution. Conversely, errors in the largest eigenvalue of the covariance $\varepsilon_{\vc{\lambda}}$ are almost constant across the three rows. This is also expected when the mean of each Gaussian kernel does not deviate much from the constant part of the associated polynomial. Indeed, this deviation is the only difference between the three solutions, as the nilpotent part of the polynomial used to compute the kernels' statistics is the same across the three methods. Finally, the best results are obtained for the \gls{meo} case, where the propagated \gls{pdf} is closer to a Gaussian distribution as demonstrated by the lower number of components generated in this scenario. Overall, these results confirm the effectiveness of the \gls{mf} method, which was able to accurately capture the mean and the covariance of the propagated \gls{pdf} with an error of less than \SI{2}{\percent} and \SI{8}{\percent}, respectively, compared to the reference \gls{mc} simulation.

Since the nonlinearity of any dynamical system is not an intrinsic attribute of the system itself, but it is rather dependent upon the mathematical description of the physical system~\cite{Junkins2004}, it is of interest to analyze the performance of the \gls{up} method for different representations of the orbit state. In the following, four different coordinates sets are considered: the aforementioned Cartesian parameters, \glspl{ee}~\cite{Broucke1972}, \glspl{mee}~\cite{Walker1985}, and \glspl{aee}~\cite{Horwood2011}. In addition, both true and mean longitudes $L,\lambda$ are used as fast variable for the last three sets. The \gls{hf} dynamics is then propagated directly in osculating \glspl{oe} using the Gauss's variational equations~\cite{Battin1999}, while the \gls{lf} estimate is obtained by transforming the output of the \gls{sgp4} propagator to the desired set.

\begin{table}[ht]
    \caption{Limit \glsfmtshort{nli} $\nu_{s}$, nonlinearity threshold $\varepsilon_\nu$ and number of \glsfmtshort{gmm} kernels for the three regimes}
    \label{tab:nu_single_nb_kernels}
    \begin{adjustwidth}{-1cm}{-1cm}
    \centering
    \small
    \sisetup{round-mode=places,round-precision=3}
    \begin{tabular}{rc S[scientific-notation=true,table-format=1.3e-2] S[scientific-notation=false,table-format=1.3] S[scientific-notation=false,table-format=4,round-precision=0] S[scientific-notation=true,table-format=1.3e-2] S[scientific-notation=false,table-format=1.3] S[scientific-notation=false,table-format=4,round-precision=0] S[scientific-notation=true,table-format=1.3e-2] S[scientific-notation=false,table-format=1.3] S[scientific-notation=false,table-format=4,round-precision=0]}
        \toprule
        & & \multicolumn{3}{c}{\Glsfmtlong{heo}} & \multicolumn{3}{c}{\Glsfmtlong{leo}} & \multicolumn{3}{c}{\Glsfmtlong{meo}} \\ \midrule
        \Glsfmtshortpl{oe} & {Fast angle} & $\nu_s$ & $\varepsilon_\nu$ & $N$ & $\nu_s$ & $\varepsilon_\nu$ & $N$ & $\nu_s$ & $\varepsilon_\nu$ & $N$ \\ \midrule
        Cartesian & $-$ & 0.065204301727872 & 0.012 & 2187 & 0.0997430092099781 & 0.024 & 2187 & 0.0238560669651193 & 0.012 & 81 \\
        \glsfmtshortpl{ee} & $L$ & 0.0066548141136642 & 0.002 & 81 & 0.0102397905910868 & 0.004 & 81 & 0.00180751551947851 & 0.002 & 1 \\
        \glsfmtshortpl{mee} & $L$ & 0.0103393929185414 & 0.002 & 2187 & 0.0178889143449065 & 0.004 & 729 & 0.00350081558730028 & 0.002 & 9 \\
        \glsfmtshortpl{aee} & $L$ & 0.00680248266455395 & 0.002 & 81 & 0.0149456676849809 & 0.004 & 81 & 0.00296886526411071 & 0.002 & 9 \\
        \glsfmtshortpl{ee} & $\lambda$ & 0.00551441923542427 & 0.002 & 27 & 0.0113992266693183 & 0.004 & 63 & 0.00210418996456917 & 0.002 & 3 \\
        \glsfmtshortpl{mee} & $\lambda$ & 0.00552283004359903 & 0.002 & 243 & 0.0136099313961443 & 0.004 & 81 & 0.00279669810226131 & 0.002 & 9 \\
        \glsfmtshortpl{aee} & $\lambda$ & 3.83131414352759e-06 & 0.002 & 1 & 0.000609131727108636 & 0.004 & 1 & 7.5662554270136e-06 & 0.002 & 1 \\ \bottomrule
    \end{tabular}
    \end{adjustwidth}
\end{table}

\begin{figure}[!ht]
    \centering
    \includegraphics[width=0.4\textwidth]{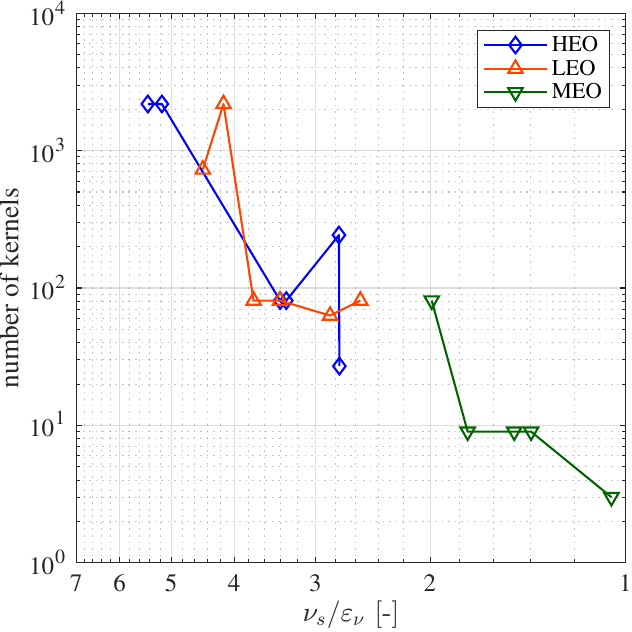}
    \caption{Number of \glsfmtshort{gmm} kernels as a function of normalized limit \glsfmtshort{nli} $\nu_{s}/\varepsilon_\nu$}
    \label{fig:nb_kernels_vs_nu_single}
\end{figure}

\Cref{tab:nu_single_nb_kernels} reports the limit \gls{nli} $\nu_s$, the nonlinearity threshold $\varepsilon_\nu$ and the final number of Gaussian kernels $N$ for the \gls{heo}, \gls{leo} and \gls{meo} cases, respectively. The parameter $\nu_s$ corresponds to the \gls{nli} of the single final domain $\left[\vb*{x}(t_f)\right]$ when mapped with an infinite threshold $\varepsilon_\nu=\infty$ and is a direct measure of nonlinearity of the different state representations. For a fixed threshold $\varepsilon_\nu$, $\nu_s$ is thus expected to be directly correlated with the final number of components $N$. This fact is demonstrated in \cref{fig:nb_kernels_vs_nu_single} in which $N$ is plotted as a function of the normalized limit \gls{nli} $\nu_s/\varepsilon_\nu$ (only the cases for which $\nu_s > \varepsilon_\nu$ and thus $N>1$ are shown in the picture). For each orbital regime, this figure shows in fact a positive correlation between $\nu_s/\varepsilon_\nu$ and the corresponding number of components, as expected by construction of the splitting algorithm. In each scenario, it is clear that Cartesian parameters are the most nonlinear since they lead to a larger number of kernels despite the higher thresholds selected in this case. Moreover, the use of \glspl{aee} and mean longitude $\lambda$ results in a quasi-linear transformation that is correctly captured by a single polynomial in all three scenarios.

Unfortunately, coordinates sets can rarely be chosen freely. Many applications require in fact the \gls{pdf} of the object's state expressed in Cartesian parameters. One example is that of \gls{poc} computations, which require a probabilistic description of the positions and velocities of the two objects at \gls{tca}. Adding a coordinate transformation at the end of the \gls{lf} propagation will be detrimental in this case, since the \gls{nli} will still be computed on Cartesian parameters, thus losing all the benefits of the alternative state parametrization. However, this is only true for analytical propagators which require a single function evaluation to compute the final state. Conversely, for semianalytical and numerical propagators splits can occur at any time during propagation. Since each split will then correspond to new trajectories to be propagated, it is of great interest to delay them as much as possible and propagate a single expansion for the entire time span. This behavior could be achieved using \glspl{aee} for the propagation itself, followed by a coordinate transformation of the single patch at final time in which all splits will occur.

\begin{figure}[!ht]
    \centering
    \includegraphics[width=0.5\textwidth]{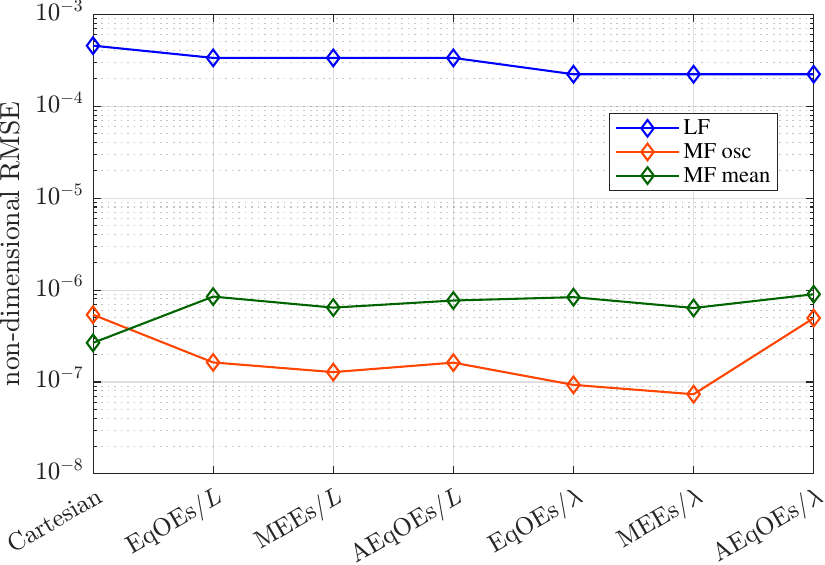}
    \caption{Non-dimensional \glsfmtshort{rmse} for \glsfmtlong{heo} case}
    \label{fig:rmse_heo}
\end{figure}

Quantitative analyses of the performance of the \gls{mf} method were then conducted for each orbital regime. An extensive discussion of these results is provided for the \gls{heo} scenario, while for conciseness the \gls{leo} and \gls{meo} cases report only the main differences with respect to the former. To start with, the accuracy of this technique is assessed by computing a non-dimensional, average \gls{rmse} between the samples output of the \gls{mc} simulation and those obtained by evaluating the polynomial manifold at final time. This \gls{rmse} is computed as
\begin{equation}
    \mathrm{RMSE}_\xi = \sqrt{\frac{1}{N_{\mathrm{MC}}\cdot n} \sum_{i=1}^{N_{\mathrm{MC}}}\sum_{j=1}^n \left(\hat{\xi}_{i,j} - \xi_{i,j}\right)^2}
    \label{eq:normalized_rmse}
\end{equation}
with $\hat{\xi}_{i,j}$ the $j^{th}$ component of the $i^{th}$ actual sample, $\xi_{i,j}$ the corresponding expected value, $N_{\mathrm{MC}}$ the number of \gls{mc} samples, and $n$ the dimension of the state space. Samples denoted with $\vb*{\xi}$ are normalized with respect to some reference quantities such that all entries are of order one. In this case, the unit of length is selected equal to the radius of pericenter, which is thus equal to one in non-dimensional units. The unit of time is instead selected such that the non-dimensional standard gravitational parameter $\mu$ is also equal to one. The \gls{rmse} is then computed for the \gls{lf} solution and the two \gls{mf} estimates with correction in osculating and mean elements respectively. The results for the \gls{heo} regime are depicted in \cref{fig:rmse_heo}. In this plot, the $x$ axis represents the \glspl{oe} and fast angle used to parameterize the orbit state, while the different colors mark the different fidelity models. The blue lines correspond to the \gls{lf} solution, while the orange and green lines correspond to the two \gls{mf} solutions.

By looking at \cref{fig:rmse_heo}, it is immediately seen that the \gls{mf} solutions are between three and four orders of magnitude more accurate than the \gls{lf} one. In particular, the latter exhibits a \gls{rmse} comprised between \num{2e-4} and \num{5e-4}, with the propagation in Cartesian parameters performing worst. Concerning the propagation in osculating elements, the error is mainly affected by the choice of the fast angle, with the mean longitude $\lambda$ guaranteeing better results compared to the true longitude $L$. The \gls{rmse} is otherwise almost independent on the selected \glspl{oe} set. Moving on to the two \gls{mf} solutions, a correction in osculating space provides a better accuracy in most cases, with the only exception of Cartesian parameters. For a correction in mean elements, the \gls{rmse} is in fact comprised between \num{6e-7} and \num{e-6} when propagating in either \glspl{ee}, \glspl{mee}, or \glspl{aee}, while it drops to \num{3e-7} for a propagation in Cartesian parameters. The opposite trend is instead observed when \cref{eq:mf_expansion_splitted_corrected} is applied in osculating space. In this case, the \gls{rmse} is about \num{5e-7} for Cartesian parameters, drops below \num{2e-7} for most \glspl{oe} sets, but it then exhibits a second peak when propagating in \glspl{aee} with $\lambda$ as fast angle.

\begin{figure}[!ht]
    \centering
    \includegraphics[width=0.5\textwidth]{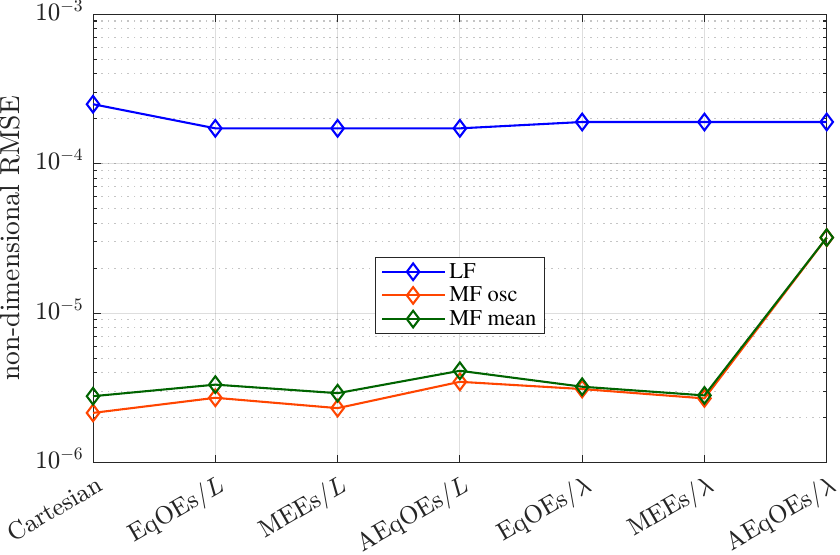}
    \caption{Non-dimensional \glsfmtshort{rmse} for \glsfmtlong{leo} case}
    \label{fig:rmse_leo}
\end{figure}

\begin{figure}[!ht]
    \centering
    \includegraphics[width=0.5\textwidth]{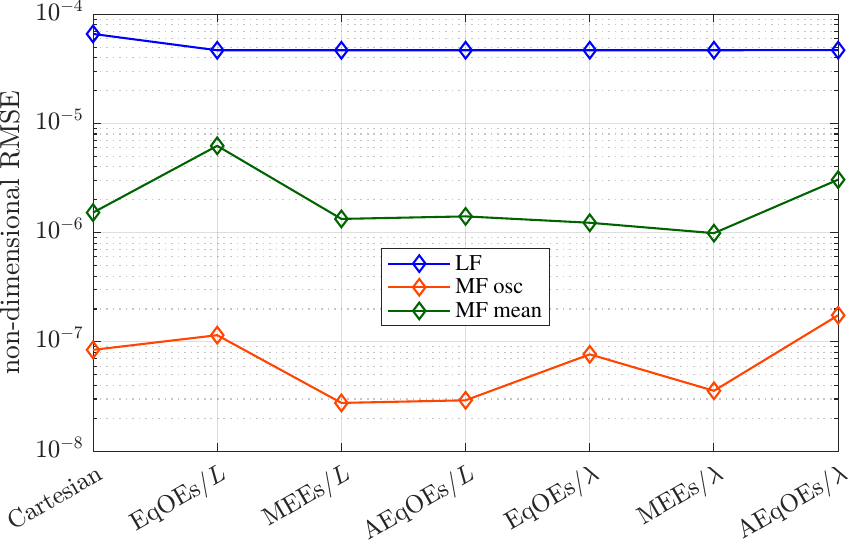}
    \caption{Non-dimensional \glsfmtshort{rmse} for \glsfmtlong{meo} case}
    \label{fig:rmse_meo}
\end{figure}

Similar conclusions can be drawn for the \gls{leo} case shown in \cref{fig:rmse_leo}, although the \glspl{rmse} of both \gls{mf} solutions are about one order of magnitude worse than those obtained for the \gls{heo} regime. Moreover, in this case a correction in osculating space always yields better results. Finally, the \gls{meo} case depicted in \cref{fig:rmse_meo} exhibits the best accuracy for a correction in osculating space, but also larger errors similar to the \gls{leo} regime when the correction is performed in mean elements.

\begin{figure}[!ht]
    \centering
    \includegraphics[width=0.6\textwidth]{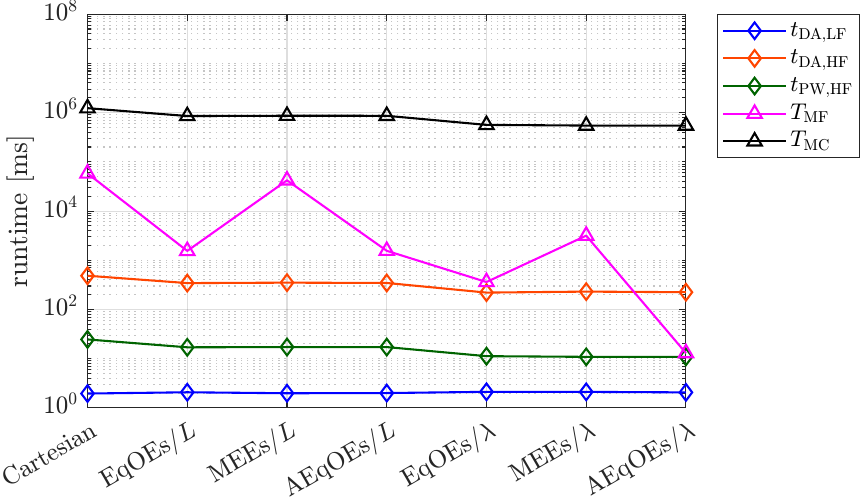}
    \caption{Runtimes for \glsfmtlong{heo} case}
    \label{fig:runtime_heo}
\end{figure}

Lastly, the computational efficiency of the \gls{mf} \gls{up} techniques is also analyzed. \Cref{fig:runtime_heo} reports the measured runtimes for the \gls{heo} test case. In this plot, the $x$ axis corresponds to a different state representation, while the $y$ axis reports the runtime in logarithmic scale expressed in milliseconds. The five curves correspond instead to five different tasks: those marked with diamonds refer to the propagation of a single trajectory (either \gls{pw} or in \gls{da}), while those marked with triangles correspond to the overall simulation. In particular, $t_{\mathrm{DA,LF}}$ (blue) and $t_{\mathrm{DA,HF}}$ (orange) represent the costs of propagating a single Taylor expansion using the \gls{lf} and the \gls{hf} propagators, respectively. Instead, $t_{\mathrm{PW,HF}}$ (green) is the cost of a single \gls{pw} propagation in \gls{hf}. Then, $t_{\mathrm{DA,LF}}$ (magenta) is the overall cost of the \gls{mf} method, while $t_{\mathrm{DA,HF}}$ (black) is that of the \gls{mc} simulation. For each case, $t_{\mathrm{DA,LF}}$ and $t_{\mathrm{DA,HF}}$ were estimated by propagating the \glspl{ic} for the selected time interval \num{2e3} times, and then taking the arithmetic mean of the single runtimes. $T_{\mathrm{MF}}$ is instead approximated as
\begin{equation}
    T_{\mathrm{MF}}\approx N\left(t_{\mathrm{DA,LF}}+t_{\mathrm{PW,HF}}\right)
\end{equation}
where $N$ is the number of kernels found in \cref{tab:nu_single_nb_kernels}. This approximation holds since the computational cost of the splitting operations and of the polynomial correction is negligible compared to both the evaluation of the \gls{sgp4} model in \gls{da} and the \gls{pw} propagation in \gls{hf} dynamics. Moreover, $T_{\mathrm{MC}}$ is the actual cost of the \gls{mc} simulation while $t_{\mathrm{PW,HF}}$ is approximated as the average runtime among \gls{mc} samples, i.e. $t_{\mathrm{PW,HF}}=T_{\mathrm{MC}}/N_{\mathrm{MC}}$.

By looking at \cref{fig:runtime_heo}, it is immediately noticed that $t_{\mathrm{DA,LF}}$ is almost constant across \glspl{oe}. This is expected since the cost of the \gls{lf} propagation is dominated by the evaluation of the \gls{sgp4} model, which processes only mean \gls{tle} elements independently on the selected coordinates. These elements are computed from the osculating parameters at the beginning of the propagation, but the cost of this operation is negligible compared to that of propagation. However, this observation does not apply to numerical propagation, whose performance is strongly affected by the choice of the state representation. This behavior is demonstrated by the orange and green curves, corresponding to $t_{\mathrm{DA,HF}}$ and $t_{\mathrm{PW,HF}}$, respectively. These results show in fact a correlation between the runtime and the type of parameters used to represent the orbit state. The propagations in Cartesian coordinates are characterized by the longest runtime, which is about \qty{30}{\percent} higher than the next worst case, i.e. the \glspl{mee} with the true longitude $L$ as fast angle. Both \gls{ee} and \gls{aee} with $L$ have instead runtimes that are similar to the latter case. Finally, the three fastest cases are those in \glspl{ee}, \glspl{mee}, and \glspl{aee} with the mean longitude $\lambda$ as fast angle. These propagations run about \qty{30}{\percent} faster compared to the ones with $L$. This behavior is explained by the average integration step size that is required to maintain the same accuracy for different parameters. Cartesian coordinates require the shortest step size, while for \glspl{oe} the time step is mainly driven by the fast angle, with $\lambda$ always resulting in larger step sizes than $L$. Choosing between \glspl{ee}, \glspl{mee}, and \glspl{aee} has otherwise almost no influence on the runtime. The magenta curve shows instead the same trends highlighted for a single trajectory, but now weighted by the number of kernels reported in \cref{tab:nu_single_nb_kernels}. The impact of the state parametrization is thus further emphasized in this case. The computational efficiency of the \gls{mf} method is instead assessed by looking at the runtime for a fix coordinates set. A \gls{hf} numerical propagation in \gls{da} is in fact about \num{20} times slower that its \gls{pw} counterpart. This number is very close to the number of coefficients of the Taylor polynomials being propagated, in this case equal to $\mathrm{dim}\left(_2D_6\right)=28$. Moreover, $t_{\mathrm{DA,LF}}$ varies between \qty{10}{\percent} and \qty{20}{\percent} of $t_{\mathrm{PW,HF}}$, thus demonstrating that the overall runtime of the \gls{mf} method is driven by the propagation of the kernels' means in \gls{hf} dynamics. Finally, $T_{\mathrm{MC}}$ is \num{20} times higher than $T_{\mathrm{MF}}$ in the worst case scenario, and up to four orders of magnitude higher in the best one.

\begin{figure}[!ht]
    \centering
    \includegraphics[width=0.6\textwidth]{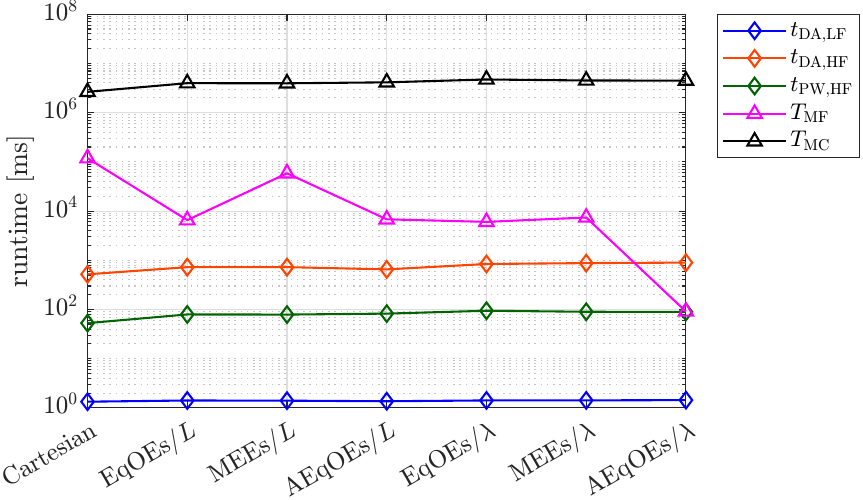}
    \caption{Runtimes for \glsfmtlong{leo} case}
    \label{fig:runtime_leo}
\end{figure}

\begin{figure}[!ht]
    \centering
    \includegraphics[width=0.6\textwidth]{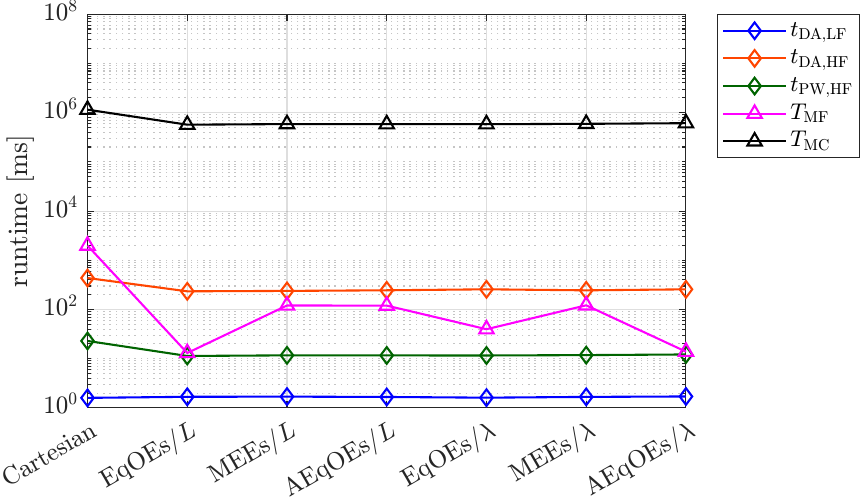}
    \caption{Runtimes for \glsfmtlong{meo} case}
    \label{fig:runtime_meo}
\end{figure}

\Cref{fig:runtime_leo,fig:runtime_meo} show the runtimes for the \gls{leo} and \gls{meo} cases, respectively. The main difference in the \gls{leo} scenario lies in the ratio between the computational costs of a \gls{da}-based and a \gls{pw} numerical propagation, which is now closer to \num{10} rather than \num{20}. Moreover, depending on the state representation, the \gls{mf} method is about \num{20} to \num{70} times faster than a \gls{mc} simulation. The reduced performance is explained by the fact that non-\gls{da} operations, such as frame transformations between \gls{gcrf} and \gls{itrf} for the computation of the gravitational and drag perturbations, are now more relevant. Conversely, the performance in the \gls{meo} scenario is close to that observed in the \gls{heo} one. However, given the quasi-circular nature of the orbit, the propagation in \glspl{oe} is not impacted by the choice of the fast angle. Finally, the weaker nonlinearities and the smaller number of kernels needed in this case result in a computational cost \num{500} times to four orders of magnitude smaller than that of a \gls{mc} simulation.

\subsection{Keplerian motion perturbed by stochastic accelerations}\label{sec:applications:keplerian_motion}

This section reproduces the scenario found in \citet{LopezYela2023}, which studies the dynamics of a planar Keplerian orbit perturbed by a stochastic acceleration inversely proportional to the instantaneous velocity of the object. The motion is studied in Cartesian coordinates projected onto the orbit's plane. The drift coefficient of the \gls{sde} is thus given by
\begin{equation}
    \vb*{u}\left(\vb*{X},t\right) = \left[v_x,\,v_y,\,-\dfrac{\mu}{r^3}x,\,-\dfrac{\mu}{r^3}y\right]^T
\end{equation}
with $\vb*{X}=\left[x\,y\,v_x\,v_y\right]^T$ the state vector of the \gls{so} and $r=\sqrt{x^2+y^2}$ the instantaneous orbit radius. The diffusion coefficient is instead given by
\begin{equation}
    \vb*{G}\left(\vb*{X},t\right)=
    \begin{bmatrix}
        \vb*{0}_{2\times 2}\\
        -\dfrac{\sigma_w}{v}\cdot\vb*{I}_{2\times 2}
    \end{bmatrix}
\end{equation}
where $\sigma_w>0$ is a constant scaling parameter and $v=\sqrt{v_x^2+v_y^2}$ is the instantaneous velocity of the object. Moreover, the random process $\vb*{W}(t)$ is a two-dimensional uncorrelated white noise with constant \gls{psd} $\vb*{Q}=h\cdot\vb*{I}_{2\times 2}$ and $h$ is the integration step size. At time $t_0$ the object is at the periapsis of a circular orbit, i.e.
\begin{equation}
    \label{eq:perturbed_kep_ic}
    \vb*{x}_0 = \left[R_\oplus+h_0,\,0,\,0,\,\sqrt{\dfrac{\mu_\oplus}{R_\oplus+h_0}}\right]^T
\end{equation}
with $R_\oplus$ the radius of the Earth, $h_0$ the initial orbit altitude, and $\mu_\oplus$ the Earth's standard gravitational parameter. These \glspl{ic} are then propagated to time $t_f$ using the methods described in \cref{sec:plasma:da,sec:multifidelity:stochastic}, respectively. \Gls{mc} simulations are also performed to verify the accuracy of the former. The latter are implemented in Julia using the \texttt{DifferentialEquations.jl} package~\cite{Rackauckas2017}. To reproduce the results presented in~\cite{LopezYela2023}, the \gls{em} method with a constant step size $h$ is used to approximate the solution to the \gls{sde}. The simulation parameters are summarized in \cref{tab:kepler_inputs}.

\begin{table}[!ht]
    \caption{Simulation parameters for perturbed Keplerian motion}
    \label{tab:kepler_inputs}
    \centering
    \small
    \begin{tabular}{cclccl}
        \toprule
        Parameter & Value & Description & Parameter & Value & Description\\ \midrule
        $\mu_\oplus$ & \qty{3.986e5}{\km\cubed\per\s\squared} & standard gravitational parameter & $t_0$ & \qty{0.0}{\s} & initial time\\
        $R_\oplus$ & \qty{6378.0}{\km} & equatorial radius & $t_f$ & \qty{1.51}{\day} & final time\\
        $h_0$ & \qty{200.0}{\km} & initial orbit altitude & $h$ & \qty{0.1}{\s} & \glsfmtshort{em} step size\\
        $\sigma_w$ & \qty{2e-4}{\km\squared\per\s\tothe{2.5}} & scaling parameter in diffusion term & $N$ & 2 & \gls{da} expansion order\\
        \bottomrule
    \end{tabular}
\end{table}

The first simulation considers deterministic \glspl{ic}, namely $\vb*{x}_0$ given by \cref{eq:perturbed_kep_ic}, and uses the \gls{plasma} method to estimate $\mathbb{E}\left[\hat{\vb*{X}}^{\vb*{r}}_{f,2}\left(\vb*{x}_0\right)\right]$ at final time $t_f$. The results are reported in \cref{tab:mean_kep_det,tab:cov_kep_det}, which also include those of a \gls{mc} simulation with $N_{\mathrm{MC}}=\num{e4}$ samples and the relative errors between the two solutions.

\begin{table}[!ht]
    \caption{Estimated mean for perturbed Keplerian motion with deterministic \glsfmtlongpl{ic}}
    \label{tab:mean_kep_det}
    \centering
    \small
    \sisetup{round-mode=places,round-precision=3}
    \begin{tabular}{r *4{S[scientific-notation=true,table-format=-1.3e-1]}}
        \toprule
        & {$x,\unit{\km}$} & {$y,\unit{\km}$} & {$v_x,\unit{\km\per\s}$} & {$v_y,\unit{\km\per\s}$} \\ \midrule
        \Glsfmtshort{plasma} & 5935.7143264343895 & -2661.641732752862 & 2.983302274246626 & 6.669976955926232 \\
        \Glsfmtlong{mc} & 5956.801345584536 & -2640.4396328693633 & 2.960111062410341 & 6.693927007247425 \\
        Relative errors & 0.0035525663787821355 & 0.00796579780914765 & 0.007773671490309048 & 0.0035907247475441707 \\ \bottomrule
    \end{tabular}
\end{table}

\begin{table}[!ht]
    \caption{Estimated covariance for perturbed Keplerian motion with deterministic \glsfmtlongpl{ic}. Units are \unit{\km\squared}, \unit{\km\squared\per\s}, and \unit{\km\squared\per\s\squared}, respectively}
    \label{tab:cov_kep_det}
    \centering
    \small
    \sisetup{round-mode=places,round-precision=3}
    \begin{tabular}{c *4{S[scientific-notation=true,table-format=-1.3e-1]}}
        \toprule
        \multirow{4}{*}{\Glsfmtshort{plasma}} & 673006.6416414595 & 1.5321654153871248e6 & -1713.156540777673 & 773.0689107862491 \\
        & 1.5321654153871248e6 & 3.489308726003569e6 & -3901.038916556096 & 1761.0107930683137 \\
        & -1713.1565407776727 & -3901.038916556096 & 4.361708803243799 & -1.9687186678950324 \\
        & 773.0689107862491 & 1761.0107930683137 & -1.9687186678950324 & 0.888966660103023 \\ \midrule
        \multirow{4}{*}{\shortstack[c]{Monte\\ Carlo}} & 747283.5878591138 & 1.3219953320955911e6 & -1475.4163279421768 & 851.8519419287477 \\
        & 1.3219953320955911e6 & 3.1802887152576386e6 & -3558.212500960482 & 1518.0735027225319 \\
        & -1475.416327942177 & -3558.212500960482 & 3.981533863013293 & -1.6940980372372914 \\
        & 851.8519419287477 & 1518.0735027225319 & -1.6940980372372914 & 0.9720889065279227 \\ \midrule
        \multirow{4}{*}{\shortstack[c]{Relative\\ errors}} & 0.11036584429017403 & 0.13717192750916585 & 0.13877319858206066 & 0.1019094546983817 \\
        & 0.13717192750916585 & 0.088561957399471 & 0.08788079866126201 & 0.13795332277464223 \\
        & 0.13877319858206041 & 0.08788079866126201 & 0.08716192606617157 & 0.13949206412075488 \\
        & 0.1019094546983817 & 0.13795332277464223 & 0.13949206412075488 & 0.09350434628815722 \\ \bottomrule
    \end{tabular}
\end{table}

These values are in agreement with those reported by \citet{LopezYela2023} and demonstrate the capability of the algorithm to accurately propagate the statistics through nonlinear \glspl{sde}. Moreover, the runtime of the \gls{plasma} method is about \qty{23}{\s}, while the \gls{mc} simulation takes around \qty{8538}{\s} to complete. This demonstrates the computational efficiency of the former which is about \num{370} times faster than the \gls{mc} simulation.

The second scenario considers random \glspl{ic} and uses the hybrid algorithm discussed in \cref{sec:multifidelity:stochastic} to estimate the moments of $\vb*{X}$ at final time $t_f$. The initial state is assumed to be normally distributed with mean $\vb*{\mu}_0=\vb*{x}_0$ given by \cref{eq:perturbed_kep_ic} and a diagonal covariance matrix $\vb*{P}_0$ equal to~\cite{LopezYela2023}
\begin{equation}
    \label{eq:perturbed_kep_cov_ic}
    \vb*{P}_0 = \mathrm{diag}\left(\qty{1e-1}{\km\squared},\,\qty{1e-1}{\km\squared},\,\qty{1e-4}{\km\squared\per\s\squared},\,\qty{1e-4}{\km\squared\per\s\squared}\right).
\end{equation}
In this case, \cref{eq:p_x0_gmm_pn} reduces to a single kernel, i.e.
\begin{equation}
    p_{\vb*{X}}(\vb*{x}_0)=p_g\left(\vb*{x}_0;\vb*{\mu}_0,\vb*{P}_0\right)
\end{equation}
with $\vb*{\mu}_0$ and $\vb*{P}_0$ given by \cref{eq:perturbed_kep_ic,eq:perturbed_kep_cov_ic}, respectively.

The adaptive \gls{gmm} method is run with $\varepsilon_\nu=0.05,\alpha_{min}=\num{e-3}$, $N_{max,e}=4$, and the same three-components splitting library used in \cref{sec:applications:deterministic}. A total of \num{81} kernels are generated by the algorithm in this case. The first two moments of the transformed \gls{pdf} output of the \gls{mf} method are then summarized in \cref{tab:mean_kep_mvn,tab:cov_kep_mvn}, respectively. These tables also include the moments estimated from the \gls{mc} samples and the relative errors between the two solutions.

\begin{table}[!ht]
    \caption{Estimated mean for perturbed Keplerian motion with random \glsfmtlongpl{ic}}
    \label{tab:mean_kep_mvn}
    \centering
    \small
    \sisetup{round-mode=places,round-precision=3}
    \begin{tabular}{r *4{S[scientific-notation=true,table-format=-1.3e-1]}}
        \toprule
        & {$x,\unit{\km}$} & {$y,\unit{\km}$} & {$v_x,\unit{\km\per\s}$} & {$v_y,\unit{\km\per\s}$} \\ \midrule
        \Glsfmtshort{plasma}-\glsfmtshort{gmm} & 5148.353360230148 & -2324.3869247483676 & 2.5978679254216206 & 5.7863241481057806 \\
        \Glsfmtlong{mc} & 5093.137660008753 & -2301.814317774163 & 2.5725520612264456 & 5.725093105034937 \\
        Relative errors & 0.010724924331714108 & 0.009711208892920516 & 0.009744861910586282 & 0.010582027813095786 \\ \bottomrule
    \end{tabular}
\end{table}

\begin{table}[!ht]
    \caption{Estimated covariance for perturbed Keplerian motion with random \glsfmtlongpl{ic}. Units are \unit{\km\squared}, \unit{\km\squared\per\s} and \unit{\km\squared\per\s\squared} respectively}
    \label{tab:cov_kep_mvn}
    \centering
    \small
    \sisetup{round-mode=places,round-precision=3}
    \begin{tabular}{c *4{S[scientific-notation=true,table-format=-1.3e-1]}}
        \toprule
        \multirow{4}{*}{\shortstack[c]{\Glsfmtshort{plasma}-\\ \glsfmtshort{gmm}}} & 3.832968233607461e6 & 3.9913498134718616e6 & -4447.531374024072 & 4290.651526314648 \\
        & 3.9913498134718616e6 & 1.0814741492866285e7 & -12124.582140011544 & 4473.349406217434 \\
        & -4447.531374024072 & -12124.582140011544 & 13.594111371869772 & -4.984348312111524 \\
        & 4290.651526314648 & 4473.349406217434 & -4.984348312111524 & 4.804415927538659 \\ \midrule
        \multirow{4}{*}{\shortstack[c]{Monte\\ Carlo}} & 4.352556700285577e6 & 3.8948857592090755e6 & -4334.897012214353 & 4866.208194372283 \\
        & 3.8948857592090755e6 & 1.0838569905397894e7 & -12156.066596619914 & 4355.265233129937 \\
        & -4334.897012214354 & -12156.066596619914 & 13.634858455883444 & -4.847058918639009 \\
        & 4866.208194372282 & 4355.265233129937 & -4.847058918639009 & 5.441825644903913 \\ \midrule
        \multirow{4}{*}{\shortstack[c]{Relative\\ errors}} & 0.13555772837415267 & 0.024168278595174612 & 0.0253251415982275 & 0.13414202121233448 \\
        & 0.024168278595174612 & 0.00220332705569769 & 0.002596745705938222 & 0.02639726128331804 \\
        & 0.025325141598227296 & 0.002596745705938222 & 0.0029974069579854884 & 0.027544101028997906 \\
        & 0.13414202121233426 & 0.02639726128331804 & 0.027544101028997906 & 0.13267163521618835 \\ \bottomrule
    \end{tabular}
\end{table}

These results demonstrate the capabilities of the proposed method to not only propagate the effects of stochastic accelerations on the moments of the transformed \gls{pdf}, but also to accurately propagate the initial uncertainty and combine the two solutions into a \gls{gmm} approximation of the transformed \gls{pdf}. As before, these results also agree with those found in~\cite{LopezYela2023}. In this case, the propagation of the initial uncertainty takes about \qty{297}{\s} to complete, while the propagation of the effective noise moments requires around \qty{23}{\s} for each kernel generated in the previous step. The overall runtime is thus about \qty{2160}{\s}, or \SI{25}{\percent} of that of a \gls{mc} simulation. The smaller gain in performance with respect to the first simulation is due to the fact that the computational cost of the \gls{mc} simulation is independent on the \glspl{ic} (deterministic or random), while that of the \gls{da}-based algorithm is strongly influenced by the number of kernels generated by the adaptive \gls{gmm} method. Indeed, the overall runtime is dominated by the computation of \gls{pn} effects, which accounts for about \qty{86}{\percent} of the execution time of the \gls{mf} solution. This cost can however be reduced by choosing a different parametrization of the orbit state, i.e. \glspl{ee} or one of its variants, such that the initial uncertainty can be accurately mapped with a lower number of kernels.

\subsection{Low-thrust orbit raising maneuver}\label{sec:applications:low_thrust}

The last test case considers a spacecraft on a \gls{gto} subject to a continuous and constant tangential thrust acceleration $\vb*{a}_t$. The drift term $\vb*{u}(\vb*{x},t)$ models the central body attraction, the gravitational perturbations due to the non-uniform Earth gravity potential, and the nominal thrust law realization. Its expression in Cartesian coordinates is thus given by
\begin{equation}
    \vb*{u}(\vb*{x},t) = -\dfrac{\mu}{r^3}\vb*{r} + \vb*{a}_{g}(\vb*{r},t) + a_t\dfrac{\vb*{v}}{v}
    \label{eq:drift_thrust}
\end{equation}
with $\vb*{r}=[x\ y\ z]^T$ and $\vb*{v}=[v_x\ v_y\ v_z]^T$ the instantaneous position and velocity vectors, $r=\norm{\vb*{r}}_2$ and $v=\norm{\vb*{v}}_2$ their respective norms, $\mu$ the standard gravitational parameter, $\vb*{a}_{g}(\vb*{r},t)$ the gravitational perturbations, and $a_t$ the thrust acceleration magnitude. Instead, the diffusion term $\vb*{G}(\vb*{x},t)$ models the effects of thrust dispersions as small deviations with respect to the nominal thrust law realization~\cite{Maestrini2023}. To obtain an explicit expression for $\vb*{G}(\vb*{x},t)$, the theoretical thrust acceleration is firstly expressed in polar coordinates $\vb*{y}=[a_t\ \alpha\ \beta]^T$, where $\alpha=\arctan{(v_y/v_x)}$ and $\beta=\arcsin{(v_z/v)}$ are the two angles that define the thrust direction in the spacecraft's \gls{tnw} frame~\cite{Vallado2013}. The Cartesian acceleration vector is thus obtained as $\vb*{a}_t=a_t\cdot\left[\cos{\alpha}\cos{\beta},\,\sin{\alpha}\cos{\beta},\,\sin{\beta}\right]^T$. The acceleration magnitude $a_t$ and the angles $\alpha,\beta$ are then modeled as independent Gaussian random variables with means coincident with their nominal values and standard deviations equal to $\sigma_{a_t},\sigma_{\alpha},\sigma_{\beta}$, respectively. The diffusion matrix $\vb*{G}(\vb*{x},t)$ is finally given by
\begin{equation}
    \vb*{G}(\vb*{x},t) = \begin{bmatrix}
        \vb*{0}_{3\times 3}\\
        \vb*{J}\sqrt{\vb*{\Sigma}}\\
    \end{bmatrix}
    \label{eq:diffusion_thrust}
\end{equation}
with $\vb*{J}=\partial\vb*{a}_t/\partial\vb*{y}$ the partials of $\vb*{a}_t$ with respect to $\vb*{y}$ evaluated along the nominal trajectory and $\vb*{\Sigma}=\text{diag}(\sigma_{a_t}^2,\sigma_{\alpha}^2,\sigma_{\beta}^2)$. Together, \cref{eq:drift_thrust,eq:diffusion_thrust} constitute the \gls{hf} stochastic dynamical model used to propagate the \gls{pn} effects with the \gls{plasma} method of \cref{sec:plasma:da}. Conversely, if the \gls{bf} formulation described in \cref{sec:multifidelity:stochastic} is exploited, the full \cref{eq:drift_thrust} is still used to obtain the dense output model $\vb*{\xi}(t)$, but the gravitational perturbations $\vb*{a}_{g}(\vb*{r},t)$ are neglected in the subsequent steps. The diffusion term given by \cref{eq:diffusion_thrust} is instead identical in the two cases. The propagation of the initial uncertainty also considers two different models of the deterministic dynamics: the \gls{hf} simulation integrates \cref{eq:drift_thrust} numerically, while the \gls{lf} propagation in the \gls{mf} method exploits the closed-form solution for the two-body problem with constant tangential thrust acceleration derived by \citet{Bombardelli2011}. In both cases, the uncertainty in the \glspl{ic} is modeled as a \gls{gmm} and propagated using the adaptive method presented in \cref{sec:background:gmm}. The simulation parameters are summarized in \cref{tab:thrust_params}, while the nominal \glspl{ic} at $t_0$ are expressed in Keplerian elements and given in \cref{tab:kepler_ics_thrust}. The initial uncertainty is instead expressed in Cartesian coordinates and assumes a diagonal covariance matrix $\vb*{P}_0$ with standard deviations $\sigma_{\vb*{r}}=\qty{1e2}{\m}$ and $\sigma_{\vb*{v}}=\qty{1e-1}{\m\per\s}$ on the components of the spacecraft's position and velocity vectors, respectively.

\begin{table}[!ht]
    \caption{Simulation parameters for low-thrust orbit raising maneuver}
    \label{tab:thrust_params}
    \begin{adjustwidth}{-1.5cm}{-1.5cm}
    \centering
    \small
    \begin{tabular}{cclccl}
        \toprule
        Parameter & Value & Description & Parameter & Value & Description \\ \midrule
        $\mu_\oplus$ & \qty{3.986004418e5}{\km\cubed\per\s\squared} & standard gravitational parameter & $t_0$ & \qty{0.0}{\s} & initial time\\
        $R_\oplus$ & \qty{6378.137}{\km} & equatorial radius & $t_f$ & \qty{2.0}{\day} & final time\\
        $J_2$ & \num{1.0826269e-3} & second zonal harmonic coefficient & $\sigma_{a_t}$ & \qty{5e-6}{\m\per\s\tothe{1.5}} & standard deviation on $a_t$\\
        $a_t$ & \qty{1e-4}{\m\per\s\squared} & thrust acceleration magnitude & $\sigma_{\alpha},\sigma_{\beta}$ & \qty{5.0}{\deg\s\tothe{0.5}} & standard deviation on $\alpha,\beta$\\
        \bottomrule
    \end{tabular}
    \end{adjustwidth}
\end{table}

\begin{table}[!ht]
    \caption{Initial Keplerian parameters for low-thrust orbit raising maneuver}
    \label{tab:kepler_ics_thrust}
    \centering
    \small
    \begin{tabular}{cccccc}
        \toprule
        $a,\unit{\km}$ & $e,-$ & $i,\unit{\deg}$ & $\Omega,\unit{\deg}$& $\omega,\unit{\deg}$& $\nu,\unit{\deg}$\\ \midrule
        \num{2.4e4} & \num{0.72} & \num{5.0} & \num{0.0} & \num{0.0} & \num{0.0}\\ \bottomrule
    \end{tabular}
\end{table}

\Glspl{ic} in \cref{tab:kepler_ics_thrust} are firstly converted to Cartesian parameters, and the initial state vector initialized in the \gls{da} framework using \cref{eq:da_mu_cov_def} with $\zeta=3$ to obtain $[\vb*{x}_0]=\mathcal{T}_{\vb*{x}_0}(\var{\vb*{x}})$. All propagations are then carried out in \glspl{mee} with the mean longitude $\lambda$ as fast variable. Since the initial uncertainty is expressed in Cartesian coordinates, $[\vb*{x}_0]$ must first be converted to \glspl{mee}. The required transformation is embedded in the adaptive \gls{gmm} algorithm, such that all higher order information is retained and there is no need to evaluate the initial covariance in \glspl{mee}. The nonlinearity threshold is then set equal to $\varepsilon_\nu=0.01$, and the same three-components splitting library of \cref{sec:applications:deterministic} is used. Two simulations are then conducted: the first one highlights the performance of the \gls{mf} formulation of the adaptive \gls{gmm} and \gls{plasma} algorithms with respect to their \gls{hf} counterparts. These methods are then compared with a reference \gls{mc} simulation in the second case.

In the first case, the perturbing acceleration $\vb*{a}_{g}(\vb*{r},t)$ models the Earth's gravitational potential with a spherical harmonics model truncated at order and degree 8~\cite{Holmes2002}. Numerical propagations use the adaptive \gls{dop853} method, and the performance of the \gls{mf} technique is compared to that of the \gls{hf} solution.

\begin{table}[!ht]
    \caption{Estimated mean and relative errors for low-thrust orbit raising maneuver with $8\times 8$ Earth's gravity model}
    \label{tab:mean_tangential_thrust}
    \centering
    \small
    \setlength{\tabcolsep}{3pt}
    \sisetup{round-mode=places,round-precision=3}
    \begin{tabular}{r *6{S[scientific-notation=true,table-format=-1.3e-2]}}
        \toprule
        & {$p,\unit{\km}$} & {$f,-$} & {$g,-$} & {$h,-$} & {$k,-$} & {$\lambda,\unit{\radian}$} \\ \midrule
        \Glsfmtshort{mf} estimate & 11711.3062650012 & 0.716764464351676 & 0.0101598249111405 & 0.0436129594368403 & -0.000602280573800692 & 4.21272980990873 \\
        \Glsfmtshort{mf}--\glsfmtshort{hf} & 1.57094465395873e-09 & 6.04686431082094e-10 & 1.93531935474256e-08 & 1.26657500255493e-09 & 5.21991063209376e-08 & 5.64345717490392e-10 \\ \bottomrule
    \end{tabular}
\end{table}

\begin{table}[!ht]
    \caption{Estimated covariance and relative errors for low-thrust orbit raising maneuver with $8\times 8$ Earth's gravity model. Units are \unit{\km\squared}, \unit{\km\radian} and \unit{\radian\squared} respectively}
    \label{tab:cov_tangential_thrust}
    \centering
    \small
    \setlength{\tabcolsep}{3pt}
    \sisetup{round-mode=places,round-precision=3}
    \begin{tabular}{S[scientific-notation=true,table-format=-1.3e-1] *5{S[scientific-notation=true,table-format=-1.3e-2]}}
        \toprule
        \multicolumn{6}{c}{\Glsfmtlong{mf} estimate in \glsfmtshortpl{mee}} \\ \midrule
        0.184414082685874 & 1.74546207623176e-05 & 4.2018170854077e-07 & -2.25918226573567e-11 & 1.27404630926359e-11 & -0.00293549341925482 \\
        1.74546207623176e-05 & 1.85869382334468e-09 & 4.32851884054799e-11 & 2.00071672705952e-15 & -2.51130348348037e-15 & -3.04481147076765e-07 \\
        4.2018170854077e-07 & 4.32851884054799e-11 & 5.10367728050375e-10 & -7.54764411235044e-16 & 3.48346064693334e-12 & -6.57209793740098e-09 \\
        -2.25918226573567e-11 & 2.00071672705952e-15 & -7.54764411235044e-16 & 2.52443165100153e-11 & -7.67681521580102e-15 & 2.67589163442728e-15 \\
        1.27404630926359e-11 & -2.51130348348037e-15 & 3.48346064693334e-12 & -7.67681521580102e-15 & 5.5668168415732e-11 & 4.85447491715538e-12 \\
        -0.00293549341925482 & -3.04481147076765e-07 & -6.57209793740098e-09 & 2.67589163442728e-15 & 4.85447491715538e-12 & 5.01632867284463e-05 \\ \midrule
        \multicolumn{6}{c}{Relative errors} \\ \midrule
        0.0064777874226779 & 0.00153804347074266 & 0.057471960215693 & 1.02329906103978 & 0.999425888520201 & 0.00606739916230672 \\
        0.00153804347074266 & 0.00405228057125205 & 0.0287311672821463 & 0.979662667047158 & 1.00121386134278 & 0.000536918612601137 \\
        0.057471960215693 & 0.0287311672821463 & 0.0024071174283305 & 0.993961704540486 & 0.0142794569347485 & 0.116025251119152 \\
        1.02329906103978 & 0.979662667047158 & 0.993961704540486 & 0.000703896517651289 & 1.01796686607799 & 1.00016426951061 \\
        0.999425888520201 & 1.00121386134278 & 0.0142794569347485 & 1.01796686607799 & 0.00089362947434442 & 1.01379487578611 \\
        0.00606739916230672 & 0.000536918612601137 & 0.116025251119152 & 1.00016426951061 & 1.01379487578611 & 0.00510194450993468 \\
        \bottomrule
    \end{tabular}
\end{table}

The estimated mean and covariance after \qty{2.0}{\day}, or about \num{4.7} revolutions, are reported in \cref{tab:mean_tangential_thrust,tab:cov_tangential_thrust}, respectively. The relative errors between the two estimates are also included. With a nonlinearity threshold equal to $\varepsilon_\nu=0.01$, a single polynomial is sufficient to describe the transformed uncertainty. \Cref{tab:mean_tangential_thrust} shows a very good agreement between the two solutions, with relative errors always below \num{6e-8}. The estimated covariance is also globally accurate, with a maximum relative error on the diagonal entries of about \qty{0.5}{\percent}. An exception is however found in the cross terms involving the $h$ and $k$ parameters, but these errors are traced back to the first step of the \gls{mf} technique.

\begin{table}[!ht]
    \caption{Runtime for low-thrust orbit raising maneuver with $8\times 8$ Earth's gravity model}
    \label{tab:runtime_tangential_thrust}
    \centering
    \small
    \sisetup{round-mode=places,round-precision=3}
    \begin{tabular}{r S[scientific-notation=false,table-format=1.3] S[scientific-notation=false,table-format=3.3] S[scientific-notation=false,table-format=1.3] S[scientific-notation=false,table-format=3.3]}
        \toprule
        & {$t_{\mathrm{MF}}$, \unit{\s}} & {$t_{\mathrm{MF}}$, \unit{\percent}} & {$t_{\mathrm{HF}}$, \unit{\s}} & {$t_{\mathrm{HF}}$, \unit{\percent}} \\ \midrule
        Adaptive \glsfmtshort{gmm} & 0.00221 & 0.387257307072265 & 0.8549 & 15.188069841315 \\
        \Glsfmtshort{plasma} & 0.56847 & 99.6127426929277 & 4.77386 & 84.811930158685 \\
        Total & 0.57068 & 100 & 5.62876 & 100 \\ \bottomrule
    \end{tabular}
\end{table}

\Cref{tab:runtime_tangential_thrust} reports the absolute runtime required to compute the former estimates and the relative cost of the adaptive \gls{gmm} and \gls{plasma} methods. The advantage of the \gls{mf} solution is readily seen by comparing the values in the third row, which correspond to the complete solutions computed with the \gls{mf} and the \gls{hf} methods, respectively. An increase by a factor of \num{10} is in fact observed in the runtime needed to carry out the two simulations. Breaking down the execution time of the \gls{mf} technique, the propagation of the initial uncertainty using the analytical model takes only \qty{0.4}{\percent} of the total run time, while the remaining \qty{99.6}{\percent} is allocated for the propagation of the \gls{pn} effects. Concerning the \gls{hf} solution, the fraction of runtime required for the propagation of the initial uncertainty grows to \qty{15}{\percent}, while that needed by the \gls{plasma} method drops to about \qty{85}{\percent}.

In the second test case, the two methods discussed above are compared with a reference \gls{mc} simulation with $N_{\mathrm{MC}}=\num{e4}$ samples. The latter is carried out using the \gls{sde} solver NON~\cite{Komori2007} implemented in the \texttt{DifferentialEquations.jl} package~\cite{Rackauckas2017,Rackauckas2017a,Rackauckas2020}. This solver was selected as one of the few methods capable of handling non-diagonal \gls{pn} models. Moreover, to match its deterministic order (equal to 4) and fixed step implementation, the adaptive \gls{gmm} and \gls{plasma} methods were coupled with a classical fourth-order \gls{rk} integrator. A fixed step size equal to \qty{60}{\s} was then selected for each numerical propagation. Finally, the gravitational perturbations were simplified to include only the effects of the Earth's oblateness modeled by the $J_2$ zonal harmonic coefficient~\cite{Curtis2021}. Also in this case, the propagation of the initial uncertainty did not require any split.

\begin{table}[!ht]
    \caption{Estimated mean and relative errors for low-thrust orbit raising maneuver with $J_2$ model}
    \label{tab:mean_tangential_thrust_mc}
    \centering
    \small
    \setlength{\tabcolsep}{3pt}
    \sisetup{round-mode=places,round-precision=3}
    \begin{tabular}{r *6{S[scientific-notation=true,table-format=-1.3e-2]}}
        \toprule
        & {$p,\unit{\km}$} & {$f,-$} & {$g,-$} & {$h,-$} & {$k,-$} & {$\lambda,\unit{\radian}$} \\ \midrule
        {\Glsfmtshort{mc} estimate} & 11711.2299463426 & 0.716737858407749 & 0.0101192991755733 & 0.0436357167799616 & -0.000601165558870795 & 4.21648712218124 \\
        {\Glsfmtshort{mc}--\glsfmtshort{gmm} \glsfmtshort{mf}} & 2.401548648068e-07 & 4.73125469635793e-07 & 4.11298155797927e-06 & 9.25375924319868e-07 & 0.000189218028208465 & 1.29051124415713e-05 \\
        {\Glsfmtshort{mc}--\glsfmtshort{gmm} \glsfmtshort{hf}} & 2.38486108764444e-07 & 4.73740488156653e-07 & 4.1410602585368e-06 & 9.26243913768576e-07 & 0.000189185211767262 & 1.29056162803024e-05 \\ \bottomrule
    \end{tabular}
\end{table}

\begin{table}[!ht]
    \caption{Estimated covariance and relative errors for low-thrust orbit raising maneuver with $J_2$ model. Units are \unit{\km\squared}, \unit{\km\radian} and \unit{\radian\squared} respectively}
    \label{tab:cov_tangential_thrust_mc}
    \centering
    \small
    \setlength{\tabcolsep}{3pt}
    \sisetup{round-mode=places,round-precision=3}
    \begin{tabular}{S[scientific-notation=true,table-format=-1.3e-1] *5{S[scientific-notation=true,table-format=-1.3e-2]}}
        \toprule
        \multicolumn{6}{c}{\Glsfmtlong{mc} estimate in \glsfmtshortpl{mee}} \\ \midrule
        0.183941574422538 & 1.72330573690184e-05 & 6.33394610776206e-07 & 3.19057760213867e-08 & 1.27035778314213e-08 & -0.00291502297348677 \\
        1.72330573690184e-05 & 1.81816613182522e-09 & 6.40174923693768e-11 & 2.31689898773625e-12 & 1.95949386390646e-12 & -2.99443112335884e-07 \\
        6.33394610776206e-07 & 6.40174923693768e-11 & 5.05575724098346e-10 & -9.06625199479893e-14 & 2.24718731224521e-12 & -1.09745366318734e-08 \\
        3.19057760213867e-08 & 2.31689898773625e-12 & -9.06625199479893e-14 & 2.52697462530889e-11 & 6.29887708229726e-13 & -4.18050959052075e-10 \\
        1.27035778314213e-08 & 1.95949386390646e-12 & 2.24718731224521e-12 & 6.29887708229726e-13 & 5.58174547661188e-11 & -3.01577222372589e-10 \\
        -0.00291502297348677 & -2.99443112335884e-07 & -1.09745366318734e-08 & -4.18050959052075e-10 & -3.01577222372589e-10 & 4.96041163206348e-05 \\ \midrule
        \multicolumn{6}{c}{Relative errors \glsfmtshort{mc}--\glsfmtshort{gmm} \glsfmtshort{mf}} \\ \midrule
        0.0025765089885667 & 0.0128544397943363 & 0.336632975333107 & 1.00083052945318 & 0.99919386065176 & 0.00702249689791772 \\
        0.0128544397943363 & 0.0222912024985818 & 0.323853538172251 & 0.999079258439641 & 1.00119054956153 & 0.0168246936267404 \\
        0.336632975333107 & 0.323853538172251 & 0.00948123780625063 & 0.99168421960887 & 0.550147732720262 & 0.401151007152774 \\
        1.00083052945318 & 0.999079258439641 & 0.99168421960887 & 0.00100340817492971 & 1.012158936377 & 1.00000245232592 \\
        0.99919386065176 & 1.00119054956153 & 0.550147732720262 & 1.012158936377 & 0.0026717647215391 & 1.01609709156858 \\
        0.00702249689791772 & 0.0168246936267404 & 0.401151007152774 & 1.00000245232592 & 1.01609709156858 & 0.0112726657047919 \\ \midrule
        \multicolumn{6}{c}{Relative errors \glsfmtshort{mc}--\glsfmtshort{gmm} \glsfmtshort{hf}} \\ \midrule
        0.0108555153547155 & 0.0151911942971494 & 0.312784930762825 & 0.894074354092646 & 0.61947644139459 & 0.0141829144315531 \\
        0.0151911942971494 & 0.0178017216422998 & 0.360537329292307 & 0.847501963887071 & 0.0288749701744656 & 0.0172591252784961 \\
        0.312784930762825 & 0.360537329292307 & 0.0118393494291306 & 0.279505908841092 & 0.564528120190932 & 0.339917565096704 \\
        0.894074354092646 & 0.847501963887071 & 0.279505908841092 & 0.00169826801868174 & 0.330791206846028 & 0.860651354383523 \\
        0.61947644139459 & 0.0288749701744656 & 0.564528120190932 & 0.330791206846028 & 0.00178338385065724 & 0.0755391567115838 \\
        0.0141829144315531 & 0.0172591252784961 & 0.339917565096704 & 0.860651354383523 & 0.0755391567115838 & 0.016600773588164 \\ \bottomrule
    \end{tabular}
\end{table}

\Cref{tab:mean_tangential_thrust_mc,tab:cov_tangential_thrust_mc} report the empirical mean and covariance matrix obtained from the \gls{mc} simulation for a propagation time of \qty{2.0}{\day}. These tables also include the relative errors between the \gls{mc} estimates and the \gls{mf} and \gls{hf} solutions obtained with the coupled adaptive \gls{gmm} and \gls{plasma} methods. \Cref{tab:mean_tangential_thrust_mc} shows a very good agreement between the reference solution and these two estimates. Indeed, the relative errors are always below \num{2e-4}, with four out of the six components exhibiting errors below \num{5e-6}. Moreover, the errors reported in the second and third columns are very similar to each others, meaning that the difference between the \gls{mf} and \gls{hf} solutions is at least one order of magnitude smaller than that between the \gls{mc} estimate and either one of these two solutions. The estimated covariances are also globally accurate, as demonstrated in \cref{tab:cov_tangential_thrust_mc}. The relative errors follow the same pattern already observed in \cref{tab:cov_tangential_thrust}, with a very good agreement on the diagonal entries and the largest errors on the cross terms involving the $h$ and $k$ parameters.

\begin{figure}[p]
    \centering
    \subcaptionbox{$p$ parameter\label{subfig:hist_thrust_p}}[0.5\textwidth]{\includegraphics[width=0.4\textwidth]{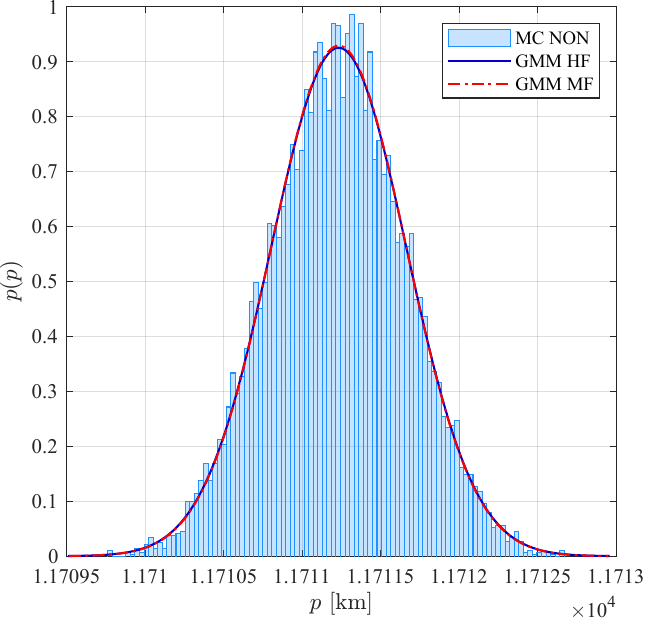}}%
    \hfill
    \subcaptionbox{$\lambda$ parameter\label{subfig:hist_thrust_l}}[0.5\textwidth]{\includegraphics[width=0.4\textwidth]{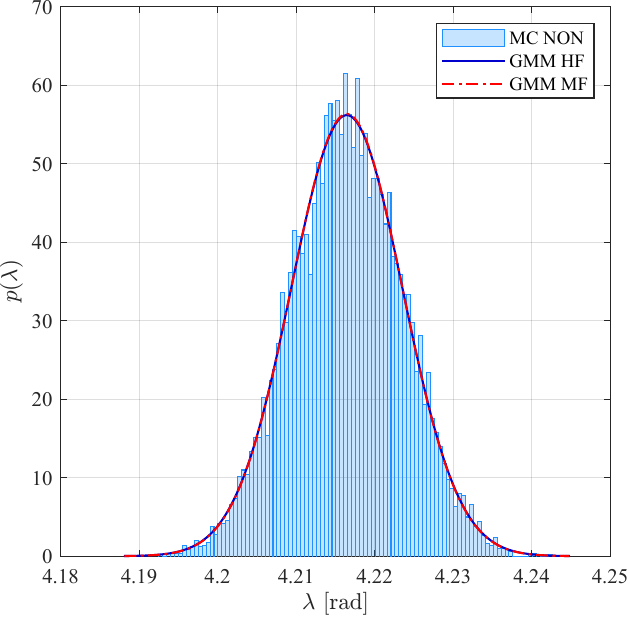}}
    \subcaptionbox{$f$ parameter\label{subfig:hist_thrust_f}}[0.5\textwidth]{\includegraphics[width=0.4\textwidth]{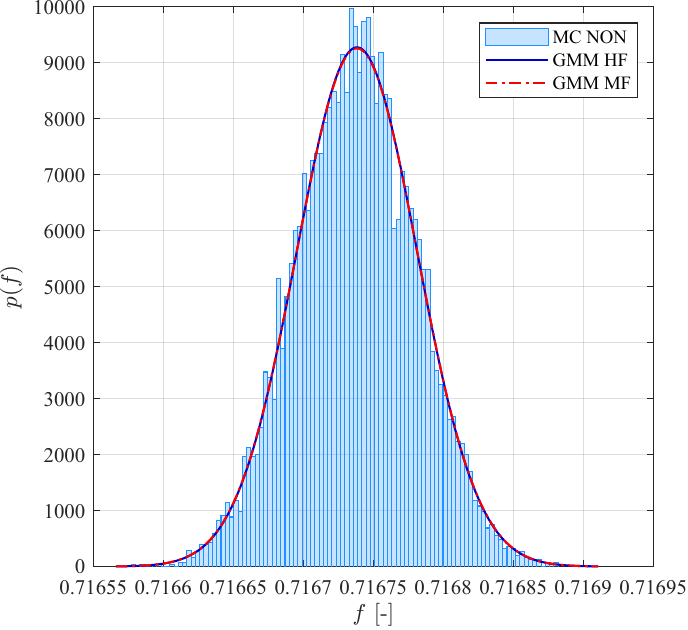}}%
    \hfill
    \subcaptionbox{$g$ parameter\label{subfig:hist_thrust_g}}[0.5\textwidth]{\includegraphics[width=0.4\textwidth]{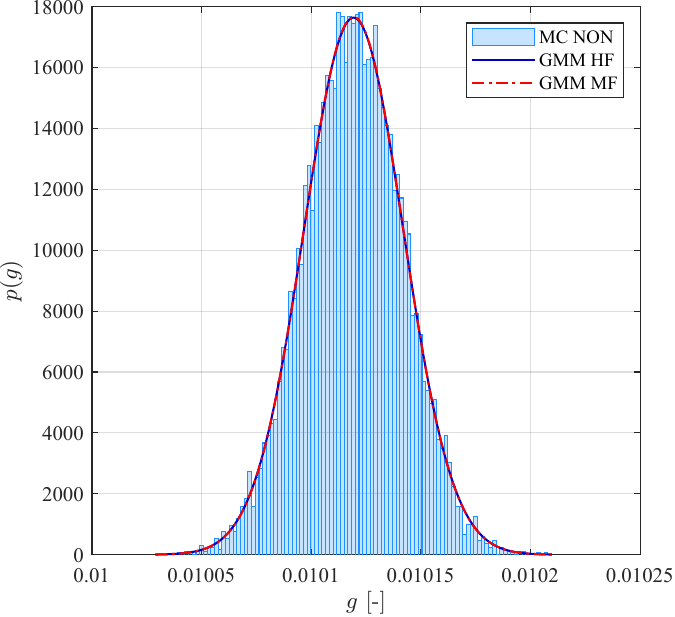}}
    \subcaptionbox{$h$ parameter\label{subfig:hist_thrust_h}}[0.5\textwidth]{\includegraphics[width=0.4\textwidth]{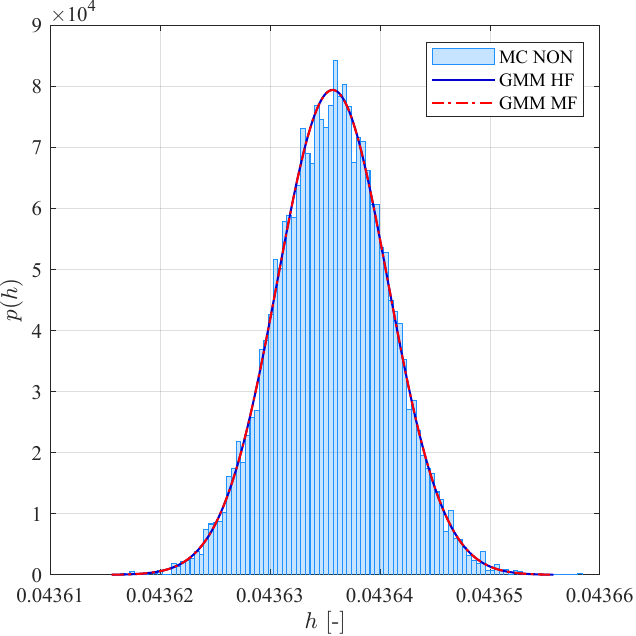}}%
    \hfill
    \subcaptionbox{$k$ parameter\label{subfig:hist_thrust_k}}[0.5\textwidth]{\includegraphics[width=0.4\textwidth]{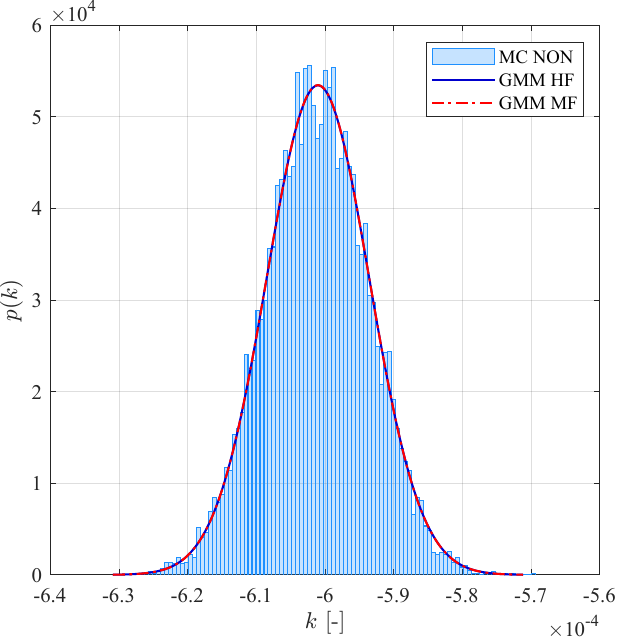}}
    \caption{Histograms and marginal \glsfmtshortpl{pdf} for low-thrust orbit raising maneuver with $J_2$ model}
    \label{fig:hist_thrust}
\end{figure}

\begin{figure}[!ht]
    \centering
    \includegraphics[width=\textwidth]{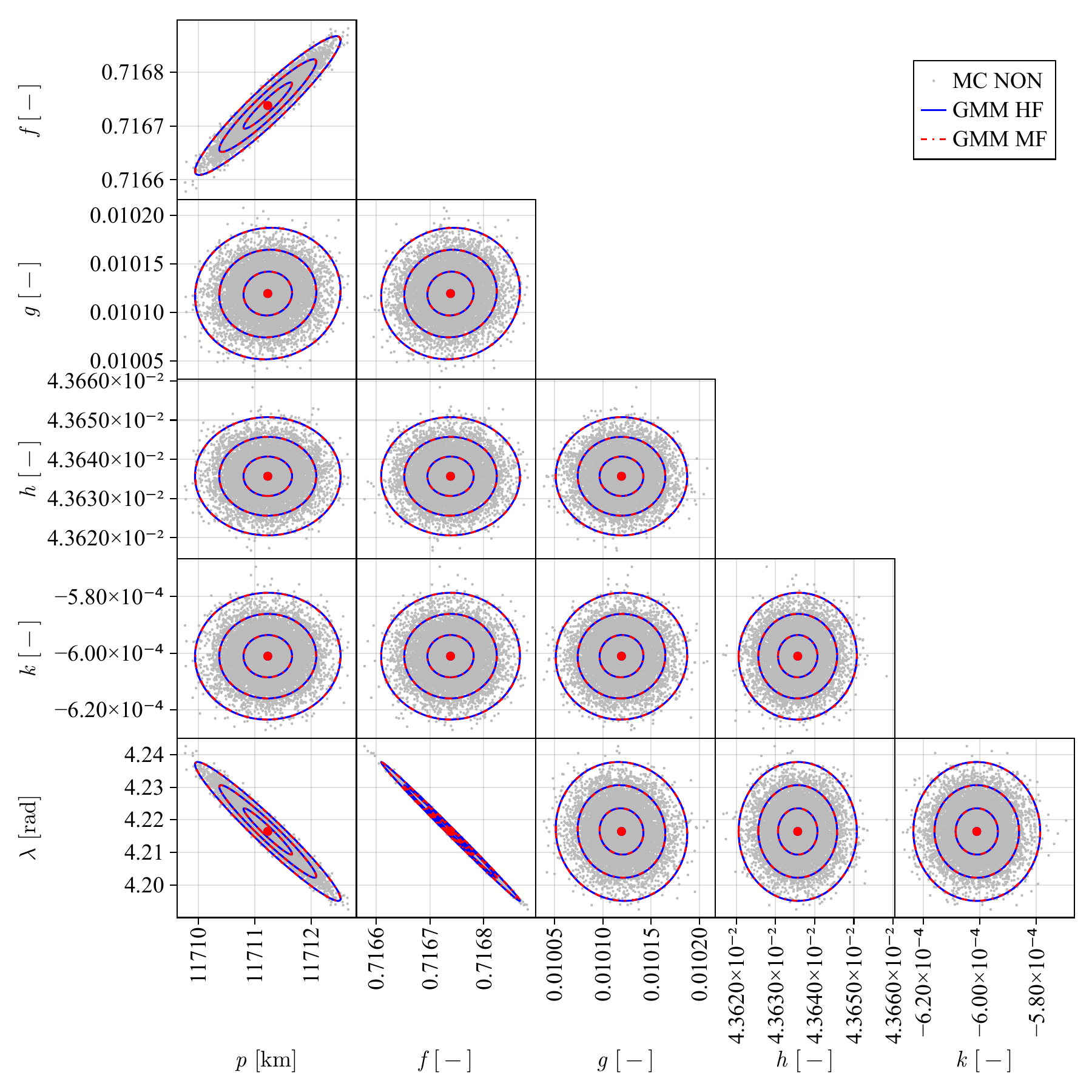}
    \caption{Projections of the \gls{gmm} and \gls{mc} covariances for low-thrust orbit raising maneuver with $J_2$ model}
    \label{fig:proj_thrust}
\end{figure}

A visual representation of the estimated statistics is then given in \cref{fig:hist_thrust,fig:proj_thrust}. \Cref{fig:hist_thrust} shows the histograms of the \gls{mc} samples (light blue bars) and the marginal \glspl{pdf} of the \gls{hf} (dark blue curves) and \gls{mf} (red curves) estimates for each \gls{mee}. In each plot, the dark blue and red curves overlap almost perfectly, thus demonstrating the capability of the \gls{mf} methods to closely reproduce the results obtained with their \gls{hf} counterparts. Moreover, these curves are in very good agreement with the \gls{mc} histograms, thus corroborating the conclusions drawn from \cref{tab:mean_tangential_thrust_mc,tab:cov_tangential_thrust_mc}. \Cref{fig:proj_thrust} shows instead a projection of the estimated covariance onto each two-dimensional subspace. The gray dots represent the \gls{mc} samples, while the dark blue and red curves are the $1\sigma,2\sigma$ and $3\sigma$ confidence ellipses of the \gls{hf} and \gls{mf} estimates, respectively. Similarly to the previous plots, the ellipses of the two solutions are almost indistinguishable, and most of the \gls{mc} samples fall within the $3\sigma$ confidence region of the two estimates. As expected, there is a strong negative correlation between the $p$ and $\lambda$ parameters in the bottom-left plot, as a smaller semilatus rectum implies a larger mean longitude. Two strong opposite correlations are also observed between the $f$ parameter and the $p,\lambda$ parameters in the top-left and bottom-second-left plots, respectively. The other projections show instead weaker correlations in agreement with \cref{tab:cov_tangential_thrust_mc}.

\begin{table}[!ht]
    \caption{Runtime for low-thrust orbit raising maneuver with $J_2$ model}
    \label{tab:runtime_tangential_thrust_mc}
    \centering
    \small
    \sisetup{round-mode=places,round-precision=3}
    \begin{tabular}{r S[scientific-notation=false,table-format=4.3] S[scientific-notation=false,table-format=1.3] S[scientific-notation=false,table-format=3.3] S[scientific-notation=false,table-format=1.3] S[scientific-notation=false,table-format=3.3]}
        \toprule
        & {$t_{\mathrm{MC}}$, \unit{\s}} & {$t_{\mathrm{MF}}$, \unit{\s}} & {$t_{\mathrm{MF}}$, \unit{\percent}} & {$t_{\mathrm{HF}}$, \unit{\s}} & {$t_{\mathrm{HF}}$, \unit{\percent}} \\ \midrule
        Adaptive \glsfmtshort{gmm} & {--} & 0.00737 & 0.503023601840097 & 0.45736 & 23.6922534992385 \\
        \Glsfmtshort{plasma} & {--} & 1.45777 & 99.4969763981599 & 1.47306 & 76.3077465007615 \\
        Total & 1364.068224 & 1.46514 & 100 & 1.93042 & 100 \\ \bottomrule
    \end{tabular}
\end{table}

\Cref{tab:runtime_tangential_thrust_mc} reports the runtime required to propagate the \gls{mc} samples and that needed by the adaptive \gls{gmm} and \gls{plasma} methods. By comparing the values on the last row, it is readily noticed that the proposed methods have a huge computational advantage over the \gls{mc} simulation. In fact, the \gls{hf} solution is about \num{710} times faster than the \gls{mc} simulation, while the \gls{mf} one is about \num{930} times faster. At the same time, the \gls{mf} solution is only about \qty{25}{\percent} faster than its \gls{hf} counterpart, which is in sharp contrast with the tenfold improvement observed in \cref{tab:runtime_tangential_thrust}. This is however expected, as the gain in performance is guaranteed by: (1) an efficient propagation of the initial uncertainty, and (2) a fast evaluation of the drift term $\vb*{u}(\vb*{x},t)$ during the computation of the \gls{pn} effects. If on one side the computational advantage of the adaptive \gls{gmm} is retained as shown in the first row of \cref{tab:runtime_tangential_thrust_mc}, the second one is instead lost due to the simplified drift model considered in this case. As a result, the \gls{bf} formulation of the \gls{sde} is less appealing, as the computational overhead introduced by the propagation of the reference trajectory $\vb*{\xi}(t)$ using the full dynamics is not compensated by the faster evaluation of the \gls{lf} drift term during the computation of the \gls{pn} effects. In other words, the \gls{bf} formulation of the \gls{plasma} algorithm is most effective when the drift term $\vb*{u}_{\mathrm{HF}}(\vb*{x},t)$ is expensive to evaluate, and a suitable \gls{lf} model $\vb*{u}_{\mathrm{LF}}(\vb*{x},t)$ is available to compensate for the computation of $\vb*{\xi}(t)$.

\afterpage{\clearpage}

\section{Conclusions}\label{sec:conclusions}

This paper presented two techniques for the efficient propagation of uncertainties in nonlinear stochastic dynamical systems. The first method, named \gls{plasma}, leverages \gls{da} to maintain a polynomial approximation of the moments of the \gls{pdf} solution to a \gls{sde} that describes the effects of unmodeled or mismodeled forces. This method was applied to the Kepler problem perturbed by a stochastic acceleration, and to an orbit raising maneuver with thrust dispersions. In both cases, the \gls{plasma} method maintained an accuracy similar to that of the reference \gls{mc} simulations, while also guaranteeing substantial savings in the overall computational cost. The second method is instead a \gls{mf} technique that aims to provide an accurate \gls{gmm} representation of the propagated \gls{pdf} while keeping the computational cost low. This objective is achieved by running the adaptive \gls{gmm} method in \citet{Losacco2024} on a \gls{lf} model of the dynamics, and then correcting its output with a \gls{pw} propagation of the initial kernels' means in \gls{hf}. Its performance was tested in three different scenarios covering the \gls{leo}, \gls{meo}, and \gls{heo} orbital regimes. In all cases, the \gls{mf} technique was able to substantially improve the estimate obtained with the \gls{lf} model, and its accuracy was similar to that of the reference \gls{mc} simulation. Moreover, its computational cost was consistently lower that that of the \gls{mc} simulation, thus making this technique suitable for resource-constrained applications. The proposed methods were finally merged to solve the problem of propagating an initial uncertainty in the presence of stochastic accelerations. The adaptive \gls{gmm} technique is still leveraged to propagate the initial uncertainty, while the \gls{plasma} method replaces the \gls{pw} propagations to provide both a \gls{hf} reference trajectory for each kernel and a quantification of the effects of \gls{pn}. This approach was applied to the aforementioned orbit raising maneuver problem, proving the effectiveness of the proposed technique in terms of both accuracy and computational efficiency.

\appendix

\section{Multi-indices}\label{apdx:indices}

A multi-index $\vc{r}$ is an $n$-vector of non-negative integers $r_i\in\mathbb{N}$ defined as
\begin{equation}
    \vc{r} = \left(r_1,\ldots,r_n\right).
    \label{eq:multi_index}
\end{equation}
Operations on multi-indices $\vc{r},\vc{r}'$ used throughout this paper are defined hereafter~\cite{LopezYela2023}
\begin{subequations}
    \allowdisplaybreaks
    \label{eq:multi_index_ops}
    \begin{align}
        \abs{\vc{r}} &= \sum_{i=1}^n r_i\\
        \vc{r}! &= \prod_{i=1}^n r_i!\\
        \dfrac{\partial^{\abs{\vc{r}}}}{\partial\vc{a}^{\vc{r}}} &= \dfrac{\partial^{\abs{\vc{r}}}}{\partial a_1^{r_1}\cdots\partial a_n^{r_n}}\\
        \vc{a}^{\vc{r}} &= \prod_{i=1}^n a_i^{r_i}\\
        \sum_{\vc{r}'=\vc{0}}^{\vc{r}} &= \sum_{r'_1=0}^{r_1} \cdots \sum_{r'_n=0}^{r_n}\\
        \binom{\vc{r}}{\vc{r}'} &= \prod_{i=1}^n \binom{r_i}{r'_i}
    \end{align}
\end{subequations}
where the convention $0! =1$ and thus $(0,\ldots,0)! = 1$ is adopted.

\section{Application to the Duffing map}\label{apdx:duffing}

The \gls{plasma} method is now applied to the Duffing map, a discrete-time dynamical system that exhibits chaotic behavior\footnote{\url{https://en.wikipedia.org/wiki/Duffing\_map}}. The main objective is to illustrate the algorithm discussed in \cref{sec:plasma:da}. The Duffing map transforms points $(x_{k-1},y_{k-1})\in\mathbb{R}^2$ into new points given by
\begin{equation}
    \begin{dcases}
        x_{k} = y_{k-1}\\
        y_{k} = -bx_{k-1}+ay_{k-1}-y_{k-1}^3
    \end{dcases}
    \label{eq:duffing_map}
\end{equation}
with $a,b\in\mathbb{R}$ constants. To analyze the properties of the \gls{plasma} method, $y_k$ is perturbed by a random sequence $w_k$ such that
\begin{equation}
    \begin{aligned}
        \mathbb{E}[w_k] &= 0\\
        \mathbb{E}[w_i w_j] &= \delta_{ij}
    \end{aligned}
\end{equation}
with $\delta_{ij}$ Kronecker delta. \Cref{eq:duffing_map} is thus rewritten as
\begin{equation}
    \begin{dcases}
        \hat{x}_k = \hat{y}_{k-1}\\
        \hat{y}_k = -b\hat{x}_{k-1}+a\hat{y}_{k-1}-\hat{y}_{k-1}^3 + \sigma w_k
    \end{dcases}
    \label{eq:duffing_map_noisy}
\end{equation}
with $\sigma\in\mathbb{R}$ noise standard deviation.

Assuming an expansion order $N=2$, explicit expressions for \cref{eq:effective_noise_expansion,eq:effective_noise_moments,eq:poly_moments} are then obtained as follows. Firstly, \cref{eq:duffing_map_noisy} is reformulated to match the form of \cref{eq:em_scheme} as
\begin{equation}
    \begin{dcases}
        \hat{x}_k = \hat{x}_{k-1} + \left[\hat{y}_{k-1}-\hat{x}_{k-1}\right]\\
        \hat{y}_k = \hat{y}_{k-1} + \left[-b\hat{x}_{k-1}+\left(a-1\right)\hat{y}_{k-1}-\hat{y}_{k-1}^3\right] + \sigma w_k
    \end{dcases}
    \label{eq:duffing_sde}
\end{equation}
where the time step has been set to $h=1$. The state vector $\hat{\vc{X}}_{k-1}$, drift coefficient $\vc{u}(\hat{\vc{X}}_{k-1})$ and diffusion coefficient $\vc{G}(\hat{\vc{X}}_{k-1})$ are thus identified as
\begin{subequations}
    \allowdisplaybreaks
    \begin{align}
        \hat{\vc{X}}_{k-1} &= \begin{bmatrix}\hat{x}_{k-1}\\ \hat{y}_{k-1}\end{bmatrix}\\
        \vc{u}(\hat{\vc{X}}_{k-1}) &= \begin{bmatrix}\hat{y}_{k-1}-\hat{x}_{k-1}\\ -b\hat{x}_{k-1}+\left(a-1\right)\hat{y}_{k-1}-\hat{y}_{k-1}^3\end{bmatrix}\\
        \vc{G}(\hat{\vc{X}}_{k-1}) &= \begin{bmatrix}0\\ \sigma\end{bmatrix}
    \end{align}
\end{subequations}
where the dependency on time has been dropped for clarity.

The polynomial expansion of the two components of the effective noise are then computed from \cref{eq:effective_noise_expansion} as
\begin{subequations}
    \label{eq:poly_exp_noise_duffing}
    \begin{align}
        \begin{split}
            \Delta\hat{W}^{(1)}_{k,2} &= \Delta\hat{W}^{(1)}_{k-1,2} + \left[\Delta\hat{W}^{(2)}_{k-1,2}-\Delta\hat{W}^{(1)}_{k-1,2}\right]\\
            &= \Delta\hat{W}^{(2)}_{k-1,2}
        \end{split}\label{eq:poly_exp_noise_duffing_1}\\
        \begin{split}
            \Delta\hat{W}^{(2)}_{k,2} &= \Delta\hat{W}^{(2)}_{k-1,2} + \left[-b\Delta\hat{W}^{(1)}_{k-1,2}+\left(a-1-3\hat{y}_{k-1}^2\right)\Delta\hat{W}^{(2)}_{k-1,2}-3\hat{y}_{k-1}\left(\Delta\hat{W}^{(2)}_{k-1,2}\right)^2\right]+\sigma\Delta W_{k}\\
            &=-b\Delta\hat{W}^{(1)}_{k-1,2}+\left(a-3\hat{y}_{k-1}^2\right)\Delta\hat{W}^{(2)}_{k-1,2}-3\hat{y}_{k-1}\left(\Delta\hat{W}^{(2)}_{k-1,2}\right)^2+\sigma\Delta W_{k}
        \end{split}\label{eq:poly_exp_noise_duffing_2}
    \end{align}
\end{subequations}
where the dependency on $\vc{x}_0$ has also been dropped for clarity. \Cref{eq:poly_exp_noise_duffing_1,eq:poly_exp_noise_duffing_2} are second-order polynomials in $\left\{\Delta\hat{W}^{(1)}_{k-1,2},\Delta\hat{W}^{(2)}_{k-1,2},\Delta W_{k}\right\}$, i.e.
\begin{subequations}
    \begin{align}
        \Delta\hat{W}^{(1)}_{k,2} &= p^{(1)}_{k,2}\left(\Delta\hat{W}^{(1)}_{k-1,2},\Delta\hat{W}^{(2)}_{k-1,2},\Delta W_{k}\right)\\
        \Delta\hat{W}^{(2)}_{k,2} &= p^{(2)}_{k,2}\left(\Delta\hat{W}^{(1)}_{k-1,2},\Delta\hat{W}^{(2)}_{k-1,2},\Delta W_{k}\right)
    \end{align}
\end{subequations}

The conditional moments of the effective noise up to second order are then obtained by computing the polynomial expansions $\Delta\hat{\vc{W}}_{k,2}^{\vc{r}}=p^{\vc{r}}_{k,2}\left(\Delta\hat{W}^{(1)}_{k-1,2},\Delta\hat{W}^{(2)}_{k-1,2},\Delta W_{k}\right)$ appearing in \cref{eq:effective_noise_moments} $\forall\vc{r}=(r_1,r_2):\abs{\vc{r}}=1,2$. These expansions are obtained by taking the products of $\Delta\hat{W}^{(1)}_{k,2}$ and $\Delta\hat{W}^{(2)}_{k,2}$ in \cref{eq:poly_exp_noise_duffing} truncated at second order, thus resulting in the following explicit expressions for $\Delta\hat{\vc{W}}_{k,2}^{\vc{r}}$
\begin{subequations}
    \label{eq:prods_poly_exp_noise_duffing}
    \allowdisplaybreaks
    \begin{align}
        \begin{split}
            \Delta\hat{\vc{W}}_{k,2}^{(1,0)} &= \Delta\hat{W}^{(1)}_{k,2}\\
            &= \Delta\hat{W}^{(2)}_{k-1,2}
        \end{split}\\
        \begin{split}
            \Delta\hat{\vc{W}}_{k,2}^{(0,1)} &= \Delta\hat{W}^{(2)}_{k,2}\\
            &= -b\Delta\hat{W}^{(1)}_{k-1,2}+\left(a-3\hat{y}_{k-1}^2\right)\Delta\hat{W}^{(2)}_{k-1,2}-3\hat{y}_{k-1}\left(\Delta\hat{W}^{(2)}_{k-1,2}\right)^2+\sigma\Delta W_{k}
        \end{split}\\
        \begin{split}
            \Delta\hat{\vc{W}}_{k,2}^{(2,0)} &= \Delta\hat{W}^{(1)}_{k,2}\Delta\hat{W}^{(1)}_{k,2}\\
            &= \left(\Delta\hat{W}^{(2)}_{k-1,2}\right)^2
        \end{split}\\
        \begin{split}
            \Delta\hat{\vc{W}}_{k,2}^{(1,1)} &= \Delta\hat{W}^{(1)}_{k,2}\Delta\hat{W}^{(2)}_{k,2}\\
            &= -b\Delta\hat{W}^{(1)}_{k-1,2}\Delta\hat{W}^{(2)}_{k-1,2}+\left(a-3\hat{y}_{k-1}^2\right)\left(\Delta\hat{W}^{(2)}_{k-1,2}\right)^2+\sigma\Delta\hat{W}^{(2)}_{k-1,2}\Delta W_{k}
        \end{split}\\
        \begin{split}
            \Delta\hat{\vc{W}}_{k,2}^{(0,2)} &= \Delta\hat{W}^{(2)}_{k,2}\Delta\hat{W}^{(2)}_{k,2}\\
            &=b^2\left(\Delta\hat{W}^{(1)}_{k-1,2}\right)^2
            -2b\left(a-3\hat{y}_{k-1}^2\right)\Delta\hat{W}^{(1)}_{k-1,2}\Delta\hat{W}^{(2)}_{k-1,2}
            +\left(a-3\hat{y}_{k-1}^2\right)^2\left(\Delta\hat{W}^{(2)}_{k-1,2}\right)^2\\
            &-2b\sigma\Delta\hat{W}^{(1)}_{k-1,2}\Delta W_{k}
            +2\sigma\left(a-3\hat{y}_{k-1}^2\right)\Delta\hat{W}^{(2)}_{k-1,2}\Delta W_{k}
            +\sigma^2\Delta W_{k}^2
        \end{split}
    \end{align}
\end{subequations}
where these relations hold from the multi-index notation introduced in \cref{apdx:indices}.

Applying the expectation operator to the \gls{rhs} of \cref{eq:prods_poly_exp_noise_duffing} then leads to the following expressions for the conditional moments of the effective noise $\ex{\Delta\hat{\vc{W}}_{k,2}^{\vc{r}}}$
\begin{subequations}
    \label{eq:noise_exp_duffing}
    \allowdisplaybreaks
    \begin{align}
        \begin{split}
            \ex{\Delta\hat{\vc{W}}^{(1,0)}_{k,2}} &= \ex{\Delta\hat{W}^{(2)}_{k-1,2}}\\
            &= \ex{\Delta\hat{\vc{W}}^{(0,1)}_{k-1,2}}
        \end{split}\\
        \begin{split}
            \ex{\Delta\hat{\vc{W}}^{(0,1)}_{k,2}}
            &=\ex{-b\Delta\hat{W}^{(1)}_{k-1,2}+\left(a-3\hat{y}_{k-1}^2\right)\Delta\hat{W}^{(2)}_{k-1,2}-3\hat{y}_{k-1}\left(\Delta\hat{W}^{(2)}_{k-1,2}\right)^2+\sigma\Delta W_{k}}\\
            &=-b\ex{\Delta\hat{\vc{W}}^{(1,0)}_{k-1,2}}
            +\left(a-3\hat{y}_{k-1}^2\right)\ex{\Delta\hat{\vc{W}}^{(0,1)}_{k-1,2}}
            -3\hat{y}_{k-1}\ex{\Delta\hat{\vc{W}}^{(0,2)}_{k-1,2}}
        \end{split}\\
        \begin{split}
            \ex{\Delta\hat{\vc{W}}^{(2,0)}_{k,2}} &= \ex{\left(\Delta\hat{W}^{(2)}_{k-1,2}\right)^2}\\
            &=\ex{\Delta\hat{\vc{W}}^{(0,2)}_{k-1,2}}
        \end{split}\\
        \begin{split}
            \ex{\Delta\hat{\vc{W}}^{(1,1)}_{k,2}} &= \ex{-b\Delta\hat{W}^{(1)}_{k-1,2}\Delta\hat{W}^{(2)}_{k-1,2}+\left(a-3\hat{y}_{k-1}^2\right)\left(\Delta\hat{W}^{(2)}_{k-1,2}\right)^2+\sigma\Delta\hat{W}^{(2)}_{k-1,2}\Delta W_{k}}\\
            &=-b\ex{\Delta\hat{\vc{W}}^{(1,1)}_{k-1,2}}
            +\left(a-3\hat{y}_{k-1}^2\right)\ex{\Delta\hat{\vc{W}}^{(0,2)}_{k-1,2}}
        \end{split}\\
        \begin{split}
            \ex{\Delta\hat{\vc{W}}^{(0,2)}_{k,2}} &=\mathbb{E}\biggl[b^2\left(\Delta\hat{W}^{(1)}_{k-1,2}\right)^2-2b\left(a-3\hat{y}_{k-1}^2\right)\Delta\hat{W}^{(1)}_{k-1,2}\Delta\hat{W}^{(2)}_{k-1,2}+\left(a-3\hat{y}_{k-1}^2\right)^2\left(\Delta\hat{W}^{(2)}_{k-1,2}\right)^2 \\
            &\qquad -2b\sigma\Delta\hat{W}^{(1)}_{k-1,2}\Delta W_{k}+2\sigma\left(a-3\hat{y}_{k-1}^2\right)\Delta\hat{W}^{(2)}_{k-1,2}\Delta W_{k}+\sigma^2\Delta W_{k}^2\biggr] \\
            &=b^2\ex{\Delta\hat{\vc{W}}^{(2,0)}_{k-1,2}}
            -2b\left(a-3\hat{y}_{k-1}^2\right)\ex{\Delta\hat{\vc{W}}^{(1,1)}_{k-1,2}}+\left(a-3\hat{y}_{k-1}^2\right)^2\ex{\Delta\hat{\vc{W}}^{(0,2)}_{k-1,2}}+\sigma^2
        \end{split}
    \end{align}
\end{subequations}
where it has been made use of the property of linearity of the expectation and of the assumption that $\Delta W_k$ is a random sequence with independent increments such that
\begin{equation}
    \begin{aligned}
        \ex{\Delta W_k} = 0\\
        \ex{\Delta W_k^2} = 1
    \end{aligned}
\end{equation}
and $\Delta\hat{\vc{W}}^{\vc{r}}_{k-1,2}$ is independent from $\Delta W_k\,\forall\,k\in\mathbb{N}^{+}$. The conditional moments of the state $\hat{\vc{X}}_k$ are then readily obtained from \cref{eq:poly_moments} as
\begin{subequations}
    \label{eq:state_exp_duffing}
    \allowdisplaybreaks
    \begin{align}
        \begin{split}
            \ex{\hat{\vc{X}}^{(1,0)}_{k,2}}
            &=\hat{\vc{x}}_k^{(1,0)}\ex{\Delta\hat{\vc{W}}^{(0,0)}_{k,2}}+\hat{\vc{x}}_k^{(0,0)}\ex{\Delta\hat{\vc{W}}^{(1,0)}_{k,2}}\\
            &=\hat{x}_k+\ex{\Delta\hat{\vc{W}}^{(1,0)}_{k,2}}
        \end{split}\\
        \begin{split}
            \ex{\hat{\vc{X}}^{(0,1)}_{k,2}}
            &=\hat{\vc{x}}_k^{(0,1)}\ex{\Delta\hat{\vc{W}}^{(0,0)}_{k,2}}+\hat{\vc{x}}_k^{(0,0)}\ex{\Delta\hat{\vc{W}}^{(0,1)}_{k,2}}\\
            &=\hat{y}_k+\ex{\Delta\hat{\vc{W}}^{(0,1)}_{k,2}}
        \end{split}\\
        \begin{split}
            \ex{\hat{\vc{X}}^{(2,0)}_{k,2}}
            &=\hat{\vc{x}}_k^{(2,0)}\ex{\Delta\hat{\vc{W}}^{(0,0)}_{k,2}}+2\hat{\vc{x}}_k^{(1,0)}\ex{\Delta\hat{\vc{W}}^{(1,0)}_{k,2}}+\hat{\vc{x}}_k^{(0,0)}\ex{\Delta\hat{\vc{W}}^{(2,0)}_{k,2}}\\
            &=\hat{x}_k^2+2\hat{x}_k\ex{\Delta\hat{\vc{W}}^{(1,0)}_{k,2}}+\ex{\Delta\hat{\vc{W}}^{(2,0)}_{k,2}}
        \end{split}\\
        \begin{split}
            \ex{\hat{\vc{X}}^{(1,1)}_{k,2}}
            &=\hat{\vc{x}}_k^{(1,1)}\ex{\Delta\hat{\vc{W}}^{(0,0)}_{k,2}}+\hat{\vc{x}}_k^{(0,1)}\ex{\Delta\hat{\vc{W}}^{(1,0)}_{k,2}}+\hat{\vc{x}}_k^{(1,0)}\ex{\Delta\hat{\vc{W}}^{(0,1)}_{k,2}}+\hat{\vc{x}}_k^{(0,0)}\ex{\Delta\hat{\vc{W}}^{(1,1)}_{k,2}}\\
            &=\hat{x}_k\hat{y}_k+\hat{y}_k\ex{\Delta\hat{\vc{W}}^{(1,0)}_{k,2}}+\hat{x}_k\ex{\Delta\hat{\vc{W}}^{(0,1)}_{k,2}}+\ex{\Delta\hat{\vc{W}}^{(1,1)}_{k,2}}
        \end{split}\\
        \begin{split}
            \ex{\hat{\vc{X}}^{(0,2)}_{k,2}}
            &=\hat{\vc{x}}_k^{(0,2)}\ex{\Delta\hat{\vc{W}}^{(0,0)}_{k,2}}+2\hat{\vc{x}}_k^{(0,1)}\ex{\Delta\hat{\vc{W}}^{(0,1)}_{k,2}}+\hat{\vc{x}}_k^{(0,0)}\ex{\Delta\hat{\vc{W}}^{(0,2)}_{k,2}}\\
            &=\hat{y}_k^2+2\hat{y}_k\ex{\Delta\hat{\vc{W}}^{(0,1)}_{k,2}}+\ex{\Delta\hat{\vc{W}}^{(0,2)}_{k,2}}
        \end{split}
    \end{align}
\end{subequations}
where $\ex{\Delta\hat{\vc{W}}^{(0,0)}_{k,2}}=1$ and $\hat{\vc{x}}_k=[\hat{x}_k\ \hat{y}_k]^T$ is the central part of $\hat{\vc{X}}_{k,2}$ which coincides with the nominal solution to \cref{eq:duffing_map}. The entries of the covariance matrix $\hat{\vc{P}}_{k,2}$ are finally computed using \cref{eq:poly_cov} as
\begin{equation}
    \label{eq:cov_duffing}
    \begin{aligned}
        \hat{P}_{k,2}^{(1,0)} &= \ex{\hat{\vc{X}}^{(2,0)}_{k,2}}-\ex{\hat{\vc{X}}^{(1,0)}_{k,2}}\cdot\ex{\hat{\vc{X}}^{(1,0)}_{k,2}}\\
        \hat{P}_{k,2}^{(1,0)}=\hat{P}_{k,2}^{(0,1)} &= \ex{\hat{\vc{X}}^{(1,1)}_{k,2}}-\ex{\hat{\vc{X}}^{(1,0)}_{k,2}}\cdot\ex{\hat{\vc{X}}^{(0,1)}_{k,2}}\\
        \hat{P}_{k,2}^{(1,1)} &= \ex{\hat{\vc{X}}^{(0,2)}_{k,2}}-\ex{\hat{\vc{X}}^{(0,1)}_{k,2}}\cdot\ex{\hat{\vc{X}}^{(0,1)}_{k,2}}.
    \end{aligned}
\end{equation}

These expression were demonstrated to provide results that agree within machine precision with those obtained with the \gls{plasma} method described above.

\section*{Acknowledgments}

This work is co-funded by the \gls{cnes} through A. Foss\`a's PhD program, and made use of the \gls{cnes} orbital propagation tools, including the \gls{pace} library.

\bibliography{references}

\end{document}

%% file: acronyms.tex
% A
\newabbreviation[type=\acronymtype,category=ls]{ads}{ADS}{automatic domain splitting}
\newabbreviation[type=\acronymtype,category=ls]{aee}{AEqOE}{alternate equinoctial element}
\newabbreviation[type=\acronymtype,category=ls]{aegis}{AEGIS}{Adaptive Entropy-based Gaussian-mixture Information Synthesis}

% B
\newabbreviation[type=\acronymtype,category=l]{bf}{BF}{bi-fidelity}
\newabbreviation[type=\acronymtype,category=ls]{bvp}{BVP}{boundary value problem}

% C
\newabbreviation[type=\acronymtype,category=ls]{cam}{CAM}{collision avoidance maneuver}
\newabbreviation[type=\acronymtype,category=s]{cnes}{CNES}{Centre National d'\'Etudes Spatiales}
\newabbreviation[type=\acronymtype,category=ls]{coe}{COE}{classical orbital element}
\newabbreviation[type=\acronymtype,category=ls]{cut}{CUT}{conjugate unscented transform}

% D
\newabbreviation[type=\acronymtype,category=ls]{da}{DA}{differential algebra}
\newabbreviation[type=\acronymtype,category=s]{de}{DE}{Developmental Ephemerides}
\newabbreviation[type=\acronymtype,category=ls]{dmc}{DMC}{dynamic model compensation}
\newabbreviation[type=\acronymtype,category=ls]{dmm}{DMM}{Dirac mixture model}
\newabbreviation[type=\acronymtype,category=ls]{dop853}{DOP853}{Dorman-Prince 8(5,3)}
\newabbreviation[type=\acronymtype,category=ls]{dsst}{DSST}{Draper Semi-analytical Satellite Theory}

% E
\newabbreviation[type=\acronymtype,category=ls]{ee}{EqOE}{equinoctial element}
\newabbreviation[type=\acronymtype,category=ls]{ekf}{EKF}{extended Kalman filter}
\newabbreviation[type=\acronymtype,category=ls]{em}{EM}{Euler-Maruyama}
\newabbreviation[type=\acronymtype,category=l]{exmx}{EM}{expectation maximization}

% F
\newabbreviation[type=\acronymtype,category=ls]{fp}{FP}{floating point}
\newabbreviation[type=\acronymtype,category=ls]{fpe}{FPE}{Fokker-Planck equation}
\newabbreviation[type=\acronymtype,category=ls]{fpke}{FPKE}{Fokker-Planck-Kolmogorov equation}

% G
\newabbreviation[type=\acronymtype,category=ls]{gcrf}{GCRF}{geocentric celestial reference frame}
\newabbreviation[type=\acronymtype,category=ls]{geo}{GEO}{geostationary equatorial orbit}
\newabbreviation[type=\acronymtype,category=ls]{gmm}{GMM}{Gaussian mixture model}
\newabbreviation[type=\acronymtype,category=ls]{gto}{GTO}{geostationary transfer orbit}

% H
\newabbreviation[type=\acronymtype,category=ls]{heo}{HEO}{high Earth orbit}
\newabbreviation[type=\acronymtype,category=l]{hf}{HF}{high-fidelity}

% I
\newabbreviation[type=\acronymtype,category=ls]{ic}{IC}{initial condition}
\newabbreviation[type=\acronymtype,category=ls]{itrf}{ITRF}{international terrestrial reference frame}
\newabbreviation[type=\acronymtype,category=ls]{ivp}{IVP}{initial value problem}

% J
\newabbreviation[type=\acronymtype,category=s]{jpl}{JPL}{Jet Propulsion Laboratory}

% L
\newabbreviation[type=\acronymtype,category=ls]{lam}{LAM}{likelihood agreement measure}
\newabbreviation[type=\acronymtype,category=ls]{lc}{LinCov}{linear covariance}
\newabbreviation[type=\acronymtype,category=ls]{leo}{LEO}{low Earth orbit}
\newabbreviation[type=\acronymtype,category=l]{lf}{LF}{low-fidelity}
\newabbreviation[type=\acronymtype,category=ls]{loads}{LOADS}{low-order automatic domain splitting}

% M
\newabbreviation[type=\acronymtype,category=ls]{mc}{MC}{Monte Carlo}
\newabbreviation[type=\acronymtype,category=ls]{mce}{MCE}{moment constraint equation}
\newabbreviation[type=\acronymtype,category=ls]{mee}{MEE}{modified equinoctial element}
\newabbreviation[type=\acronymtype,category=ls]{meo}{MEO}{medium Earth orbit}
\newabbreviation[type=\acronymtype,category=l]{mf}{MF}{multifidelity}

% N
\newabbreviation[type=\acronymtype,category=ls]{nli}{NLI}{nonlinearity index}

% O
\newabbreviation[type=\acronymtype,category=ls]{od}{OD}{orbit determination}
\newabbreviation[type=\acronymtype,category=ls]{ode}{ODE}{ordinary differential equation}
\newabbreviation[type=\acronymtype,category=ls]{oe}{OE}{orbital element}

% P
\newabbreviation[type=\acronymtype,category=sl]{pace}{PACE}{Polynomial Algebra Computational Engine}
\newabbreviation[type=\acronymtype,category=ls]{pce}{PCE}{polynomial chaos expansion}
\newabbreviation[type=\acronymtype,category=ls]{pde}{PDE}{partial differential equation}
\newabbreviation[type=\acronymtype,category=ls]{pdf}{pdf}{probability density function}
\newabbreviation[type=\acronymtype,category=ls]{plasma}{PLASMA}{PoLynomial Algebra Stochastic Moments Analysis}
\newabbreviation[type=\acronymtype,category=l]{pn}{PN}{process noise}
\newabbreviation[type=\acronymtype,category=ls]{poc}{PoC}{probability of collision}
\newabbreviation[type=\acronymtype,category=ls]{psd}{PSD}{power spectral density}
\newabbreviation[type=\acronymtype,category=l]{pw}{PW}{point-wise}

% R
\newabbreviation[type=\acronymtype,category=s]{ram}{RAM}{random access memory}
\newabbreviation[type=\acronymtype,category=ls]{rhs}{RHS}{right-hand side}
\newabbreviation[type=\acronymtype,category=ls]{rk}{RK}{Runge-Kutta}
\newabbreviation[type=\acronymtype,category=ls]{rmse}{RMSE}{root mean square error}

% S
\newabbreviation[type=\acronymtype,category=ls]{sc}{SC}{stochastic collocation}
\newabbreviation[type=\acronymtype,category=ls]{sde}{SDE}{stochastic differential equation}
\newabbreviation[type=\acronymtype,category=ls]{sgp4}{SGP4}{Simplified General Perturbations 4}
\newabbreviation[type=\acronymtype,category=ls]{snc}{SNC}{state noise compensation}
\newabbreviation[type=\acronymtype,category=ls]{so}{SO}{space object}
\newabbreviation[type=\acronymtype,category=ls]{sode}{SODE}{stochastic ordinary differential equation}
\newabbreviation[type=\acronymtype,category=ls]{srp}{SRP}{solar radiation pressure}
\newabbreviation[type=\acronymtype,category=ls]{ssa}{SSA}{space situational awareness}
\newabbreviation[type=\acronymtype,category=ls]{sst}{SST}{space surveillance and tracking}
\newabbreviation[type=\acronymtype,category=ls]{stela}{STELA}{Semi-analytic Tool for End of Life Analysis}
\newabbreviation[type=\acronymtype,category=ls]{stm}{STM}{state transition matrix}
\newabbreviation[type=\acronymtype,category=ls]{stt}{STT}{state transition tensor}

% T
\newabbreviation[type=\acronymtype,category=ls]{tca}{TCA}{time of closest approach}
\newabbreviation[type=\acronymtype,category=ls]{teme}{TEME}{true equator mean equinox}
\newabbreviation[type=\acronymtype,category=ls]{tle}{TLE}{two-line element}
\newabbreviation[type=\acronymtype,category=ls]{tnw}{TNW}{tangential, normal, $\omega$}

% U
\newabbreviation[type=\acronymtype,category=ls]{ukf}{UKF}{unscented Kalman filter}
\newabbreviation[type=\acronymtype,category=ls]{up}{UP}{uncertainty propagation}
\newabbreviation[type=\acronymtype,category=ls]{ut}{UT}{unscented transform}
\newabbreviation[type=\acronymtype,category=s]{utc}{UTC}{coordinated universal time}